\documentclass[11pt]{article}
\usepackage{graphicx}
\usepackage[margin=1truein]{geometry}
\usepackage[dvipsnames,svgnames]{xcolor}
\usepackage{setspace}
\usepackage{amsmath,amsfonts,amssymb,amsthm,mathtools}  
\usepackage{enumerate}  
\usepackage[noprefix]{nomencl}
\makenomenclature
\usepackage[boxruled,vlined]{algorithm2e}
\usepackage{autobreak}
\usepackage{float}
\usepackage{indentfirst}
\SetKwInOut{Input}{Input}
\SetKwInOut{Output}{Output}
\usepackage{etoolbox}
\usepackage{url}
\usepackage[colorlinks,linkcolor=red,citecolor={blue}]{hyperref}
\usepackage[square,numbers]{natbib}
\def\X{\mathbb{X}}
\def\A{\mathbb{A}}
\def\R{\mathbb{R}}
\def\N{\mathbb{N}}

\def\P{\mathbb{P}}

\def\F{\mathbb{F}}
\def\I{\mathbb{I}}
\newtheorem*{abc0}{Condition $\mathrm{C}$}
\newtheorem*{abc1}{Condition $\mathrm{A_-}$}
\newtheorem*{abc2}{Condition $\mathrm{A_+}$}  
\newtheorem*{con3}{Condition $\mathrm{A}$}
\newtheorem*{abc3}{Condition $\mathrm{B_-}$} 
\newtheorem*{abc4}{Condition $\mathrm{B_+}$}

\newtheorem{thm}[algocf]{Theorem}
\newtheorem{pro}[algocf]{Proposition}
\newtheorem{lem}[algocf]{Lemma}
\newtheorem{cor}[algocf]{Corollary}
\newtheorem{rem}[algocf]{Remark}
\newtheorem{exa}[algocf]{Example}
\newtheorem{defi}[algocf]{Definition}
\DeclarePairedDelimiter\abs{\lvert}{\rvert}%
\makeatletter
\let\c@algocf\c@subsection
\makeatother

\numberwithin{equation}{section}
\allowdisplaybreaks[4]

\newcommand{\norm}[1]{\left\lVert#1\right\rVert}
\newcommand{\argmax}{\operatornamewithlimits{argmax}}

\begin{document}
	\title{Properties of Turnpike Functions for Discounted  Finite Markov Decision Processes
    \footnote{
       This research was partially supported  by the U.S. Office of Naval Research (ONR) under Grant N000142412608. An extended abstract presenting some of the results
of this submission were published without proofs in conference proceedings \cite{PC}. This extended abstract does not contain proofs. }}
    
	\author{
    Eugene A. Feinberg \\
    \small{Department of Applied Mathematics and Statistics}\\
\small{Stony Brook University}\\
    \small{Stony Brook, 11733}
    \and
    Gaojin He\\
    \small{Department of Mathematics}\\
\small{University of California San Diego}\\
\small{La Jolla, 92093}
}
\author{Eugene A. Feinberg\thanks{Department of Applied Mathematics and Statistics,
        Stony Brook University,
        Stony Brook, 11733}       
    \and    
    Gaojin He\thanks{Department of Mathematics,
        University of California San Diego,
        La Jolla, 92093}}
	\date{}
	\maketitle
	\begin{abstract}
		This paper studies convergence times of the Value Iteration Algorithm (VIA) for discounted discrete-time Markov Decision Processes (MDPs) with finite state and action sets. For each discount factor, starting from a finite number of iterations, which is called the turnpike integer, the VIA generates deterministic optimal policies for infinite-horizon problems. Turnpike integers are viewed as functions of discount factors called turnpike functions. We design an algorithm to detect in strongly polynomial time whether a discount factor is irregular -- at which the sets of optimal policies change. We then prove that a turnpike function is upper-semicontinuous and finitely-piecewise constant on any closed subinterval of $[0,1)$ that does not contain irregular points. In particular, we prove that for small discount factors a turnpike function is bounded by the number of states and design an algorithm based on the VIA to solve an MDP for all small discount factors in strongly polynomial time.\\
		
		\noindent\begin{keywords}
			Markov Decision Process, Value Iteration, Policy, Turnpike Theorems, Rolling Horizon
		\end{keywords}
	\end{abstract}

\section{Introduction}\label{sec:Intro}
This paper deals with discounted Markov Decision Processes (MDPs) with finite state and action sets. Policy and value iterations are two classic methods for solving MDPs. For a deterministic policy $\phi,$ sometimes also called stationary, deterministic stationary, or nonrandomized stationary,  the main step of the policy iteration algorithm (PIA) either detects that $\phi$ is optimal or improves it. Starting from a deterministic policy $\phi^{(0)},$ the PIA generates a finite sequence of deterministic policies $\phi^{(0)},\phi^{(1)},\ldots,\phi^{(n)},$ such that $v_\alpha^{\phi^{(i+1)}}(x) \ge v_\alpha^{\phi^{(i)}}(x)$ for all $x\in\X$  and $v_\alpha^{\phi^{(i+1)}}(x) > v_\alpha^{\phi^{(i)}}(x)$  for some $x\in\X,$  where $\phi^{(n)}$ is an optimal policy, $\X$ is the state space, $i=0,\ldots,n-1,$ and $v_\alpha^\phi(x)$ is the infinite-horizon expected total discounted reward earned by the policy $\phi,$ if $x$ is the initial state, and $\alpha\in [0,1)$ is the discount factor. After PI reaches an optimal policy $\phi^{(n)},$ in principle, it can be continued, but it will generate optimal policies. PIA is strongly polynomial with the bound depending on the discount factor \cite{PIS}. This means that the number of iterations to find an optimal policy is bounded as \begin{equation}\label{eq*} \eta(\alpha,m,q)\le f(\alpha)g(m,q),\qquad \alpha\in [0,1),\end{equation} where $m$ is the number of states, $q$ is the number of state-action pairs, $f$ and $g$ are nonnegative continuous real-valued functions, $f$ is increasing, and $g$ is polynomial. In particular, for a fixed problem, \begin{equation}\label{eq**} \eta(\alpha)\le F(\alpha),\qquad \alpha\in [0,1),{\rm \ and \ } m  {\rm \ and\ }     q {\rm\ are\ fixed,}\end{equation} where $\eta(\alpha):= \eta(\alpha,m,q),$ and $F(\alpha):=f(\alpha)g(m,q)$ is a continuous increasing function. MDPs can also be solved by linear programming (LP). PIA is an implementation of the simplex method with the block pivoting rule. Different pivoting rules define different versions of PIAs. For the version defined by Dantzig's pivoting rule, the above estimations also hold \cite{PIS}; see also~\cite{Scherrer}. For deterministic MDPs, PIA with Dantzig's pivoting rule is strongly polynomial for all discount factors \cite{SPD}, that is, $f$ and $F$ can be substituted with constants in (\ref{eq*}) and (\ref{eq**}). 

The value iteration algorithm (VIA) sequentially computes finite-horizon value functions $V_{n,\alpha},$ which converge to the infinite-horizon value function, and this convergence is geometrically fast. VIA also can be used to construct a sequence of deterministic policies $\phi^{(1)},\phi^{(2)},...,\phi^{(n)},$ where $\phi^{(i)}$ is some policy from the output of the $i$-th iteration of the VIA, such that $\phi^{(n-i+1)}$ is an optimal first-step decision rule for the problem with horizon $i,$ where $i=1,2,\cdots,n.$ Here a decision rule is viewed as implementing a deterministic policy for one horizon. VIA is weakly polynomial \cite{TSENG1990287}, but it is not strongly polynomial \cite{FH}. In particular, (\ref{eq*}) and  (\ref{eq**}) do not hold for VIA because for some $\alpha\in (0,1)$ the function $\eta(\alpha)$ can be unbounded in a  neighborhood of $\alpha$ \cite[Example 3]{doi:10.1287/moor.2017.0912}.  Therefore, it is important to study $\eta(\alpha,m,q)$ for VIA. For a given $\epsilon>0,$  VIA is a strongly polynomial algorithm for finding an $\epsilon$-optimal deterministic policy, and (\ref{eq*}), (\ref{eq**}) hold \cite{CB}. In addition, for deterministic MDPs, VIA is also weakly polynomial \cite{FH}, and (\ref{eq*}), (\ref{eq**}) do not hold even for deterministic MDPs. Also, if VIA is continued after it finds an optimal policy for the first time, it may generate suboptimal policies. However, after a finite number of iterations, VIA starts to generate only optimal policies \cite{doi:10.1287/mnsc.14.5.292}. This number (for a fixed problem), denoted by $N(\alpha),$ is called the turnpike integer of the discount factor $\alpha,$ and $N(\cdot)$ is called a turnpike function. This paper describes properties of the turnpike function. This function and its upper bounds define the complexity of the value iteration algorithm. 


Turnpike functions were studied by Lewis and Paul \cite{doi:10.1287/moor.2017.0912},  where break points and degenerate points were defined. Here we say that the value of a discount factor is irregular if it is either a break point or a degenerate point. Irregular points define partition intervals on which the value function is analytic. It is proved in \cite{doi:10.1287/moor.2017.0912} that a turnpike function is bounded in any closed interval without irregular points. It is also claimed there that any closed interval without irregular points can be partitioned into finitely many intervals with the turnpike function on each of them being constant. Such intervals are called turnpike intervals. \cite{doi:10.1287/moor.2017.0912} also provided necessary conditions for a turnpike function to be bounded near an irregular point. 

In this paper, we introduce new methods for studying properties of turnpike functions and develop a more general theory. We do not assume that terminal rewards equal $0.$ Section~\ref{secpre} provides preliminary facts. In particular, Theorem~\ref{coreq} implies that, by comparing finite-step transition probabilities, it can be detected in strongly polynomial time whether two deterministic policies have equal values for all discount factors. Section~\ref{secclass} provides classification of discount factors viewed as points from the interval $[0,1),$ which is slightly different from the classification in~\cite{doi:10.1287/moor.2017.0912}. We consider regular and irregular points, and irregular points can be break points or touching points. In particular, an irregular point can be break and touching simultaneously. Algorithm~\ref{algequi} describes the strongly polynomial procedure for detecting whether a discount factor is regular.  

Section~\ref{sectp} studies properties of turnpike functions. 
Theorem~\ref{uppersemi} shows that a turnpike function is upper semicontinuous on each partition interval. Theorem~\ref{finitetpi} states that an interval can be partitioned into finitely many turnpike intervals iff the turnpike function is bounded on that interval. 
 Section~\ref{sectbb} studies conditions for boundedness of turnpike functions near irregular points. Theorem~\ref{tbnecessary} combines and rephrases \cite[Theorems 3,4]{doi:10.1287/moor.2017.0912} stating necessary conditions for boundedness, for which we provide a shorter proof. Example~\ref{exatbbounded} shows that these necessary conditions are not sufficient. Theorem~\ref{tbsufficient} provides sufficient conditions. 

Section~\ref{secsmalltp} studies VIA for MDPs with small discount factors.
Theorem~\ref{ndboundm}(a) shows that a turnpike function for small discount factors is bounded by the number of states. Following this result, within at most $O(m^2q)$ arithmetic operations, Algorithm~\ref{algsmall} computes a number $\tilde{\Delta}_L\in(0,1]$ and a set of deterministic policies $\F_L,$ such that for any policy $\phi\in\F_L$ and any $\alpha\in(0,\tilde{\Delta}_L],$ $\phi$ is optimal for the infinite-horizon problem with discount factor $\alpha.$ 

Section~\ref{secexamples} provides examples illustrating the results of this paper and counterexamples to some possible generalizations. Section~\ref{seccon} is the conclusion summarizing the results of this paper as well as subsequent directions for future studies. A summary of the nomenclature is also provided before the list of references.

\section{Preliminary}\label{secpre}
Let $\N,$ $\N^+,$ and $\R$ be the set of natural numbers including $0$, positive integers, and real numbers, respectively. We consider a Markov decision process (MDP) with finite state space $\X$ consisting of $m$ states $\X:=\{x^i\}_{i=1}^m$ and finite action space $\A:=\bigcup_{x\in\X}A(x)$, where $A(x)\ne\emptyset$ represents the set of available actions at each state $x\in\X$. Each action set $A(x)$ consists of $q_x$ actions. Thus, the total number of actions is $q:=\sum_{x=1}^m q_x,$ and this number can be interpreted as the total number of state-action pairs. For each $x\in\X$ and $a\in A(x)$ there is a \textbf{one-step reward} $r(x,a)$. Starting from some initial state $x\in\X$, a decision-maker selects an action and collects the corresponding reward, and the state is moved to the next state $y\in\X$ with probability distribution $p(y|x,a).$ If this process stops after $n$ steps at state $x_n$, then there is a \textbf{terminal reward} $s(x_n)$ (sometimes also called a salvage value), where $s:\X\to\R$. This process continues over a finite or infinite planning horizon. A \textbf{decision rule} is a mapping $\phi:\X\to\A$ such that $\phi(x)\in A(x)$ for each $x\in\X$, and a \textbf{Markov policy} $\pi$ is a sequence of decision rules $(\phi_0,\phi_1,\ldots)$, which means the decision maker applies $\phi_n$ as the decision rule at the $n$-th horizon. There are more general decision rules and policies, but for finite-state MDPs with expected total rewards, it is sufficient to consider the policies in this form; see e.g., \cite[p.~154]{Puterman:1994:MDP:528623}. Therefore, in this paper we only consider Markov policies. A policy $\pi$ is \textbf{deterministic} if $\pi=(\phi,\phi,\ldots).$ We denote such policies by $\phi$ and identify a decision rule $\phi$ with a deterministic policy $\phi.$ Throughout this paper, we do not distinguish between decision rules and deterministic policies. The term ``decision rule" is used when ``one-step policy" is emphasized, while the term ``deterministic policy" refers to the policy that applies the same one-step policy at all steps. The objective is to find a policy maximizing the expected discounted total reward. We denote the set of all Markov policies by $\Pi$, and the set of all decision rules by $\F,$ which is also the set of deterministic policies. Let $\alpha\in[0,1)$ be a \textbf{discount factor}. The total $\alpha$-discounted expected rewards for an $n$-horizon and infinite-horizon MDP under the initial state $x$ and the policy $\pi=(\phi_0,\phi_1,\ldots),$ which are also called the \textbf{objective functions} of the policy $\pi,$ are
    \begin{equation}\label{eqexpectation}
        \begin{aligned}
			v_{n,\alpha}^{\pi}(x) & :=\mathbb{E}_x^\pi\left[\sum_{t=0}^{n-1}\alpha^tr(x_t,\phi_t(x_t))+\alpha^ns(x_n)\right],\ n\in\N,\\
			v_\alpha^{\pi}(x) & :=\mathbb{E}_x^\pi\left[\sum_{t=0}^{\infty}\alpha^tr(x_t,\phi_t(x_t))\right],
		\end{aligned}
    \end{equation}
	where $x_0,\phi_0(x_0),x_1,\phi_1(x_1),\ldots,x_{t-1},\phi_{t-1}(x_{t-1}),x_t,\ldots$ is a trajectory of the process, $t$ is the integer time parameter, and $\mathbb{E}_x^\pi$ is the expectation defined by the initial state $x$ and under the policy $\pi$. A policy $\pi$ is \textbf{optimal} for an $n$-horizon problem, if $v_{n,\alpha}^\pi(x)=\sup_{\pi\in\Pi}v_{n,\alpha}^\pi(x)$, and for an infinite-horizon problem if $v_\alpha^\pi(x)=\sup_{\pi\in\Pi}v_\alpha^\pi(x)$ for all $x\in\X$. It is well-known that there are Markov optimal policies for finite-horizon problems, and deterministic optimal policies for infinite-horizon problems; see e.g., \cite[p.~154]{Puterman:1994:MDP:528623}. \textbf{Value functions} are defined as
	\begin{equation}\label{maxvaluef}
		V_{n,\alpha}(x) : = \max_{\pi\in\Pi}v_{n,\alpha}^\pi(x),\qquad V_\alpha(x) := \max_{\phi\in\F}v_\alpha^\phi(x),\quad n\in\N,\ x\in\X.
	\end{equation}
    The \textbf{optimality operators} $T_\alpha^\phi$ and $T_\alpha,$ where $\phi\in\F,\alpha\in[0,1),$ are defined as the mappings
    \begin{equation}
		\begin{aligned}\label{eqTalpha}
			& T_\alpha^\phi: \R^m\to\R^m \quad s.t. \quad T_\alpha^\phi V(x)=r(x,\phi(x))+\alpha\sum_{y\in\X}p(y|x,\phi(x))V(y), &V\in\R^m,\ x\in\X;\\
			& T_\alpha: \R^m\to\R^m \quad s.t. \quad T_\alpha V(x)=\max_{\phi\in\F}T_\alpha^\phi V(x), &V\in\R^m,\ x\in\X.
		\end{aligned}
    \end{equation}
	   Denote by $v_{n,\alpha}^\pi, v_\alpha^\pi, V_{n,\alpha}, V_\alpha, s$ the corresponding column vectors in $\R^m.$ Let $\pmb{0},\pmb{1}\in\R^m$ be the zero vector and the vector with all entries being $1$ respectively. For each $\phi\in\F$, let $P_{i,j}(\phi):=p(x^j|x^i,\phi(x^i))$ for $x^i,x^j\in\mathbb{X}$ and $r(\phi)(x):=r(x,\phi(x))$ for $x\in\mathbb{X}$. Denote by $P(\phi):=\left[P_{i,j}(\phi)\right]_{1\le i,j\le m}\in\R^{m\times m}$ the stochastic matrix and by $r(\phi)$ the column vector of one-step rewards under the decision rule $\phi$. Note that $v_{0,\alpha}^\pi=V_{0,\alpha}=s$ for all $\alpha\in[0,1)$ and for all $\pi\in\Pi$. For a policy $\pi=(\phi_0,\phi_1,\ldots),$ denote $P_t(\pi):=P(\phi_0)P(\phi_1)\ldots P(\phi_{t-1})$ for $t>0,$ $r_t(\pi):=r(\phi_t)$ for $t\ge0$, and set $P_0(\pi)=I$ to be the $m\times m$ identity matrix. Note that $P_t(\pi)$ is always a stochastic matrix for any $\pi\in\Pi$ and any $t\in\N.$ Under these notations, it is well-known that objective functions and value functions satisfy
        \begin{align}
		    v_{n,\alpha}^\pi&=T_\alpha^{\phi_0}T_\alpha^{\phi_1}\ldots T_\alpha^{\phi_{n-1}}s=\sum_{t=0}^{n-1}\alpha^tP_t(\pi)r_t(\pi)+\alpha^nP_n(\pi)s, &\text{for all } n\in\N^+;\label{eqvaluesum1}\\ 
            v_\alpha^\pi&=\lim_{n\to\infty}T_\alpha^{\phi_0}T_\alpha^{\phi_1}\ldots T_\alpha^{\phi_{n-1}}V=\sum_{t=0}^\infty\alpha^tP_t(\pi)r_t(\pi), &\text{for all } V\in\R^m;\label{eqvaluesum2}\\
            v_\alpha^\phi&=T_\alpha^\phi v_\alpha^\phi=\sum_{t=0}^\infty\alpha^tP^t(\phi)r(\phi)={\left[I-\alpha P(\phi)\right]}^{-1}r(\phi);\label{eqvaluesum3}\\
            V_{n,\alpha} & = T_\alpha V_{n-1,\alpha}=\left(T_\alpha\right)^2V_{n-2,\alpha}=\ldots=\left(T_\alpha\right)^ns,\quad V_\alpha = T_\alpha V_\alpha, &\text{for all } n\in\N^+,\label{eqoptimal}
	    \end{align}
    where equations $\eqref{eqoptimal}$ are called the \textbf{optimality equations}. For each discount factor $\alpha\in[0,1)$ and for each $n\in\N,$ a deterministic policy $\phi$ is first-step-optimal for the $n$-horizon problem if and only if $V_{n,\alpha}=T_\alpha^\phi V_{n-1,\alpha},$ and it is optimal for the infinite-horizon problem if and only if $V_\alpha=T_\alpha^\phi V_\alpha.$ 
    If decision rules $\phi$ and $\psi$ are both optimal for an infinite-horizon problem, then they can be applied interchangeably to attain the optimal discounted expected total reward for the infinite-horizon problem in view of \eqref{eqvaluesum2} and the optimality equations \eqref{eqoptimal}. However, this does not mean that they can be applied interchangeably for optimizing finite-horizon problems. 
    
The following theorem provides necessary and sufficient conditions for two deterministic policies to have equal expected total discounted infinite-horizon rewards for all initial states $x\in\X$ and for all discount factors $\alpha\in [0,1).$ Theorem~\ref{coreq} implies that its conditions (a) and (b) are equivalent.

\begin{thm}\label{coreq} Let $\{\phi,\psi\}\subset\F$. Each of the following two conditions is necessary and sufficient for $v_\alpha^{\phi}=v_\alpha^{\psi}$ for all $\alpha\in[0,1):$
	\begin{enumerate}[(a)]
		\item $P^t(\phi)r(\phi)= P^t(\psi)r(\psi)$ for all $t\in\N;$
		\item $P^t(\phi)r(\phi)=P^t(\psi)r(\psi)$ for all $t=0,1,\ldots,\min\{G(\phi,r(\phi)),G(\psi,r(\psi))\}-1.$
	\end{enumerate}
\end{thm}

To prove Theorem~\ref{coreq} we need the following constructions. Let $P^0:=I$ for $P\in\R^{m\times m}.$ Then $P^0\mathbf{v}=\mathbf{v}.$ For
$\phi\in\F$ and $\mathbf{v}\in\mathbb{R}^m,$ let $G(\phi,\mathbf{v}):=1,$ if $\mathbf{v}=c\pmb{1}$ for some $c\in\R,$ and
	\begin{equation}\label{eqG}
		G(\phi,\mathbf{v}):=\max\left\{t\in\N^+: \left\{\pmb{1},P^0(\phi)\mathbf{v},P^1(\phi)\mathbf{v},\ldots,P^{t-2}(\phi)\mathbf{v}\right\}\ \text{are linearly independent}\right\},
	\end{equation}
    if  $\mathbf{v}\ne c\pmb{1}$ for all $c\in\R.$ 
Clearly $2\le G(\phi,\mathbf{v})\le m$ iff  $\mathbf{v}\ne c\pmb{1}$ for all $c\in\R.$ In particular, $G(\phi,\mathbf{v})=1$ if $m=1.$ Thus, $1\le G(\phi,\mathbf{v})\le m.$ The following proposition, which follows from linear algebra arguments, is used in the proof of Theorem~\ref{coreq}.
 
\begin{pro}\label{thmla} Let $\phi^{(1)},\phi^{(2)}\in\F$ and $\mathbf{v}_1,\mathbf{v}_2\in\mathbb{R}^m.$ Then $P^t(\phi^{(1)})\mathbf{v}_1=P^t(\phi^{(2)})\mathbf{v}_2$ for all $t\in\N$ iff $P^t(\phi^{(1)})\mathbf{v}_1=P^t(\phi^{(2)})\mathbf{v}_2$ for all $t=0,1,\ldots,\min\{G(\phi^{(1)},\mathbf{v}_1),G(\phi^{(2)},\mathbf{v}_2)\}-1.$
\end{pro}
\begin{proof}The necessity is true because, if the equality holds for all natural numbers, then it also holds for any finite set of natural numbers. Let us prove the sufficiency. Without loss of generality, we assume $G:=G(\phi^{(1)},\mathbf{v}_1)=\min\{G(\phi^{(1)},\mathbf{v}_1),G(\phi^{(2)},\mathbf{v}_2)\}$. We observe that $\mathbf{v}_1=\mathbf{v}_2$ since $G\ge 1.$ If $G=1$, then $\mathbf{v}_1=c\pmb{1}$ for some $c\in\R$, which together with $\mathbf{v}_1=\mathbf{v}_2$ implies $P^t(\phi^{(1)})\mathbf{v}_1=P^t(\phi^{(2)})\mathbf{v}_2=c\pmb{1}$ for all $t\in\N$.  If $G\ge 2,$ let 
	\[B:=\begin{bmatrix} \pmb{1} & P^0(\phi^{(1)})\mathbf{v}_1 & \ldots & P^{G-2}(\phi^{(1)})\mathbf{v}_1 \end{bmatrix}=\begin{bmatrix} \pmb{1} & P^0(\phi^{(2)})\mathbf{v}_2 & \ldots & P^{G-2}(\phi^{(2)})\mathbf{v}_2 \end{bmatrix}.\]
	Note that $B\in\R^{m\times G}$. Let us prove that 
	\begin{equation}\label{eqPTphipPTpsi} P^t(\phi^{(1)})\mathbf{v}_1=P^t(\phi^{(2)})\mathbf{v}_2=B\mathbf{c}_t\qquad \text{for all } t\in\N,
	\end{equation}
	where $\mathbf{c}_t\in\R^G$ such that 
	\begin{equation}\label{eqg-1GHEF}
		\mathbf{c}_{G-1}={\bf \lambda}\quad  {\rm for\ some\ } \mathbf{\lambda}\in\R^G,
	\end{equation}
	and
	\begin{equation}\label{eqctGHEF}\mathbf{c}_t=\begin{cases}
		\mathbf{e}_{t+2},&\ \text{if $t\le G-2,$}\\
		\begin{bmatrix} \mathbf{e}_1 & \mathbf{e}_3 & \ldots & \mathbf{e}_G & \mathbf{c}_{G-1}\\ \end{bmatrix}^{t-G+1}\mathbf{c}_{G-1},&\ \text{if $t\ge G-1.$}
	\end{cases}\end{equation}
	Here each $\mathbf{e}_i\in\R^G, i=1,\ldots,G,$ is the column vector with the $i$-th entry being $1$ and all other entries being $0$. We observe that $\pmb{1}=B\mathbf{e}_1.$ Equation \eqref{eqg-1GHEF} and  equation \eqref{eqctGHEF} for $t=0,1,\ldots,G-2$ follow from the definition of $G$ and from the assumptions of the theorem. Equation \eqref{eqctGHEF} is an identity when $t=G-1$ since $t-G+1=0.$ Suppose that \eqref{eqctGHEF} holds for some $t\ge G-1$. Then 
    \begin{equation*}
	\begin{aligned}
	P^{t+1}(\phi^{(1)})\mathbf{v}_1=P(\phi^{(1)})P^t(\phi^{(1)})\mathbf{v}_1=P(\phi^{(1)})B\mathbf{c}_t&=\begin{bmatrix} \pmb{1} & P(\phi^{(1)})\mathbf{v}_1 & \ldots & P^{G-1}(\phi^{(1)})\mathbf{v}_1 \end{bmatrix}\mathbf{c}_t\\
		&=B\begin{bmatrix} \mathbf{e}_1 & \mathbf{e}_3 & \ldots & \mathbf{e}_G & \mathbf{c}_{G-1}\\ \end{bmatrix}\mathbf{c}_t=B\mathbf{c}_{t+1},
		\end{aligned}
        \end{equation*}
	where $\mathbf{c}_{t+1}:=\begin{bmatrix} \mathbf{e}_1 & \mathbf{e}_3 & \ldots & \mathbf{e}_G & \mathbf{c}_{G-1}\\ \end{bmatrix}\mathbf{c}_t,$  and the same formulae hold if $\phi^{(1)}$ is substituted with $\phi^{(2)}.$ Thus, equation~\eqref{eqPTphipPTpsi} is proved by induction.
\end{proof}


\begin{proof}{Proof of Theorem~\ref{coreq}.} Necessity of (a,b) follows from taking derivatives with respect to $\alpha$ in \eqref{eqvaluesum3} and letting $\alpha=0.$ Sufficiency of (a) follows from \eqref{eqvaluesum3}, and sufficiency of (b) is implied by the sufficiency of (a) and by Proposition~\ref{thmla}.
\end{proof}

\section{Classification of Discount Factors}\label{secclass}
In this section, we classify different types of discount factors based on the properties of sets of optimal policies and provide relevant results. 
For each discount factor $\alpha\in [0,1),$  the following definition introduces the set $M_n(\alpha)$ of optimal Markov policies for an $n$-horizon problem, the set $D_n(\alpha)$ of first-step-optimal decision rules for an $n$-horizon problem, and the set $D(\alpha)$  of deterministic optimal policies for the infinite-horizon problem, where $n\in\N^+.$ It also defines the corresponding sets  $M_n(E),$  $D_n(E),$ and $D(E),$ where $E\subset [0,1].$

\begin{defi}\label{defMnDnD} \emph{For $\alpha\in[0,1),$  $E\subset[0,1),$ and $n\in\N^+,$ let $M_0(\alpha):=\Pi$ and define
		$M_n(\alpha)  := \{\pi\in\Pi:v_{n,\alpha}^\pi=V_{n,\alpha}\},$ 
		$D_n(\alpha)  := \{\phi\in\F:V_{n,\alpha}=T_\alpha^\phi V_{n-1,\alpha}\}$, 
		$D(\alpha)  := \{\phi\in\F:V_\alpha=v_\alpha^\phi\},
	$ and
    \[M_n(E):=\bigcup_{\alpha\in E}M_n(\alpha),\qquad D_n(E):=\bigcup_{\alpha\in E}D_n(\alpha),\qquad D(E):=\bigcup_{\alpha\in E}D(\alpha).\]}
\end{defi}

If $E\ne\emptyset,$ where $E\subset[0,1),$ then $M_n(E)\ne\emptyset$ and $D_n(E)\ne\emptyset$ for all $n\in\N^+,$ and $D(E)\ne\emptyset.$ 
We observe that, if $E_1\subset E_2$, then $M_n(E_1)\subset M_n(E_2),$ $D_n(E_1)\subset D_n(E_2)$ and $D(E_1)\subset D(E_2), $ which means $M_n,$ $D_n,$ and $D$ are non-decreasing set-valued functions. Note that $D_n(\alpha)$ is also the set of policies output by the VIA at the $n$-th iteration. According to \cite[Theorem 4]{doi:10.1287/mnsc.14.5.292},  $\limsup_{n\to\infty}D_n(\alpha)\subset D(\alpha)$ for all $\alpha\in[0,1),$ which allows us to define integers $N(\alpha)$ such that starting from the $N(\alpha)$-th iteration, the VIA generates deterministic optimal policies for infinite-horizon problems.

\begin{defi}\label{deftp} \emph{The \textbf{turnpike integer} of discount factor $\alpha$ is the smallest positive integer $N(\alpha)$ such that $D_n(\alpha)\subset D(\alpha)$ for all $n\ge N(\alpha)$. We say $N(\alpha)$ is the \textbf{turnpike function}. We also define  $N^*(E):=\sup_{\alpha\in E}N(\alpha)$ for $E\subset [0,1).$}
\end{defi}

The ``smallest integer" in the definition above implies that, if $N(\alpha)\ge2,$ then $D_{N(\alpha)-1}(\alpha)\setminus D(\alpha)\ne\emptyset$. Note that $N(0)=1,$ and the turnpike integer for a fixed nonzero discount factor may depend on terminal rewards. For a fixed discount factor $\alpha,$ it is possible that the sequence $\{D_n(\alpha)\}_{n=1}^\infty$ does not converge; \cite[Appendix C]{TSENG1990287}. 

Let $\alpha\in(0,1)$ and $n\in\N^+.$ The monotonicity of the set-valued functions $M_n,$ $D_n,$ and $D$ allows us to define $M_n(\alpha-)$, $D_n(\alpha-),$ and $D(\alpha-),$ which are the set of left optimal Markov policies for the $n$-horizon problem, the set of first-step-left optimal decision rules for the $n$-horizon problem, and the set of left optimal deterministic policies, respectively, which means they are optimal if the discount factor belongs to an interval $(\alpha-\epsilon,\alpha)$ for some $\epsilon>0.$ Similarly, $M_n(\alpha+)$, $D_n(\alpha+),$ and $D(\alpha+)$ are the set of right optimal Markov policies for the $n$-horizon problem, the set of first-step-right optimal decision rules for the $n$-horizon problem, and the set of right optimal deterministic policies, respectively, which means they are optimal if the discount factor belongs to an interval $(\alpha,\alpha+\epsilon)$ for some $\epsilon>0.$ The following definition formally defines these sets.


\begin{defi}\label{def+-}\emph{
	For $n\in\N^+$ and for $\alpha\in[0,1)$, define $M_n(0-)=D_n(0-)=D(0-):=\emptyset,$ and
    \begin{equation*}
	\begin{aligned}
		&M_n(\alpha-) := \bigcap_{i=1}^\infty M_n\left(\left(\alpha-\frac{1}{i},\alpha\right)\bigcap\left[0,1\right)\right), &M_n(\alpha+) := \bigcap_{i=1}^\infty M_n\left(\left(\alpha,\alpha+\frac{1}{i}\right)\bigcap\left[0,1\right)\right),\\
		&D_n(\alpha-) := \bigcap_{i=1}^\infty D_n\left(\left(\alpha-\frac{1}{i},\alpha\right)\bigcap\left[0,1\right)\right), &D_n(\alpha+) := \bigcap_{i=1}^\infty D_n\left(\left(\alpha,\alpha+\frac{1}{i}\right)\bigcap\left[0,1\right)\right),\\
		&D(\alpha-) := \bigcap_{i=1}^\infty D\left(\left(\alpha-\frac{1}{i},\alpha\right)\bigcap\left[0,1\right)\right), &D(\alpha+) := \bigcap_{i=1}^\infty D\left(\left(\alpha,\alpha+\frac{1}{i}\right)\bigcap\left[0,1\right)\right).
	\end{aligned}
    \end{equation*}}
\end{defi}

\begin{figure}[ht!]
	\centerline{\includegraphics[scale=0.35]{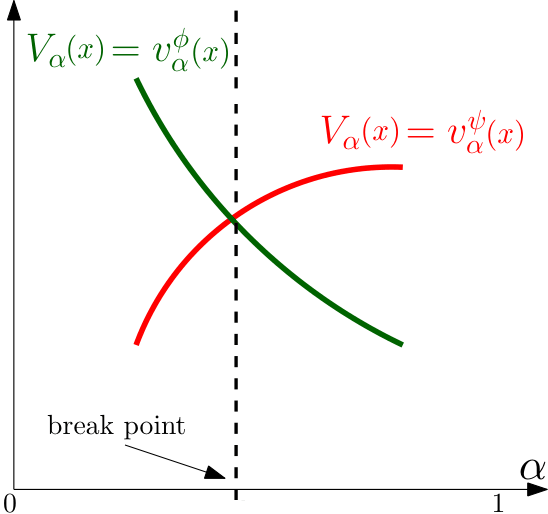}\qquad\includegraphics[scale=0.38]{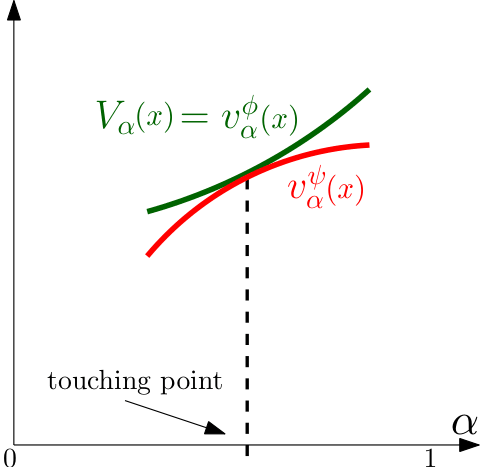}\qquad\includegraphics[scale=0.38]{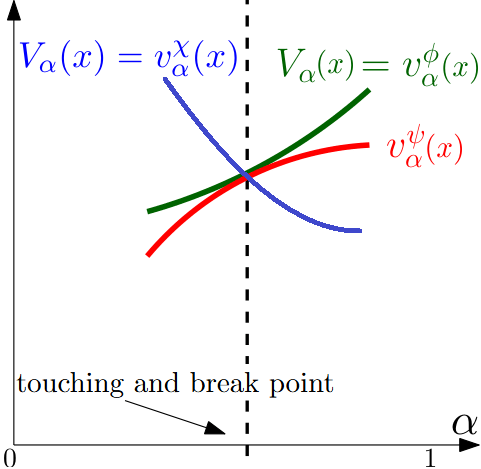}}
	\caption{Three possible cases of irregular points: a break point, a touching point, and a point which is touching and break.}\label{fig:btg}
\end{figure}

The sets of left and right optimal policies may be different for some discount factor $\alpha\in (0,1)$, which implies that the corresponding set-valued function is discontinuous at the point $\alpha.$  There are two types of discontinuity: break and touching, and they are described in Definition~\ref{defirr}. At a break point, the set of left optimal policies is different from the set of right optimal policies; at a touching point, the set of optimal policies does not equal the union of the set of left optimal policies and the set of right optimal policies; see the first two graphs in Figure~\ref{fig:btg}. 


\begin{defi}\label{defirr} \emph{We say that $\alpha\in(0,1)$ is a \textbf{break point} if $D(\alpha-)\ne D(\alpha+)$. We say $\alpha\in[0,1)$ is a \textbf{touching point} if $D(\alpha)\ne D(\alpha-)\cup D(\alpha+)$. If $\alpha\in[0,1)$ is neither a break point nor a touching point, then it is called a \textbf{regular point}. If a point is not regular, then it is an \textbf{irregular point}.}
\end{defi}

Similarly, the following definition formally introduces two types of discontinuous points of the set-valued functions $D_n$ and $M_n.$

\begin{defi}\label{defirrfinite} \emph{We say that $\alpha\in(0,1)$ is an \textbf{$\pmb{n}$-horizon break point} if $M_n(\alpha-)\ne M_n(\alpha+)$. We say that $\alpha\in[0,1)$ is an \textbf{$\pmb{n}$-horizon touching point} if $M_n(\alpha)\ne M_n(\alpha-)\cup M_n(\alpha+)$. If $\alpha\in[0,1)$ is neither an $n$-horizon break point nor an $n$-horizon touching point, then it is called an \textbf{$\pmb{n}$-horizon regular point}. If a point is not $n$-horizon regular, then it is an \textbf{$\pmb{n}$-horizon irregular point}. Similarly, we define \textbf{$\pmb{n}$-horizon-first-step break points}, \textbf{$\pmb{n}$-horizon-first-step touching points}, \textbf{$\pmb{n}$-horizon-first-step regular points}, and \textbf{$\pmb{n}$-horizon-first-step irregular point} by replacing the function $M_n$ with $D_n$.}
\end{defi}
\medskip

In Definitions~\ref{defirr} and \ref{defirrfinite}, regular points and irregular points are the points where the corresponding set-valued function is continuous and discontinuous, respectively. Note that a break point, an $n$-horizon break point, or an $n$-horizon-first-step break point may also be a touching point, an $n$-horizon touching point, or an $n$-horizon-first-step touching point, respectively; see the third graph in Figure~\ref{fig:btg}. According to Definition~\ref{defirr}, point 0 is irregular if it is touching, and it is regular otherwise. A point $\alpha\in(0,1)$ is regular iff $D(\alpha)=D(\alpha-)=D(\alpha+),$ and $0$ is a regular point iff $D(0+)=D(0)$. Similar relations also hold if the set-valued function $D$ is replaced with $M_n$ or $D_n$.

For a fixed $n\in\N$, the following proposition shows that the set of irregular points, the set of $n$-horizon-first-step irregular points, and the set of $n$-horizon irregular points, are finite.
\begin{pro}\label{partition} (a) For each $x\in\X,$ $\pi\in\Pi,$ and $n\in\N$, the function $v_{n,\alpha}^\pi(x)$ is  polynomial in $\alpha$, and the function $V_{n,\alpha}(x)$ is  continuous and piecewise polynomial in $\alpha$. There is a unique finite set of points $\{a_i^{(n)}\}_{i=0}^{l(n)}\subset[0,1]$ such that $a_0^{(n)}=0,$  $a_{l(n)}^{(n)}=1,$  $a_i^{(n)}<a_{i+1}^{(n)}$ for $0\le i\le l(n)-1$, the function $M_n(\alpha)$ is constant on each open interval $\mathcal{I}_i^{(n)}:=(a_i^{(n)},a_{i+1}^{(n)})$ for $0\le i\le l(n)-1,$ but it has jumps at the points $\{a_i^{(n)}\}_{i=1}^{l(n)-1}$ and  is not constant on each of half-open intervals $(a_i^{(n)},a_{i+1}^{(n)}],$ $[a_{i+1}^{(n)},a_{i+2}^{(n)})$ for $0\le i\le l(n)-2.$ In addition, $V_{n,\alpha}=v_{n,\alpha}^\pi$ for all $\pi\in M_n(\mathcal{I}_i^{(n)}),\alpha\in\mathcal{I}_i^{(n)},$ and $0\le i\le l(n)-1.$\smallskip\\
(b) For each $x\in\X,$ $\phi\in\F,$ and $n\in\N^+$, the function $T_\alpha^\phi V_{n-1,\alpha}(x)$ is  continuous and piecewise polynomial in $\alpha$. There is a unique finite set of points $\{b_i^{(n)}\}_{i=0}^{\tilde{l}(n)}\subset[0,1]$ such that $b_0^{(n)}=0,$  $b_{\tilde{l}(n)}^{(n)}=1,$  $b_i^{(n)}<b_{i+1}^{(n)}$ for $0\le i\le \tilde{l}(n)-1$, and the function $D_n(\alpha)$ is constant on each open interval $\tilde{\mathcal{I}}_i^{(n)}:=(b_i^{(n)},b_{i+1}^{(n)})$ for $0\le i\le \tilde{l}(n)-1,$ but it has jumps at the points $\{b_i^{(n)}\}_{i=1}^{\tilde{l}(n)-1}$ and is not constant on each of half-open intervals $(b_i^{(n)},b_{i+1}^{(n)}],$ $[b_{i+1}^{(n)},b_{i+2}^{(n)})$ for $0\le i\le \tilde{l}(n)-2.$ In addition, $T_\alpha^\phi V_{n-1,\alpha}=T_\alpha^\psi V_{n-1,\alpha}$ for all $\phi,\psi\in D_n(\tilde{\mathcal{I}}_i^{(n)}),\alpha\in\tilde{\mathcal{I}}_i^{(n)},$ and $0\le i\le \tilde{l}(n)-1.$\smallskip\\
(c) For each $x\in\X$ and $\phi\in\F$, the function $v_\alpha^\phi(x)$ is rational  in $\alpha,$ and the function $V_\alpha(x)$ is continuous and piecewise rational in $\alpha$. There is a unique finite set of points $\{a_i\}_{i=0}^l\subset[0,1]$ such that $a_0=0,$  $a_l=1,$  $a_i<a_{i+1}$ for $0\le i\le l-1$, the function $D(\alpha)$ is constant on each open interval $\mathcal{I}_i:=(a_i,a_{i+1})$ for $0\le i\le l-1,$ but it has jumps at the points $\{a_i\}_{i=1}^l$ and is not constant on each of half-open intervals $(a_i,a_{i+1}],$ $[a_{i+1},a_{i+2})$ for $0\le i\le l-2.$ In addition, $V_\alpha=v_\alpha^\phi$ for all $\phi\in D(\mathcal{I}_i),\alpha\in\mathcal{I}_i,$ and $0\le i\le l-1.$
\end{pro}

\begin{proof}
(a) follows from \eqref{eqvaluesum1}, \cite[Proposition 3.1.4]{sennott1999stochastic}, and the finiteness of the number of Markov policies for a finite-horizon problem; (b) follows from (a) and \eqref{eqTalpha}; (c) follows from the same arguments as (a) because, in view of \eqref{eqvaluesum3}, the functions $v_\alpha^\phi(x)$ are rational in $\alpha.$ In addition, the  intervals ${\mathcal{I}}_i^{(n)}$
 and $\mathcal{I}_i$ indicated in statements (a) and (c) are the intervals on which the corresponding functions are polynomial or rational. This fact for statement (a) implies that the function $T_\alpha^\phi V_{n-1,\alpha}(x)$ mentioned in (b) is also polynomial on a finite set of intervals covering  $[0,1).$
\end{proof}

Since $D(\alpha)=\times_{x\in\X}A_\alpha(x),$ where $A_\alpha (x)=\{a\in A(x):r(x,a)+\alpha\sum_{y\in\X}p(y|x,a)V(y)=V(x)\},$ the set $D(\alpha)$ can be represented by the sets $A_\alpha(x)$ for all $x\in\X.$ Thus, each set $D(\alpha)$ can be represented by no more than $q$ symbols if $A(x)\cap A(z)=
\emptyset$ for all distinct states $x, z\in\X$. In view of Proposition~\ref{partition}(c), $\alpha\in[0,1)$ is a regular point iff $v_\beta^\phi=v_\beta^\psi$ for all $\beta\in[0,1)$ and all $\phi,\psi\in D(\alpha).$ In general, finding an irregular point is difficult. However, for a given discount factor $\alpha,$ it can be verified in strongly polynomial time whether $\alpha$ is an irregular point. First, $D(\alpha)$ can be computed using PIA in strongly polynomial time \cite{PIS}. Then, with $D(\alpha)$ being known, the following algorithm, which takes at most $m$ iterations and each iteration takes at most $O(mq)$ arithmetic operations, verifies whether $P^t(\phi)r(\phi)$ is the same for all $\phi\in D(\alpha)$ for each $t=0,1,\cdots,m-1,$ which by Theorem~\ref{coreq} implies $v_\beta^\phi=v_\beta^\psi$ for all $\beta\in[0,1)$ and all $\phi,\psi\in D(\alpha).$ In view of Proposition~\ref{partition}(c), this algorithm also verifies whether $\alpha$ is an irregular point.

\begin{algorithm}[H]
\SetAlgoLined
\LinesNumbered
\caption{Detect whether a discount factor $\alpha\in [0,1)$ is an irregular point, assuming that $D(\alpha)$ is known.}\label{algequi}
\KwIn{$\X,$ $\{A(x)\}_{x\in\X},$ $\{r(x,a)\}_{x\in\X,\ a\in A(x)},$ $\{p(y|x,a)\}_{x,y\in\X,\ a\in A(x)},$ $D(\alpha)$}
\KwOut{whether or not $\alpha$ is an irregular point.}
Set $n=1,$ $A_\alpha(x)=\{a\in A(x):a=\phi(x)\ \text{for some $\phi\in D(\alpha)$}\}$  and $u(x,a)=r(x,a)$ for all $x\in\X$ and for all $a\in A_\alpha(x)$\;
 \While{$n<m$ and $\{u(x,a)\}_{a\in A_\alpha(x)}$ is a singleton for every $x\in\X$}{
$n\leftarrow n+1,$ $u(x)\leftarrow u(x,a)$ for any $a\in A_\alpha(x)$\;
 \lFor{$x\in\X$ and $a\in A_\alpha(x)$}{
     $u(x,a)\leftarrow\sum_{y\in\X}p(y|x,a)u(y)$
  }
 }
 \leIf{$n=m$}{
   Output: $\alpha$ is a regular point
   }
   {Output: $\alpha$ is an irregular point
  }
\end{algorithm}

\begin{rem}\label{remirr}
\emph{This remark deals with an infinite-horizon problem. }
	\smallskip
    
\noindent\emph{(i) Though the same notation $\{a_i\}_{i=1}^{l-1}$ is used in this paper and in \cite[p.~1148]{doi:10.1287/moor.2017.0912}, this notation is used for different sets.  In \cite{doi:10.1287/moor.2017.0912} this notation denotes the set of break points.  In this paper it denotes the set of irregular points, which is also the union of the sets of break and touching points.}
\smallskip

\noindent\emph{(ii) If $\alpha\in[0,1)$ is a touching point, then there exists a deterministic policy optimal at $\alpha$, but it is strictly suboptimal at all other points in some neighborhood of $\alpha$. Since the value function is greater than or equal to the objective function for a policy, this means that the graphs of the value function and the objective function for a deterministic policy touch at the point $\alpha.$}
\smallskip

\noindent\emph{(iii) The definition of a degenerate point in \cite[p.~1148]{doi:10.1287/moor.2017.0912} corresponds to the definition of a non-break touching point in this paper. However, the definition of a degenerate point in \cite[p.~1148]{doi:10.1287/moor.2017.0912} is misstated because it claims that for a degenerate point $\alpha$ from an interval $I_i,$ on which some deterministic policy is optimal, there is another deterministic policy, which is optimal for the discount factor $\alpha,$ and it is suboptimal for all other points in $I_i$.  In other words, the point $\alpha$ is touching, and there is no other touching point in $I_i$ corresponding to the same policy. However, an interval $I_i$ described above can have multiple touching points corresponding to the same deterministic policy.  The correct version of the definition of a degenerative point in \cite[p.~1148]{doi:10.1287/moor.2017.0912} should assume the existence of the interval $I_i$ with the described properties, rather than claiming $I_i$ is an interval on which some deterministic policy is optimal.} 
\smallskip

\noindent\emph{(iv) In the classical definition, a policy $\phi\in\F$ is said to be Blackwell optimal, if there exists $\gamma'\in[0,1)$ such that $\phi$ is optimal for all discount factors $\gamma\in[\gamma',1)$, and let the smallest such $\gamma'$ be denoted by $\gamma(\phi)$. \cite{BWD} defines the Blackwell discount factor $\gamma_{\mathrm{bw}}\in[0,1)$ to be the smallest discount factor such that every $\phi\in D((\gamma_{\mathrm{bw}},1))$ is Blackwell optimal. For a Blackwell optimal policy $\phi$, it is possible that $\gamma(\phi)<\gamma_{\mathrm{bw}}$ because there may exist some touching point $\gamma_*\in (\gamma(\phi),1);$ see \cite[Example 3.6]{BWD}. Note that, according to the definitions in this paper, $\gamma_{\mathrm{bw}}$ is the largest irregular point.}
\end{rem}

\begin{defi}\label{defpartition} \emph{The partition of the interval $[0,1)$ into the sets of open intervals $\{\mathcal{I}_i\}_{i=0}^{l-1}$ and the sets of points $\{a_i\}_{i=0}^{l-1}$ introduced in Proposition~\ref{partition}(c) is called the \textbf{canonical partition}. We call the partition of the interval $[0,1)$ into the sets of open intervals $\{\mathcal{I}_i^{n}\}_{i=0}^{l(n)-1}$ and the set of points $\{a_i^{(n)}\}_{i=0}^{l(n)-1}$ in Proposition~\ref{partition}(a) as the \textbf{$\pmb{n}$-horizon canonical partition}. Each open interval $\mathcal{I}_i$ or $\mathcal{I}_i^{n}$ is called a \textbf{partition interval} or an \textbf{$\pmb{n}$-horizon partition interval} respectively.}
\end{defi}

The points $\{a_i\}_{i=1}^{l-1}$, $\{b_i^{(n)}\}_{i=1}^{\tilde{l}(n)-1}$, and $\{a_i^{(n)}\}_{i=1}^{l(n)-1}$ are irregular points, $n$-horizon-first-step irregular points, and $n$-horizon irregular points, respectively. For each $\alpha\in[0,1)$, we have $D(\alpha)\supset D(\alpha-)\cup D(\alpha+)$, $D_n(\alpha)\supset D_n(\alpha-)\cup D_n(\alpha+)$, and $M_n(\alpha)\supset M_n(\alpha-)\cup M_n(\alpha+)$. If $\alpha$ is a break point, then $D(\alpha-)\cap D(\alpha+)=\emptyset,$ and, if $\alpha$ is also non-touching, then $D(\alpha)=D(\alpha-)\cup D(\alpha+).$  If $\alpha$ is an $n$-horizon break point, then $M_n(\alpha-)\cap M_n(\alpha+)=\emptyset,$ and $M_n(\alpha)=M_n(\alpha-)\cup M_n(\alpha+)$ if $\alpha$ is also non-touching. Let $\alpha$ be an $n$-horizon-first-step break point. According to Definition~\ref{defirrfinite}, $\alpha$ is also an $n$-horizon break point, and therefore $\{b_i^{(n)}\}_{i=1}^{\tilde{l}(n)-1}\subset\{a_i^{(n)}\}_{i=1}^{l(n)-1}.$ However, if $\alpha$ is also an $(n-1)$-horizon irregular point, that is, if $\alpha\in\{a_i^{(n-1)}\}_{i=1}^{l(n-1)-1},$ then, according to Example~\ref{exaDn+-}, it is possible that $D_n(\alpha-)\cap D_n(\alpha+)\ne\emptyset.$  The following proposition is relevant to these observations.

\begin{pro}\label{upperhemi} Let $n\in\N^+$ and $\alpha\in[0,1).$ Then there exists $\delta>0$ such that $D_n(\beta)\subset D_n(\alpha),$ $M_n(\beta)\subset M_n(\alpha),$ and $D(\beta)\subset D(\alpha)$ for all $\beta\in(\alpha-\delta,\alpha+\delta)\cap[0,1)$. 
\end{pro}
\begin{proof} We prove this statement for the set-valued function $D(\alpha);$ the proofs of the other two are similar. Suppose this is not true, that is, there exists a sequence of discount factors $\{\alpha_i\}_{i=1}^\infty\subset[0,1)$ such that $\alpha_i\to\alpha$ as $i\to\infty$ and $D(\alpha_i)\not\subset D(\alpha)$ for all $i\in\N^+.$ Then for each $i\in\N^+$ there exists $\phi^{(i)}\in D(\alpha_i)\setminus D(\alpha).$ Since $D(\alpha_i)\setminus D(\alpha)\subset\F$ is finite, there exists a subsequence $\phi^{(i_j)}=\phi\in D(\alpha_i)\setminus D(\alpha).$ Then $v_{\alpha_{i_j}}^\phi=V_{\alpha_{i_j}}.$ By taking $j\to\infty$ and by the continuity of the function $V_\alpha$ in $\alpha$ from Proposition~\ref{partition}(c), we have $v_{\alpha}^\phi=V_{\alpha},$ which implies $\phi\in D(\alpha).$ This is a contradiction.
\end{proof}

\begin{rem}
    \emph{We consider the discrete topologies on $\F.$  We view $D_n(\alpha)$ and $D(\alpha)$ as functions from $[0,1)$ to $2^\F$  for $n\in\N^+,$ where $2^X$ denotes the set of subsets of the set $X.$ Proposition~\ref{upperhemi} implies that $D_n,$ $M_n,$ and $D$ are upper hemicontinuous set-valued functions.}
\end{rem}

\begin{defi}\label{defsubint} Let $\I$ be the set of all nonempty subintervals of $[0,1)$ including singletons. That is, $\I=$ 
\begin{equation*}\left\{(a,b):0\le a<b\le 1\right\}\cup\left\{[a,b]:0\le a\le b< 1\right\}\cup\left\{(a,b]:0\le a<b< 1\right\}\cup\left\{[a,b):0\le a<b\le 1\right\}.
\end{equation*}
\end{defi}

\begin{defi}\label{deftpi}
\emph{We say $\mathcal{I}\in\I$ is a \textbf{turnpike interval} if $N(\alpha)$ is constant on $\mathcal{I}.$ }
\end{defi}

\begin{defi}\label{defdiscon}
    \emph{For $\mathcal{I}\in\I$ and for each $n\in\N^+,$ let us define the following sets:
\begin{equation*}
\begin{aligned}
	\mathbb{D}^-(\mathcal{I})&:=\left\{\alpha\in\mathcal{I}\setminus\left\{0\right\}: \text{$N(\alpha)$ is not left-continuous at $\alpha$} \right\},\\
	\mathbb{D}^+(\mathcal{I})&:=\left\{\alpha\in\mathcal{I}: \text{$N(\alpha)$ is not right-continuous at $\alpha$}\right\},\\
    \mathbb{D}(\mathcal{I})&:=\mathbb{D}^-(\mathcal{I})\cup \mathbb{D}^+(\mathcal{I}),\quad\Hat{\mathbb{D}}(\mathcal{I}):=\mathbb{D}^-(\mathcal{I})\cap \mathbb{D}^+(\mathcal{I}),\\
	\P(\mathcal{I})&:=\left\{\alpha\in\mathcal{I}: \text{$\alpha$ is an irregular point}\right\},\\
	\P_n(\mathcal{I})&:=\left\{\alpha\in\mathcal{I}: \text{$\alpha$ is an $n$-horizon-first-step irregular point}\right\},\\
	\P^T_n(\mathcal{I})&:=\left\{\alpha\in\mathcal{I}: \text{$\alpha$ is an $n$-horizon-first-step touching point}\right\},\\
	\P^B_n(\mathcal{I})&:=\left\{\alpha\in\mathcal{I}: \text{$\alpha$ is an $n$-horizon-first-step break point}\right\}.
\end{aligned}
\end{equation*}}
\end{defi}

We use the notation $\P$ for the sets of irregular points because irregular points can also be called pathological. For each $\mathcal{I}\in\I,$ note that $\mathbb{D}(\mathcal{I})$ is the set of points in $\mathcal{I}$ at which the turnpike function $N(\alpha)$ is discontinuous, and $\Hat{\mathbb{D}}(\mathcal{I})$ is the set of points in $\mathcal{I}$ at which $N(\alpha)$ is neither left-continuous nor right-continuous. By Proposition~\ref{partition}(b,c), for each $n\in\N^+$ and each $\mathcal{I}\in\I$, the sets $\P(\mathcal{I}),$ $\P_n(\mathcal{I}),$ $\P^T_n(\mathcal{I}),$ and $\P^B_n(\mathcal{I})$ are finite. Note that $\P_n(\mathcal{I})=\P^T_n(\mathcal{I})\cup\P^B_n(\mathcal{I})$ for each $\mathcal{I}\in\I$. 

\section{Properties of the Turnpike Function}\label{sectp}
This section studies properties of the turnpike function. Theorem~\ref{utt} shows that the turnpike function is bounded on each closed interval without irregular points. Proposition~\ref{uttirr} extends this fact to all closed intervals.  Theorem~\ref{uppersemi} shows that the turnpike function is upper semicontinuous at each regular point. Proposition~\ref{Ndiff} provides necessary and sufficient conditions for the turnpike function to be discontinuous at a regular point. Proposition~\ref{Nd} shows that at a regular point $\alpha$, at which the turnpike function is not continuous, there exists an optimal decision rule such that it is also $(N(\alpha)-1)$-horizon-first-step optimal for some discount factors close to $\alpha.$ Proposition~\ref{ndclass} shows that a point, at which the turnpike function is discontinuous, must be either an irregular point or an $n$-horizon irregular point for some $n\in\N^+$. As a result, any subinterval of $[0,1)$, on which the turnpike function is bounded, can be partitioned into finitely many turnpike intervals; Theorem~\ref{finitetpi}. 
Proposition~\ref{Lebesgue} states the existence of a closed subset of the interval $[0,1)$  on which the turnpike function is bounded, and the Lebesgue measures of the points in $[0,1)$ outside such a set can be chosen arbitrarily small. 

Let $\norm{\cdot}$ be the norms in $\R^m$ or $\R^{m\times m}$ defined as
\[\norm{V}:=\max_{x\in\X}\left|V(x)\right|,\ V\in\R^m;\qquad \norm{B}:=\max_{1\le i\le m}\sum_{j=1}^m\left|B_{i,j}\right|,\ B\in\R^{m\times m}.\]
	Note that $P\pmb{1}=\pmb{1}$ and $\norm{P}=1$ for any stochastic matrix $P$, and $\norm{Bw}\le\norm{B}\norm{w}$ for $B\in\R^{m\times m},w\in\R^m$. Let us define the maximum absolute values
    
	\begin{equation}\label{eqR}
		R_1:= \max_{x\in\X,a\in A(x)}\abs*{r(x,a)},\qquad R_2:=\max_{x\in\X}\abs*{s(x)},\qquad R:=\max\left\{R_1,R_2\right\}.
	\end{equation}
We begin with the following lemma which provides bounds of objective functions and value functions.

\begin{lem}\label{lembound} (a) $\max\{\norm{V_\alpha},\norm{V_{n,\alpha}},\norm{v_\alpha^\pi},\norm{v_{n,\alpha}^\pi}\}\le\frac{R}{1-\alpha}$ for all $\alpha\in[0,1),n\in\N,\pi\in\Pi;$\smallskip\\
(b) $\norm{V_\alpha-V_{n,\alpha}}\le\frac{(2-\alpha)\alpha^nR}{1-\alpha}$ for all $n\in\N,\alpha\in[0,1)$.
\end{lem}
\begin{proof}
(a). This statement follows from taking the norm on both sides of \eqref{eqvaluesum1} and \eqref{eqvaluesum3}, from the triangle inequality,  \eqref{eqR}, and because $P_t(\pi)$ are stochastic matrices. 
    
    (b). It is well-known that the optimality operator $T_\alpha$ is a contraction mapping with contracting coefficient $\alpha$ \cite[Proposition 6.2.4, p.~163]{Puterman:1994:MDP:528623}. Then
    \begin{equation*}
    \begin{aligned}
        \norm{V_\alpha-V_{n,\alpha}}&=\norm{T_\alpha V_\alpha-T_\alpha V_{n-1,\alpha}}\le\alpha\norm{V_\alpha-V_{n-1,\alpha}}\le\alpha^2\norm{V_\alpha-V_{n-2,\alpha}}
        \le\cdots\le\alpha^n\norm{V_\alpha-V_{0,\alpha}}\\
        &=\alpha^n\norm{V_\alpha-s}\le\alpha^n(\norm{V_\alpha}+\norm{s})\le\alpha^n\left(\frac{R}{1-\alpha}+R\right)=\frac{(2-\alpha)\alpha^nR}{1-\alpha},
    \end{aligned}
    \end{equation*}
    where the last inequality is by \eqref{eqR} and by (a).
\end{proof}

The following lemma shows that, the value functions $V_\alpha$ and $V_{n,\alpha},$ $n\in\N,$ are  Lipschitz continuous. 

\begin{lem}\label{lemlipschitz} Let $0\le a<b<1$. Then $\norm{V_{n,\alpha_1}-V_{n,\alpha_2}}\le\frac{R}{(1-b)^2}\abs{\alpha_1-\alpha_2}$ and $\norm{V_{\alpha_1}-V_{\alpha_2}}\le\frac{R}{(1-b)^2}\abs{\alpha_1-\alpha_2}$ for all $n\in\N$ and for all $\alpha_1,\alpha_2\in[a,b]$.
\end{lem}
\begin{proof} Let $n\in\N$ be fixed. By Proposition~\ref{partition}(a), on each $n$-horizon partition interval $\mathcal{I}_i^{(n)},$ where partition intervals are defined in Definition~\ref{defpartition}, for some policy $\pi\in\Pi,$ we have $V_{n,\alpha}=v_{n,\alpha}^{\pi}$ for all $\alpha\in \mathcal{I}_i^{(n)}.$ By formula \eqref{eqvaluesum1},  for $x\in\X$ and $\alpha\in [0,b],$
\begin{equation*}
\begin{aligned}
	\abs*{\frac{d}{d\alpha}V_{n,\alpha}(x)} & = \abs*{\frac{d}{d\alpha}v_{n,\alpha}^\pi(x)}  = \abs*{\sum_{t=0}^{n-1}t\alpha^{t-1}P_t(\pi)r_t(\pi)(x)+n\alpha^{n-1}P_n(\pi)s(x)}\\ &\le R\sum_{t=0}^{\infty}t\alpha^{t-1} = \frac{R}{(1-\alpha)^2}\le \frac{R}{(1-b)^2},
\end{aligned}
\end{equation*}
and this inequality implies the first inequality stated in the lemma since, in view of Proposition~\ref{partition}(a), the set of $n$-horizon irregular points is finite. The second inequality follows by taking $n\to\infty$ and applying Lemma~\ref{lembound}(b).
\end{proof}

\begin{thm}\hspace*{-0.2cm}{\rm (cp. \cite[Theorem 1]{doi:10.1287/moor.2017.0912})}\label{utt} Let $\mathcal{I}:=[a,b]\subset [0,1)$ and all $\alpha\in\mathcal{I}$ are regular points, that is, $\P(\mathcal{I})=\emptyset$. Then there exists $K<\infty$ such that $D_n(\mathcal{I})\subset D(\mathcal{I})=D(\alpha)$ for all $n\ge K$ and for all $\alpha\in\mathcal{I},$ and therefore $N^*(\mathcal{I})\le K.$
\end{thm}

The difference between Theorem~\ref{utt} and \cite[Theorem 1]{doi:10.1287/moor.2017.0912}  is that Theorem~\ref{utt} does not assume that the terminal reward $s$ identically equals to $\mathbf{0}.$ The proof of Theorem~\ref{utt}  provided below is significantly shorter than the proof of Theorem 1 in \cite{doi:10.1287/moor.2017.0912}.

\begin{proof}{Proof of Theorem~\ref{utt}.} By Proposition~\ref{partition}(c), $D(\mathcal{I})=D(\alpha)$ for all $\alpha\in\mathcal{I}.$ Suppose that Theorem~\ref{utt} is not true. Then there exist an increasing sequence of positive integers $\{n_i\}_{i=1}^\infty\subset\N^+$ and a sequence of discount factors $\{\alpha_i\}_{i=1}^\infty\subset\mathcal{I}$ such that $D_{n_i}(\alpha_i)\not\subset D(\mathcal{I})$ for all $i\ge1$. This means for each $i\in\mathbb{N^+}$, there exists $\phi^{(i)}\in\F\setminus D(\mathcal{I})$ such that $V_{n_i,\alpha_i}=T_{\alpha_i}^{\phi^{(i)}} V_{n_i-1,\alpha_i}=r(\phi^{(i)})+\alpha_iP(\phi^{(i)})V_{n_i-1,\alpha_i}.$ Since $\F$ is finite and $\mathcal{I}$ is 
closed, there exist $\phi\in\F\setminus D(\mathcal{I})$, $\alpha'\in\mathcal{I}$, and an increasing sequence $\{i_j\}_{j=1}^\infty$ of positive integers such that  $\alpha_{i_j}\to\alpha'$ as $j\to\infty$ and
\begin{equation}\label{nij}
	V_{n_{i_j},\alpha_{i_j}}=r(\phi)+\alpha_{i_j}P(\phi)V_{n_{i_j-1},\alpha_{i_j}}.
\end{equation}
We claim that $\lim_{j\to\infty}V_{n_{i_j},\alpha_{i_j}}=V_{\alpha'}$ and $\lim_{j\to\infty}V_{n_{i_j}-1,\alpha_{i_j}}=V_{\alpha'}$. The first limit follows from
\begin{equation}\label{limit}
	\norm{V_{\alpha'}-V_{n_{i_j},\alpha_{i_j}}}\le\norm{V_{\alpha'}-V_{n_{i_j},\alpha'}}+\norm{V_{n_{i_j},\alpha'}-V_{n_{i_j},\alpha_{i_j}}}\le\frac{(2-\alpha')(\alpha')^{n_{i_j}}R}{1-\alpha'}+\frac{R\left|\alpha'-\alpha_{i_j}\right|}{(1-b)^2}\to 0
\end{equation}
as $j\to\infty$, where the first inequality follows from the triangle inequality, and the second inequality follows from Lemma~\ref{lembound}(b) and Lemma~\ref{lemlipschitz}. The second limit follows similarly by replacing the $n_{i_j}$ in \eqref{limit} with $n_{i_j}-1.$   Thus, by taking $j\to\infty$, \eqref{nij} becomes 
\begin{equation*}
	V_{\alpha'}=r(\phi)+\alpha'P(\phi)V_{\alpha'}=T_{\alpha'}^\phi V_{\alpha'},
\end{equation*}
which by \eqref{eqoptimal} implies $\phi\in D(\mathcal{I})$. This is a contradiction.
\end{proof}

It is possible that the turnpike function $N(\alpha)$ is unbounded near an irregular point $\alpha$ or near $\alpha=1;$ see \cite[Example 1]{doi:10.1287/moor.2017.0912} for $N(\alpha)$ being unbounded only at one side of a non-touching break point, and \cite[Example 3]{doi:10.1287/moor.2017.0912} for $N(\alpha)$ being unbounded near a touching point. In Example~\ref{exatbbounded},  $N(\alpha)$ is unbounded on both sides of a break point. In Example~\ref{exaub1},  $N(\alpha)$ is unbounded near $\alpha=1.$ In Example~\ref{exatbb},  $N(\alpha)$ is bounded near a break point.

Statement (a) of the following proposition can be viewed as an extension of Theorem~\ref{utt} to intervals containing irregular points. This statement shows that for a closed subinterval of $[0,1)$ that contains at least one irregular point, and for a finite-horizon problem with a sufficiently large horizon, every first-step-optimal policy for a discount factor from this subinterval must also be optimal for an infinite-horizon problem at one of the irregular points of this interval. Statement (b) shows that for each $\mathcal{I}\in\I$ the turnpike function is bounded on some closed subsets of $\mathcal{I}$ whose Lebesgue measures can be arbitrarily close to the Lebesgue measure of $\mathcal{I}$.

\begin{pro}\label{uttirr} Let $\mathcal{I}:=[a,b]\subset [0,1)$ such that $\mathcal{I}$ contains at least one irregular point, that is, $\P(\mathcal{I})\ne\emptyset.$ Then 
\begin{enumerate}[(a)]
	\item there exists $K<\infty$ such that $D_n(\mathcal{I})\subset D(\P(\mathcal{I}))$ for all $n\ge K$;
	\item for any $\epsilon>0$, there exists a closed set $E\subset\mathcal{I}$ such that $N^*(E)<\infty,$ and $\mu(\mathcal{I}\setminus E)<\epsilon,$ where $\mu$ is the Lebesgue measure on $[0,1].$
   \end{enumerate}
\end{pro}
\begin{proof}
 (a). Let $\P(\mathcal{I})=\{c_1,\ldots,c_k\}\ne\emptyset$, where $c_1<\ldots<c_k$. For each $1\le j\le k$, we choose  $\epsilon_j>0$ such that $c_k+\epsilon_k<1$ and $\{[c_j-\epsilon_j,c_j+\epsilon_j]\}_{j=1}^k$ are disjoint intervals. Then for each $1\le j\le k$, the interval $[c_j-\epsilon_j,c_j+\epsilon_j]$ contains exactly one irregular point $c_j$. We claim that, for each $1\le j\le k$, there exists some $K_j\in\N^+$ such that $D_n(\alpha)\subset D(c_j)$ for all $\alpha\in[c_j-\epsilon_j,c_i+\epsilon_j]\cap[0,1)$. Otherwise, let us fix some $j=1,\ldots,k$ for which this is not true. Then closedness of $[c_j-\epsilon_j,c_j+\epsilon_j]\cap[0,1)$ and the finiteness of $\F$ imply that there exist a sequence of positive integers $\{n_i\}_{i=1}^\infty\subset\N^+$ such that $n_i\to\infty$ as $i\to\infty$, a sequence of discount factors $\{\alpha_i\}_{i=1}^\infty\subset[c_j-\epsilon_j,c_j+\epsilon_j]\cap[0,1)$ such that $\alpha_i\to\alpha'$ for some $\alpha'\in[c_j-\epsilon_j,c_j+\epsilon_j]\cap[0,1)$, and a decision rule $\phi\in\F$ such that $V_{n_i,\alpha_i}=T_{\alpha_i}^\phi V_{n_i-1,\alpha_i}$ for all $i\ge1$ but $\phi\notin D(c_j).$ By taking $i\to\infty,$ we have $V_{\alpha'}=T_{\alpha'}^\phi V_{\alpha'}$ in view of~\eqref{limit}, which implies $\phi\in D(\alpha')$ by~\eqref{eqoptimal}. However, Proposition~\ref{partition}(c) and Proposition~\ref{upperhemi} imply that $D(\alpha')\subset D(c_j),$ which is a contradiction. Hence, the claim is true, and we set $\tilde{K}:=\max_{1\le j\le k}K_j$. Next, we observe that
\begin{equation}\label{irrset}
    \mathcal{I}\setminus\left(\cup_{j=1}^k(c_j-\epsilon_j,c_j+\epsilon_j)\right)=[a,c_1-\epsilon_1]\cup\left(\cup_{j=1}^{k-1}[c_j+\epsilon_j,c_{j+1}-\epsilon_{j+1}]\right)\cup[c_k+\epsilon_k,b],
\end{equation}
where $[a,c_1-\epsilon_1]:=\emptyset$ for $a>c_1-\epsilon_1$ and $[c_k+\epsilon_k,b]:=\emptyset$ for $c_k+\epsilon_k>b.$ Note that each of these finitely many closed sets on the right-hand side of \eqref{irrset} is a subset of $[0,1)$ and does not contain irregular points. 
Thus, by Theorem~\ref{utt}, there exists $K'\in\N^+$ such that $D_n(\alpha)\subset D(\alpha)\subset D(\P(\mathcal{I}))$ for all $\alpha\in\mathcal{I}\setminus\left(\cup_{j=1}^k(c_j-\epsilon_j,c_j+\epsilon_j)\right)$ and for all $n\ge K'$. Therefore, (a) is proved with $K=\max\{\tilde{K},K'\}$.

(b). In the proof of part (a), for each $i=1,\ldots,k$   let us choose  $\epsilon_i\le\epsilon/(2k).$ Then set $E=\mathcal{I}\setminus\left(\cup_{j=1}^k(c_j-\epsilon_j,c_j+\epsilon_j)\right).$
\end{proof}

Proposition~\ref{uttirr} is strengthened later in this paper. By Corollary~\ref{uttirrno0}, the set $D(\P(\mathcal{I}))$ can be replaced with a smaller or equal set $D(\P(\mathcal{I})\setminus\{0\})$ in Proposition~\ref{uttirr}(a). By Proposition~\ref{Lebesgue}, the closed set $E$ in Proposition~\ref{uttirr}(b) can be chosen as a finite union of disjoint closed turnpike intervals.  

\begin{cor}\label{absorb} For each $\alpha\in[0,1)$ there exist $\epsilon>0$ and  $K\in\N^+$ such that $D_n((\alpha-\epsilon,\alpha+\epsilon)\cap[0,1))\subset D(\alpha)$  for all $n\ge K.$
\end{cor}
\begin{proof}
    Recall from Proposition~\ref{partition}(c) that all irregular points are isolated. If $\alpha$ is a regular point, then this corollary follows from Theorem~\ref{utt}. Otherwise, it follows from Proposition~\ref{uttirr}(a).
\end{proof}

\begin{thm}\label{uppersemi} The turnpike function $N(\alpha)$ is upper semicontinuous at each regular point $\alpha\in [0,1)$, and therefore it is upper semicontinuous on each partition interval.
\end{thm}
\begin{proof} Suppose there exists some regular point $\alpha\in[0,1)$ such that $N(\alpha)$ is not upper semicontinuous at $\alpha.$ By Proposition~\ref{partition}(c) there exist $a<\alpha$ and $b\in(\alpha,1)$ such that all $\alpha\in\mathcal{I}:=[a,b]\cap[0,1)$ are regular points. Theorem~\ref{utt} implies that there exists $K<\infty$ such that $N^*(\mathcal{I})\le K.$ Since $N(\alpha)$  takes integer values, there exist an increasing sequence of positive integers $\{n_i\}_{i=1}^\infty\subset\N^+$ and a sequence of discount factors $\{\alpha_i\}_{i=1}^\infty\subset\mathcal{I}$ such that $\alpha_i\rightarrow\alpha$ as $i\to\infty$, and $N(\alpha)+1\le N(\alpha_i)\le K$ for all $i\ge1$. Then boundedness of $\{N(\alpha_i)\}_{i=1}^\infty$ implies that there exists a subsequence $N(\alpha_{i_j})$ consisting of constant numbers equal to some $n\in\N^+$ such that $N(\alpha)+1\le n\le K$. So,   $N(\alpha_{i_j})=n$ and $N(\alpha)\le n-1$. Then by the definition of a turnpike integer, for each $j\in\N^+$ there exists $\phi^{(j)}\in D_{n-1}(\alpha_{i_j})\setminus D(\alpha_{i_j})$. By Theorem~\ref{utt}, $D(\alpha)=D(\mathcal{I})$ for all $\alpha\in\mathcal{I}$, and thus $\phi^{(j)}\in D_{n-1}(\alpha_{i_j})\setminus D(\alpha)$. By Proposition~\ref{upperhemi}, $\phi^{(j)}\in D_{n-1}(\alpha)\setminus D(\alpha)$ when $j$ is large enough, which implies $N(\alpha)\ge n.$ This is a contradiction.
\end{proof}

\begin{cor}\label{0regular} If 0 is a regular point, then there exists $\delta\in(0,1)$ such that $N(\alpha)=1$ for all $\alpha\in[0,\delta).$
\end{cor}
\begin{proof}
By Theorem~\ref{uppersemi} $N(\alpha)$ is upper semicontinuous at $0$. Since $N(\alpha)$ is integer-valued, there exists $\delta\in(0,1)$ such that $N(\alpha)\le N(0)=1$ for all $\alpha\in[0,\delta).$ Since $N(\alpha)\ge1$, we have $N(\alpha)=1$ for all $\alpha\in[0,\delta).$
\end{proof}

A concrete value of the $\delta$ in Corollary~\ref{0regular} is provided by Corollary~\ref{0regularb} and by \eqref{FXdef}, \eqref{Cdef}, \eqref{Deldef}. 

\begin{cor}\label{Ndge2} Let $\mathcal{I}$ be a partition interval. If $\alpha\in \mathbb{D}(\mathcal{I}),$ then $N(\alpha)\ge2.$
\end{cor}
\begin{proof}
    Suppose this is not true, that is, $N(\alpha)=1$ where $\alpha\in \mathbb{D}(\mathcal{I}).$ By Theorem~\ref{uppersemi} there exists $\delta>0$ such that $N(\beta)=1$ for all $\beta\in(\alpha-\delta,\alpha+\delta),$ which implies continuity of $N(\alpha)$  at $\alpha$. This contradicts the assumption $\alpha\in \mathbb{D}(\mathcal{I}).$
\end{proof}

The following proposition shows that the turnpike function $N(\alpha)$ is discontinuous at some regular point $\alpha$ if and only if there exists $\delta >0$ such that, for the discount factor $\alpha$, if a nonoptimal decision rule is first-step-optimal for the $(N(\alpha)-1)$-horizon problem, then, either for all discount factors $\beta\in (\alpha-\delta,\alpha)$ or for all discount factors $\beta\in (\alpha,\alpha+\delta),$  it is not first-step-optimal for the $(N(\alpha)-1)$-horizon problem. Furthermore, the former case corresponds to $N(\alpha)$ not being left-continuous at $\alpha$, and the latter case corresponds to $N(\alpha)$ not being right-continuous at $\alpha$.

\begin{pro}\label{Ndiff} Let $\mathcal{I}$ be a partition interval. Then
\begin{enumerate}[(a)]
	\item $\alpha\in \mathbb{D}^-(\mathcal{I})$ iff there exists $\delta >0$ such that for all $\beta\in (\alpha-\delta,\alpha),$ if $\psi\in D_{N(\alpha)-1}(\alpha)\setminus D(\mathcal{I}),$ then $\psi\notin D_{N(\alpha)-1}(\beta);$ 
	\item $\alpha\in \mathbb{D}^+(\mathcal{I})$ iff there exists $\delta >0$ such that for all $\alpha\in (\alpha,\alpha+\delta),$ if $\psi\in D_{N(\alpha)-1}(\alpha)\setminus D(\mathcal{I}),$ then $\psi\notin D_{N(\alpha)-1}(\beta).$ 
\end{enumerate}
\end{pro}
\begin{proof} Let $\alpha\in \mathbb{D}(\mathcal{I}).$ Then $N(\alpha)\ge2$ by Corollary~\ref{Ndge2}. For simplicity, we only prove (a), as the proof of (b) is similar.

Let us prove the necessity. Let $\alpha\in \mathbb{D}^-(\mathcal{I})$. Recall that $D_{N(\alpha)-1}(\alpha)\setminus D(\mathcal{I})\ne\emptyset$ by the definition of turnpike integers. Suppose the necessity is not true, that is, there exist $\delta_1>0$ and $\psi\in D_{N(\alpha)-1}(\alpha)\setminus D(\mathcal{I})$ such that $\psi\in D_{N(\alpha)-1}(\beta)$ for all $\beta\in(\alpha-\delta_1,\alpha),$ where the existence of such $\delta$ is by Proposition~\ref{partition}(b). This implies $N(\beta)\ge N(\alpha)$ for all $\beta\in(\alpha-\delta_1,\alpha).$ By Theorem~\ref{uppersemi} there exists $\delta_2>0$ such that $N(\beta)\le N(\alpha)$ for all $\beta\in(\alpha-\delta_2,\alpha).$ Let $\delta=\min\{\delta_1,\delta_2\}.$ Then $N(\beta)=N(\alpha)$ for all $\beta\in(\alpha-\delta,\alpha),$ which contradicts the assumption $\alpha\in \mathbb{D}^-(\mathcal{I})$. Therefore, the necessity is true.

Let us prove the sufficiency. Let $\delta_1>0$ such that for all $\beta\in (\alpha-\delta_1,\alpha),$ if $\psi\in D_{N(\alpha)-1}(\alpha)\setminus D(\mathcal{I}),$ then $\psi\notin D_{N(\alpha)-1}(\beta).$ By Proposition~\ref{upperhemi} there exists $\delta_2>0$ such that $D_{N(\alpha)-1}(\beta)\subset D_{N(\alpha)-1}(\alpha)$ for all $\beta\in(\alpha-\delta_2,\alpha).$ By Theorem~\ref{uppersemi} there exists $\delta_3>0$ such that $N(\beta)\le N(\alpha)$ for all $\beta\in(\alpha-\delta_3,\alpha)$, which also implies $D_{n}(\beta)\subset D(\mathcal{I})$ for all $\beta\in(\alpha-\delta_3,\alpha)$ if $n\ge N(\alpha)$. Let $\delta=\min\{\delta_1,\delta_2,\delta_3\}.$ For each $\beta\in(\alpha-\delta,\alpha)$, we claim that $D_{N(\alpha)-1}(\beta)\setminus D(\mathcal{I})=\emptyset.$ Otherwise, if $\psi\in D_{N(\alpha)-1}(\beta)\setminus D(\mathcal{I})$, then $\psi\in D_{N(\alpha)-1}(\alpha)\setminus D(\mathcal{I})$. Then $\psi\notin D_{N(\alpha)-1}(\beta)$ by our assumption, which is a contradiction. Thus, the claim is true, that is, $D_{N(\alpha)-1}(\beta)\subset D(\mathcal{I})$ for all $\beta\in(\alpha-\delta,\alpha)$. Hence, $D_{n}(\beta)\subset D(\mathcal{I})$ for all $\beta\in(\alpha-\delta,\alpha)$ and for all $n\ge N(\alpha)-1$, which implies $N(\beta)\le N(\alpha)-1$ for all $\beta\in(\alpha-\delta,\alpha)$. Therefore, $\alpha\in \mathbb{D}^-(\mathcal{I})$.
\end{proof}

For a regular point $\alpha,$ at which the turnpike function $N(\alpha)$ is not continuous, statements (b) and (c) of the following proposition show that there exist $\delta>0$ and an optimal decision rule $\phi$ such that $\phi$ is also $(N(\alpha)-1)$-horizon-first-step optimal for all discount factors in either $[\alpha-\delta,\alpha]$ or $[\alpha,\alpha+\delta].$

\begin{pro}\label{Nd} Let $\mathcal{I}$ be a partition interval, and $\alpha\in \mathbb{D}(\mathcal{I})$. Then
\begin{enumerate}[(a)]
    \item $D_{N(\alpha)-1}(\alpha)$ contains at least two decision rules: one from $D(\mathcal{I})$ and another one from $\F\setminus D(\mathcal{I});$
    \item if $\alpha\in \mathbb{D}^-(\mathcal{I})$, then there exist $\delta >0$ and $\phi\in D(\mathcal{I})$ such that $\phi\in D_{N(\alpha)-1}(\beta)$ for all $\beta\in[\alpha-\delta,\alpha]$;
    \item if $\alpha\in \mathbb{D}^+(\mathcal{I})$, then there exist $\delta >0$ and $\phi\in D(\mathcal{I})$ such that $\phi\in D_{N(\alpha)-1}(\beta)$ for all $\beta\in[\alpha,\alpha+\delta]$;
\end{enumerate}
\end{pro}
\begin{proof} Recall that $N(\alpha)\ge2$ by Corollary~\ref{Ndge2}. Let us prove (b). Let $\alpha\in \mathbb{D}^-(\mathcal{I}).$ By Proposition~\ref{upperhemi} there exists $\delta_1>0$ such that $D_{N(\alpha)-1}(\beta)\subset D_{N(\alpha)-1}(\alpha)$ for all $\beta\in(\alpha-\delta_1,\alpha).$ Then by Proposition~\ref{Ndiff}(a) there exists $\delta_2\in(0,\delta_1)$ such that $D_{N(\alpha)-1}(\beta)\subset D_{N(\alpha)-1}(\alpha)\cap D(\mathcal{I})$ for all $\beta\in (\alpha-\delta_2,\alpha).$ Also, by Proposition~\ref{partition}(b), there exist $\delta_3\in(0,\delta_2)$ and $\phi\in D_{N(\alpha)-1}(\alpha)\cap D(\mathcal{I})$ such that $\phi\in D_{N(\alpha)-1}(\beta)$ for all $\beta\in[\alpha-\delta_3,\alpha].$ Therefore, (b) is proved. The proof of (c) is similar to the proof of (b), and (a) follows from  (b) and (c). 
\end{proof}

Note that the converse of Proposition~\ref{Nd}(b) or \ref{Nd}(c) may not be true, which does not conflict with Proposition~\ref{Ndiff}. For example, for $\alpha\in \mathbb{D}(\mathcal{I}),$ where $\mathcal{I}$ is a partition interval, there may exist $\delta>0,$ $\phi\in D(\mathcal{I})$ and $\psi\in D_{N(\alpha)-1}(\alpha)\setminus D(\mathcal{I})$ such that $\{\phi,\psi\}\subset D_{N(\alpha)-1}(\beta)$ for all $\beta\in[\alpha,\alpha+\delta]$, and therefore the necessary condition of Proposition~\ref{Nd}(c) is satisfied, but $\alpha\notin \mathbb{D}^+(\mathcal{I})$ in view of Proposition~\ref{Ndiff}(b); see Example~\ref{exaDn+-}.

The following proposition shows that, at a regular point $\alpha$, if $N(\alpha)$ is neither left-continuous nor right-continuous, then $\alpha$ must be an $(N(\alpha)-1)$-horizon-first-step touching point; if $N(\alpha)$ is only continuous on one side of $\alpha$, then $\alpha$ must be an $(N(\alpha)-1)$-horizon-first-step break point. In particular, if $\alpha$ is a regular point at which $N(\alpha)$ is discontinuous, then it must be an $(N(\alpha)-1)$-horizon-first-step irregular point.

\begin{pro} \label{ndclass}
Let $\mathcal{I}$ be a partition interval. Then
\begin{enumerate}[(a)]
	\item if $\alpha\in\Hat{\mathbb{D}}(\mathcal{I}),$ then $\alpha\in\P^T_{N(\alpha)-1}(\mathcal{I});$
	\item if $\alpha\in \mathbb{D}(\mathcal{I})\setminus\Hat{\mathbb{D}}(\mathcal{I}),$ then $\alpha\in\P^B_{N(\alpha)-1}(\mathcal{I}).$
\end{enumerate}
In particular, if $\alpha\in \mathbb{D}(\mathcal{I}),$ then $\alpha\in\P_{N(\alpha)-1}(\mathcal{I}).$ Therefore, the set $\mathbb{D}(\mathcal{I})$ is at most countable, and so are $\mathbb{D}^-(\mathcal{I}),$ $\mathbb{D}^+(\mathcal{I}),$ and $\Hat{\mathbb{D}}(\mathcal{I}).$
\end{pro}
\begin{proof}
 (a). Let $\alpha\in \Hat{\mathbb{D}}(\mathcal{I}).$ Recall that $N(\alpha)\ge2$ by Corollary~\ref{Ndge2}. Let $\psi\in D_{N(\alpha)-1}(\alpha)\setminus D(\mathcal{I})$ which is not empty by the definition of $N(\alpha)$. Since $\alpha\in \mathbb{D}^-(\mathcal{I})\cap \mathbb{D}^+(\mathcal{I})$, by Proposition~\ref{Ndiff} there exists $\delta>0$ such that $\psi\notin D_{N(\alpha)-1}(\beta)$ for all $\beta\in(\alpha-\delta,\alpha)\cup(\alpha,\alpha+\delta)$. In view of Definitions~\ref{def+-} and \ref{defirrfinite}, this implies $D_{N(\alpha)-1}(\alpha)\ne D_{N(\alpha)-1}(\alpha-)\cup D_{N(\alpha)-1}(\alpha+)$, which is $\alpha\in\P^T_{N(\alpha)-1}(\mathcal{I})$. 

 (b). Let $\alpha\in \mathbb{D}^-(\mathcal{I})\setminus \Hat{\mathbb{D}}(\mathcal{I})$. Then $N(\alpha)\ge2$ by Corollary~\ref{Ndge2}, and $\alpha\notin \mathbb{D}^+(\mathcal{I}).$ Then by Proposition~\ref{partition}(b) and by Proposition~\ref{Ndiff}(b) there exist $\delta_1>0$ and $\psi\in D_{N(\alpha)-1}(\alpha)\setminus D(\mathcal{I})$ such that $\psi\in D_{N(\alpha)-1}(\beta)$ for all $\beta\in(\alpha,\alpha+\delta_1)$. Since $\alpha\in \mathbb{D}^-(\mathcal{I})$, by Proposition~\ref{Ndiff} there exists $\delta_2>0$ such that $\psi\notin D_{N(\alpha)-1}(\beta)$ for all $\beta\in (\alpha-\delta_2,\alpha).$ In view of Definition~\ref{def+-}, this implies $\psi\in D_{N(\alpha)-1}(\alpha+)$ but $\psi\notin D_{N(\alpha)-1}(\alpha-)$, and thus $D_{N(\alpha)-1}(\alpha-)\ne D_{N(\alpha)-1}(\alpha+).$ Therefore, $\alpha\in\P^B_{N(\alpha)-1}(\mathcal{I})$ by Definition~\ref{defirrfinite}. The case $\alpha\in \mathbb{D}^+(\mathcal{I})\setminus \Hat{\mathbb{D}}(\mathcal{I})$ is similar.
\end{proof}

Recall that if the function $N(\alpha)$ is discontinuous at $\alpha\in [0,1),$ then $\alpha\in \mathbb{D}(\mathcal{I}):=\mathbb{D}^-(\mathcal{I})\cup \mathbb{D}^+(\mathcal{I}).$ If $\mathcal{I}$ is a partition interval, then the sets $\mathbb{D}^-(\mathcal{I})$ and $\mathbb{D}^+(\mathcal{I})$ can be infinite; \cite[Example 3]{doi:10.1287/moor.2017.0912}. However, in view of Proposition~\ref{ndclass}, these sets are at most countable. For a partition interval $\mathcal{I}$, Example~\ref{exaisolate} shows that $\Hat{\mathbb{D}}(\mathcal{I})\ne\emptyset$ is possible, and Example~\ref{exalirb} shows that $\mathbb{D}^-(\mathcal{I})\ne\emptyset,$ $\mathbb{D}^+(\mathcal{I})\ne\emptyset$ are both possible.


According to Definition~\ref{defdiscon}, if $\mathcal{I}\in\I,$ then $\P(\mathcal{I})$ is the set of irregular points from $\mathcal{I},$ which is finite by Proposition~\ref{partition}(c). For $\mathcal{I}\in\I$, define
\begin{equation*}
    \overline{\P}(\mathcal{I}):=\P(\mathcal{I})\bigcup\left(\bigcup_{1\le n< N^*(\mathcal{I})}\P_n(\mathcal{I})\right).
\end{equation*}
In other words, $\overline{\P}(\mathcal{I})$ is the set of points from $\mathcal{I}$ that are either irregular or first-step-irregular for some finite-horizon problem with the horizon smaller than $N^*(\mathcal{I}).$

\begin{cor}\label{Ndinclu} Let $\mathcal{I}\in\I.$ Then $\mathbb{D}(\mathcal{I})\subset \overline{\P}(\mathcal{I}).$
\end{cor}
\begin{proof}
    Let $\alpha\in \mathbb{D}(\mathcal{I}).$ If $\alpha\notin\P(\mathcal{I})$, that is, $\alpha$ is a regular point, then $\alpha\in\P_{N(\alpha)-1}(\mathcal{I})$ by Proposition~\ref{ndclass}.
\end{proof}

\cite[Theorem 2]{doi:10.1287/moor.2017.0912} states that a closed subinterval of $[0,1)$ without irregular points can be partitioned into finitely many turnpike intervals, where turnpike intervals are defined in Definition~\ref{deftpi}.  The following theorem generalizes this fact in view of Theorem~\ref{utt}.

\begin{thm}\label{finitetpi} Let $\mathcal{I}\in\I$. Then $\mathcal{I}$ can be partitioned into finitely many turnpike intervals iff $N^*(\mathcal{I})<\infty$.
\end{thm}
\begin{proof} Let us prove the necessity. If $\mathcal{I}=\cup_{i=1}^j\mathcal{I}_i,$ where each $\mathcal{I}_i$ is a turnpike interval, then $N^*(\mathcal{I})=\max_{1\le i\le k}N^*(\mathcal{I}_i)<\infty,$ where each $N^*(\mathcal{I}_i)$ is finite in view of the definition of turnpike interval. Let us prove the sufficiency. Suppose $N^*(\mathcal{I})<\infty$. Then $\mathcal{I}$ can be partitioned into finitely many turnpike intervals is equivalent to $\mathbb{D}(\mathcal{I})$ is finite, which is implied by Corollary~\ref{Ndinclu} and the finiteness of $\P_n(\mathcal{\I})$ for each $n\in\N^+$ by Proposition~\ref{partition}(b).
\end{proof}

\begin{figure}[ht!]
	\centerline{\includegraphics[scale=0.8]{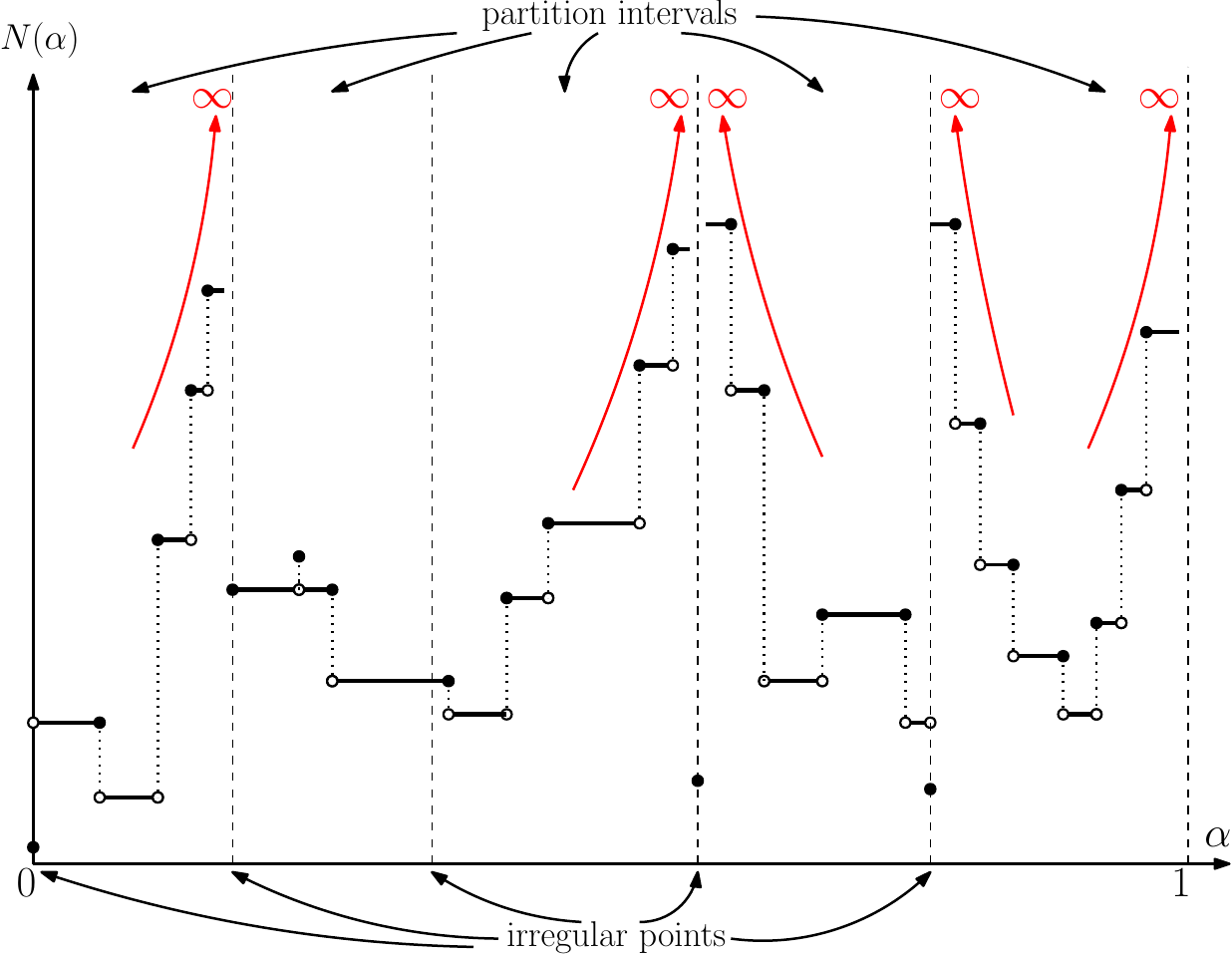}}
	\caption{Graph of a turnpike function consistent with Theorems~\ref{utt},~\ref{uppersemi}, and~\ref{finitetpi}. The red arrows indicate that the turnpike function goes to infinity when the discount factor $\alpha$ goes to the corresponding sides of irregular points.}
    \label{fig:NB}
\end{figure}

Figure~\ref{fig:NB} illustrates the graph of a turnpike function. Example~\ref{exaisolate} shows that the turnpike function can be neither left-continuous nor right-continuous at a regular point. At an irregular point, Example~\ref{exalirb} shows that the turnpike function can be unbounded on only one side, either left or right, and Example~\ref{exatbbounded} shows that the turnpike function can be unbounded on both sides. Example~\ref{exaub1} shows that the turnpike function can be unbounded near $\alpha=1,$ while all other examples in Section~\ref{secexamples} have bounded turnpike functions near $\alpha=1.$ Example~\ref{exatbb} shows that the turnpike function can be bounded on both sides of an irregular point. Example~\ref{exanbm} shows that the turnpike function can be discontinuous at $\alpha=0,$ i.e., it is a touching point in this example.

\begin{rem}\label{remwrong}\emph{
Alternatively, Theorem~\ref{finitetpi} can be derived by only using Proposition~\ref{partition}. Indeed, let $0\le a<b\le 1,$ where $a,b$ are the end points of $\mathcal{I}$. Since $N^*(\mathcal{I})<\infty,$ Proposition~\ref{partition}(b,c) imply that $\overline{\P}(\mathcal{I})$ is finite. Let us order the points in $\overline{\P}(\mathcal{I})$ in the increasing order $c_1,c_2,\ldots,c_j$, and set $c_0:=a$, $c_{j+1}:=b$. Then $a=c_0<c_1<\ldots<c_j<c_{j+1}=b$. For each $i=0,1,\ldots,j$, the interval $(c_{i},c_{i+1})$ is a subinterval of some partition interval, as well as a subinterval of some $n$-horizon-first-step partition interval for each $n=1,2,\ldots, N^*(\mathcal{I})-1$, where partition intervals are defined in Definition~\ref{defpartition}. Again, by Proposition~\ref{partition}(b,c), this implies that for each $n=1,2,\ldots,N^*(\mathcal{I})-1$ and each $i=0,1,\ldots,j,$ the set-valued functions $D_n(\alpha)$ and $D(\alpha)$ are both constants in $(c_i,c_{i+1}),$  from which we conclude that $(c_i,c_{i+1})$ is a turnpike interval.}
\end{rem}

The following corollary shows that, if $\alpha$ is a limit point of a set of discount factors where the turnpike function is discontinuous, then either $\alpha$ is an irregular point or $\alpha=1,$ and the latter case is possible in view of Example~\ref{exaub1}.

\begin{cor}\label{ndlimit}
If $\alpha$ is a limit point of $\mathbb{D}([0,1))$, then either $\alpha\in\P([0,1))$ or $\alpha=1.$
\end{cor}
\begin{proof}
    This statement follows from Proposition~\ref{partition}(c), Theorem~\ref{utt}, and Theorem~\ref{finitetpi}.
\end{proof}

Corollary~\ref{ndlimit} is strengthened in Corollary~\ref{ndlimitg}, which excludes the possibility that $0$ is a limit point of $\mathbb{D}([0,1))$. The following theorem strengthens the result of Proposition~\ref{uttirr}(b). We notice that it is possible to set $\mathcal{I}=[0,1)$ in Proposition~\ref{Lebesgue}.

\begin{pro}\label{Lebesgue} Let $\mathcal{I}\in\I$. For any $\epsilon>0$, there exists a nonempty closed set $E\subset\mathcal{I}$ such that $E$ is a finite union of disjoint closed turnpike intervals, and $\mu(\mathcal{I}\setminus E)<\epsilon$, where $\mu$ is the Lebesgue measure on $[0,1]$.
\end{pro}
\begin{proof}
    Let $\epsilon>0.$ Since $\mathcal{I}$ is a nonempty interval, there exist $0<a<b<1$ such that $E_1:=[a,b]\subset\mathcal{I}$ and $\mu(\mathcal{I}\setminus E_1)=\epsilon/3.$ By Proposition~\ref{uttirr}(b) and its proof, there exists a closed set $E_2\subset E_1$ such that $N^*(E_2)<\infty,$ $\mu(E_1\setminus E_2)<\epsilon/3,$ $\P(E_2)=\emptyset,$ and $E_2$ is a finitely disjoint union of closed intervals. If $\mathbb{D}(E_2)=\emptyset,$ then $E:=E_2\subset\mathcal{I}$ is a finitely disjoint union of closed turnpike intervals satisfying $\mu(\mathcal{I}\setminus E)=\mu(\mathcal{I}\setminus E_1)+\mu(E_1\setminus E_2)<2\epsilon/3<\epsilon.$ Otherwise, let $\mathbb{D}(E_2)=\{c_1,\ldots,c_k\},$ where the finiteness of $\mathbb{D}(E_2)$ is by Corollary~\ref{Ndinclu} and by Proposition~\ref{partition}(b). For each $1\le j\le k$, we choose some $\epsilon_j>0$ such that $\max_{1\le j\le k}\epsilon_j<\epsilon/(6k)$ and $\{[(c_j-\epsilon_j,c_j+\epsilon_j])\}_{j=1}^k$ are disjoint intervals. Let $E:=E_2\setminus\left(\cup_{i=1}^k(c_j-\epsilon_j,c_j+\epsilon_j)\right).$ Then $E\subset\mathcal{I}$ is closed,
    \[\mu(\mathcal{I}\setminus E)\le\mu(\mathcal{I}\setminus E_1)+\mu(E_1\setminus E_2)+\sum_{j=1}^k\mu((c_j-\epsilon_j,c_j+\epsilon_j))<\epsilon/3+\epsilon/3+\epsilon/3=\epsilon,\]
    and $\mathbb{D}(E)=\emptyset.$ Therefore, $E$ is a finite union of disjoint closed turnpike intervals.
\end{proof}


\section{Boundedness of the Turnpike Function Near Irregular Points}\label{sectbb}

In view of Theorem~\ref{utt}, if $\alpha$ is a regular point, then $N^*(\mathcal{I})<\infty$ for some interval $\mathcal{I}$ containing $\alpha.$ However, this can fail if $\alpha$ is an irregular point; \cite[Example 1]{doi:10.1287/moor.2017.0912}, Example~\ref{exatbbounded}. Throughout this section,  we fix an irregular point $\alpha\in\P((0,1))$ (unless otherwise stated) and investigate conditions for the boundedness of the turnpike function near this point. Boundedness of the turnpike function near $\alpha$ is described by the following equivalent statements: 

\begin{enumerate}[(i)]
    \item there exists $\delta>0$ such that $N^*((\alpha-\delta,\alpha+\delta))<\infty;$ 
    \item $\limsup_{\beta\to\alpha}N(\beta)<\infty.$
\end{enumerate}
Equivalence of these statements follows from Theorem~\ref{finitetpi}. When $\alpha$ is a non-touching break point, the negation of Condition~\hyperref[conC]{$\mathrm{C}$} below is a necessary condition of statement (i); \cite[Theorem 3]{doi:10.1287/moor.2017.0912}, but not a sufficient condition; Example~\ref{exatbbounded}. Moreover, \cite[Theorem 4]{doi:10.1287/moor.2017.0912} provides the necessary condition for statement (i), when $\alpha$ is a non-break touching point. The main result of this section, Theorem~\ref{tbsufficient}, provides the sufficient condition for the equivalent statements (i) and (ii), when $\alpha$ is an irregular point. In addition, Theorem~\ref{tbnecessary} provides a weaker necessary condition for a statement that is weaker than the equivalent statements (i) and (ii).

\begin{abc0}\hspace*{-0.2cm}{\rm (\cite[Condition C]{doi:10.1287/moor.2017.0912}})\label{conC} For any two distinct decision rules $\phi,\psi\in D(\alpha),$ there exists $x\in\X$ such that $T_\alpha^\phi V_{n-1}(x)\ne T_\alpha^\phi V_{n-1}(x)$ for infinitely many $n\in\N^+.$
\end{abc0}

\begin{con3}\label{conA}
There exists $K\in\N^+$ such that $D_n(\alpha)=D(\alpha)$ for all $n\ge K.$
\end{con3}

Condition \hyperref[conA]{$\mathrm{A}$} means that all infinite-optimal policies are also first-step-optimal for finite-horizon problems if the horizon is large enough. Furthermore, for $S\in 2^\F,$ define
\begin{equation}\label{piD}
    \Pi(S):=\left\{(\phi_0,\phi_1,\ldots)\in\Pi:\ \phi_t\in S\ \text{for all $t=0,1,\ldots$}\right\}.
\end{equation}
For $\phi\in\F$ and $\pi=(d_0,d_1,\ldots)\in\Pi,$ define $\phi^\pi:=(\phi,d_0,d_1,\ldots).$ In other words, $\Pi(D)$ is the set of policies whose decision rules at each step are from $D,$ and $\phi^\pi$ is the policy that uses the decision rule $\phi$ at the initial step and uses the decision rules $d_{t-1}$ at the steps $t=1,2,\ldots.$ We observe that, in view of \eqref{eqvaluesum2}, for every policy $\pi\in\Pi,$ the objective function $v_\beta^{\pi}$ is analytic in $\beta$ for $\beta\in[0,1).$ We consider the following conditions.

\begin{abc3}\label{thm:BB-} For all $\phi\in D(\alpha-)$ and $\psi\in D(\alpha)\setminus D(\alpha-)$, there exists $x\in\X$ such that
\[\sup_{\pi\in\Pi(D(\alpha))}\frac{d}{d\beta}\left(v_\beta^{\phi^\pi}-v_\beta^{\psi^\pi}\right)\Big|_{\beta=\alpha}(x)<0.\]\end{abc3}

\begin{abc4}\label{thm:BB+} For all $\phi\in D(\alpha+)$ and $\psi\in D(\alpha)\setminus D(\alpha+),$ there exists $x\in\X$ such that
\[\inf_{\pi\in\Pi(D(\alpha))}\frac{d}{d\beta}\left(v_\beta^{\phi^\pi}-v_\beta^{\psi^\pi}\right)\Big|_{\beta=\alpha}(x)>0.\]\end{abc4}

Conditions \hyperref[thm:BB-]{$\mathrm{B_-}$} and \hyperref[thm:BB+]{$\mathrm{B_+}$} are equivalent, if $\alpha$ is a non-touching break point, because $D(\alpha+)=D(\alpha)\setminus D(\alpha-)$ and  $D(\alpha-)=D(\alpha)\setminus D(\alpha+).$ 
Note that, if $\pi_1,\pi_2\in\Pi(D(\alpha)),$ then $v_{\alpha}^{\pi_1}=v_{\alpha}^{\pi_2}=V_{\alpha}$ by the optimality equation \eqref{eqoptimal}. Hence, Condition \hyperref[thm:BB-]{$\mathrm{B_-}$} and the finiteness of $D(\alpha)$ imply that there exist $x\in\X,$ $C>0,$ and $\delta>0$ such that
\begin{equation}\label{vsub}
v_{\beta}^{\phi^\pi}(x)-v_{\beta}^{\psi^\pi}(x)>C(\beta-\alpha)\qquad {\rm for}\ \beta\in(\alpha-\delta,\alpha),\ \ \pi\in\Pi(D(\alpha)).
\end{equation}
In other words, Condition \hyperref[thm:BB-]{$\mathrm{B_-}$} implies that, for discount factors in a left-neighborhood of $\alpha,$ 
a policy playing a decision rule from $D(\alpha-)$ at the first step and then playing a policy from $\Pi(D(\alpha))$ strictly dominates a policy playing a decision rule from $D(\alpha)\setminus D(\alpha-)$  at the first step and then playing the same policy from $\Pi(D(\alpha)).         
                          $ Similar argument holds for Condition \hyperref[thm:BB+]{$\mathrm{B_+}$}.

The following lemma is a well-known fact in dynamic programming.

\begin{lem}\label{lemrl} Let $n\in\N^+,$ $K\in\N^+,$ and $K\le n.$ Let $\pi=(\phi_0,\phi_1,\ldots)$ satisfies $\phi_t\in D_{n-K-t+1}(\alpha)$ for all $t=0,1,\ldots,n-K.$ Then 
\[V_{n,\alpha}=\sum_{t=0}^{n-K}\alpha^tP_t(\pi)r(\pi)+\alpha^{n-K+1}P_{n-K+1}(\pi)V_{K-1,\alpha}.\]
In particular, if $K=1,$ then $\pi\in M_n(\alpha).$
\end{lem}
\begin{proof}
    This statement follows from Definition~\ref{defMnDnD} and formulae~\eqref{maxvaluef}$-$\eqref{eqvaluesum1} and \eqref{eqoptimal}.
\end{proof}

The following theorem provides sufficient conditions for boundedness of the turnpike function near an irregular point; see Example~\ref{exatbb} for its application.

\begin{thm}\label{tbsufficient} Let $\alpha\in \P((0,1)).$ Then
\begin{enumerate}[(a)]
	\item if both Conditions \hyperref[conA]{$\mathrm{A}$} and \hyperref[thm:BB-]{$\mathrm{B_-}$} hold at $\alpha,$ then there exists $\delta>0$ such that $N^*((\alpha-\delta,\alpha))<\infty;$
	\item  if both Conditions \hyperref[conA]{$\mathrm{A}$} and \hyperref[thm:BB+]{$\mathrm{B_+}$} hold at $\alpha,$ then there exists $\delta>0$ such that $N^*((\alpha,\alpha+\delta))<\infty.$
\end{enumerate}
\end{thm}
\begin{proof} For simplicity, we only prove (a), as the proof of (b) is similar.\medskip

\noindent\textbf{Step 1.} Let
\begin{equation}\label{defiS}
    S_{k,\alpha}(\gamma):=(\alpha+\gamma)^{k+1}-\alpha^{k+1}-(k+1)\alpha^k\gamma,\qquad k\in\N,\ \gamma\in(-\alpha,0).
\end{equation}
By the Maclaurin expansion in $\gamma$ with the Lagrange remainder,  $S_{k,\alpha}(\gamma)=k(k+1){\lambda}^{k-1}\gamma^2/2$ for some $\lambda\in(\gamma,0).$ Since $\abs*{\lambda}<\abs*{\gamma}<\alpha$, 
\begin{equation}\label{taylor}
	\abs{S_{k,\alpha}(\gamma)}\le\frac{k(k+1)\alpha^{k-1}}{2}\gamma^2, \qquad  k\in\N,\ \gamma\in(-\alpha,0).
\end{equation}

\noindent\textbf{Step 2.} In view of Condition \hyperref[thm:BB-]{$\mathrm{B_-}$}, there exist $x\in\X$ and $C>0$ such that
\begin{equation}\label{selectx}
    \sup_{\pi\in\Pi(D(\alpha))}\frac{d}{d\beta}\left(v_\beta^{\phi^\pi}-v_\beta^{\psi^\pi}\right)\Big|_{\beta=\alpha}(x)<-C\quad{\rm for\ all}\ \phi\in D(\alpha-),\ \psi\in D(\alpha)\setminus D(\alpha-).
\end{equation}
By Corollary~\ref{absorb}, there exist $\delta_1>0$ and $K_1\in\N^+$ such that $D_n(\beta)\subset D(\alpha)$ for all $n\ge K_1$ and for all $\beta\in(\alpha-\delta_1,\alpha).$ By Proposition~\ref{uttirr}(a), $D(\beta)=D(\alpha-)$ for all $\beta\in(\alpha-\delta_1,\alpha),$ meaning that $(\alpha-\delta_1,\alpha)$ must be a subinterval of the partition interval whose right endpoint is $\alpha.$ Let $K_2$ be the constant $K$ from Condition \hyperref[conA]{$\mathrm{A}$}. Let $K_3=\max\{K_1,K_2\}.$ Let $\delta_2>0$ and $K_4\in\N$ be such that, for all $\beta\in(\alpha-\delta_2,\alpha)$ and for all $ n\ge K_4,$
\begin{equation}\label{complicated}
	\frac{4\alpha^{n+1}}{(1-\alpha)^2}+\frac{4(n+1)\alpha^n}{1-\alpha}+\left[\frac{2}{(1-\alpha)^3}+\frac{n(n+1)\alpha^{n-1}}{1-\alpha}\right](\alpha-\beta)<\frac{C}{R}.
\end{equation}
By Proposition~\ref{partition}(a) and by Lemma~\ref{lemlipschitz}, there exist $\delta_3>0$ and a vector function $u_\alpha$ such that 
\begin{equation}\label{vecextract}
\norm{u_\alpha}\le\frac{R}{(1-\alpha)^2},\quad {\rm and}\quad	V_{K_3-1,\beta}-V_{K_3-1\alpha}=(\beta-\alpha)u_\alpha \quad  {\rm for\ all}\  \beta\in(\alpha-\delta_3,\alpha). 
\end{equation}

\noindent\textbf{Step 3.} Let $\delta=\min\{\delta_1,\delta_2,\delta_3\}.$ Let $\beta\in(\alpha-\delta,\alpha)$ everywhere until the end of this proof, except formula \eqref{Vn} where the interval is broader, and let $n\ge K_3+K_4.$ Let $\pi=(\phi_0,\phi_1,\ldots)\in\Pi(D(\alpha))$ such that $\phi_t\subset D_{n-t-1}(\beta)$ for all $t=0,1,\ldots,n-2.$ Then $v_{n-1,\beta}^\pi=V_{n-1,\beta}$ by Lemma~\ref{lemrl}. Since $(n-1)-(n-K_3-1)=K_3=\max\{K_1,K_2\},$ we have $\phi_t\subset D_{n-t-1}(\beta)\subset D(\alpha)=D_{n-t-1}(\alpha)$ for all $t=0,1,\ldots,n-K_3-1,$ where the second inclusion is by the definition of $K_1,$ and the last equality is by the definition of $K_2.$ By Lemma~\ref{lemrl}, this implies
\begin{equation}\label{Vn}
V_{n-1,\beta}=\sum_{t=0}^{n-K_3-1}\beta^tP_t(\pi)r(\pi)+\beta^{n-K_3}P_{n-K_3}(\pi)V_{K_3-1,\beta}, \qquad \beta\in (\alpha-\delta,\alpha].
\end{equation}

Let us define
\begin{align*}
    I_1 & := r(\phi)-r(\psi)+\alpha\left[P(\phi)-P(\psi)\right]V_{n-1,\alpha};\\
	I_2 & := \alpha^{n-K_3}\left[P(\phi)-P(\psi)\right]P_{n-K_3}(\pi)\left[\alpha\left(V_{K_3-1,\beta}-V_{K_3-1,\alpha}\right)+(n-K_3+1)(\beta-\alpha)V_{K_3-1,\beta}\right];\\
	I_3 & := (\beta-\alpha)\left[P(\phi)-P(\psi)\right]\sum_{t=0}^{n-K_3-1}(t+1)\alpha^tP_t(\pi)r_t(\pi);\\
	I_4 
  & := [P(\phi)-P(\psi)](\beta V_{n-1,\beta} -\alpha V_{n-1,\alpha})-I_2-I_3 \\
    & = \left[P(\phi)-P(\psi)\right]\left[\sum_{t=0}^{n-K_3-1}S_{t,\alpha}(\beta-\alpha)P_t(\pi)r_t(\pi)+S_{n-K_3,\alpha}(\beta-\alpha)P_{n-K_3}(\pi)V_{K_3-1,\beta}\right],
 \end{align*}
where, to verify the last equality for $I_4,$ one can drop the common multiplier $[P(\phi)-P(\psi)],$ write the terms $\left\{S_{t,\alpha}(\beta-\alpha)\right\}_{t=0}^{n-K_3}$ explicitly in the form of \eqref{defiS}, drop the common summands on both sides, and write $V_{n-1,\alpha},$ $V_{n-1,\beta}$ in the form of \eqref{Vn}. Then
\[T_\beta^\phi V_{n-1,\beta}-T_\beta^\psi V_{n-1,\beta}=r(\phi)-r(\psi)+\beta[P(\phi)-P(\psi)]V_{n-1,\beta}=I_1+I_2+I_3+I_4.\]

\noindent\textbf{Step 4.} Since $\phi\in D_n(\alpha),$
\begin{equation}\label{errorS1}
	I_1(x)=T_{\alpha}^\phi V_{n-1,\alpha}(x)-T_{\alpha}^\psi V_{n-1,\alpha}(x) \ge0.
\end{equation}
Since $\norm{P(\phi)-P(\psi)}\le2,$ $\norm{P_{n-K_3}(\pi)}\le1,$ and $\beta<\alpha,$ by \eqref{vecextract} and by Lemma~\ref{lembound}(a) we have
\begin{align}\label{errorS2}
	\norm{I_2}&\le2\alpha^{n-K_3}\left[\alpha\norm{V_{K_3-1,\beta}-V_{K_3-1,\alpha}}+(n-K_3+1)(\alpha-\beta)\norm{V_{K_3-1,\beta}}\right]\\\nonumber
    &\le 2(\alpha-\beta)\alpha^{n-K_3}\left[\frac{\alpha}{(1-\alpha)^2}+\frac{(n-K_3+1)}{1-\alpha}\right]R.
\end{align}
Since $\norm{\left[P_t(\phi^\pi)r_t(\phi^\pi)-P_t(\psi^\pi)r_t(\psi^\pi)\right]}\le 2R$ for all $t\in\N,$ by \eqref{selectx} and \eqref{eqvaluesum2},
\begin{equation}
\begin{aligned}\label{errorS3}
    &I_3(x)=(\beta-\alpha)\left[P(\phi)\sum_{t=0}^{n-K_3-1}(t+1)\alpha^tP_t(\pi)r_t(\pi)(x)-P(\psi)\sum_{t=0}^{n-K_3-1}(t+1)\alpha^tP_t(\pi)r_t(\pi)(x)\right]\\
    &=(\beta-\alpha)\left\{\sum_{t=0}^{n-K_3}t\alpha^{t-1}\left[P_t(\phi^\pi)r_t(\phi^\pi)-P_t(\psi^\pi)r_t(\psi^\pi)\right](x)\right\}\\
    &=(\beta-\alpha)\left\{\frac{d}{d\alpha}\left(v_\alpha^{\phi^\pi}-v_\alpha^{\psi^\pi}\right)(x)-\sum_{t=n-K_3+1}^\infty t\alpha^{t-1}\left[P_t(\phi^\pi)r_t(\phi^\pi)-P_t(\psi^\pi)r_t(\psi^\pi)\right](x)\right\}\\
    & >(\alpha-\beta)\left\{C-2R\sum_{t=n-K_3+1}^\infty t\alpha^{t-1}\right\} =(\alpha-\beta)\left\{C-2\alpha^{n-K_3}\left[\frac{\alpha}{{(1-\alpha)}^2}+\frac{n-K_3+1}{1-\alpha}\right]R\right\}.
\end{aligned}
\end{equation}
Similarly, by \eqref{taylor}, by Lemma~\ref{lembound}(a) and since $\beta<\alpha,$
\begin{equation}
\begin{aligned}
	\norm{I_4} & \le 2\left[\sum_{t=0}^{n-K_3+1}\frac{t(t+1)\alpha^{t-1}R}{2}{(\beta-\alpha)}^2+\frac{(n-K_3)(n-K_3+1)\alpha^{n-K_3-1}R}{2(1-\alpha)}{(\beta-\alpha)}^2\right]\\
	& \le {(\alpha-\beta)}^2\left[\frac{2}{(1-\alpha)^3}+\frac{(n-K_3)(n-K_3+1)\alpha^{n-K_3-1}}{1-\alpha}\right]R.\label{errorS4}
\end{aligned}
\end{equation}
\noindent\textbf{Step 5.} Since $\beta\in(\alpha-\delta,\alpha)$ and $n-K_3\ge K_4$, by \eqref{complicated}, \eqref{errorS1}$-$\eqref{errorS4}, and by the choice of $I_j$, $j=1,2,3,4,$
\begin{align*}
	&T_\beta^\phi V_{n-1,\beta}(x)-T_\beta^\psi V_{n-1,\beta}(x)=I_1(x)+I_2(x)+I_3(x)+I_4(x)\ge I_3(x)-\norm{I_2}-\norm{I_4}&&\nonumber\\
	&\qquad\qquad\qquad\qquad> C(\alpha-\beta)-(\alpha-\beta)\left\{\frac{4\alpha^{n-K_3+1}}{(1-\alpha)^2}+\frac{4(n-K_3+1)\alpha^{n-K_3}}{1-\alpha}+(\alpha-\beta)\left[\frac{2}{(1-\alpha)^3}\right.\right.&&\nonumber\\
	& \qquad\qquad\qquad\qquad\qquad\qquad\ +\left.\left.\frac{(n-K_3)(n-K_3+1)\alpha^{n-K_3-1}}{1-\alpha}\right]\right\}R> C(\alpha-\beta)-C(\alpha-\beta)=0.\nonumber
\end{align*}
Since $\phi\in D(\alpha-), \psi\in D(\alpha)\setminus D(\alpha-), \beta\in(\alpha-\delta,\alpha), n\ge K_3+K_4$ are arbitrary and by the definition of $K_1,$ we have $D_n(\beta)\subset D(\alpha-)=D(\beta)$ for all $\beta\in(\alpha-\delta,\alpha)$ and for all $n\ge K_3+K_4.$ Therefore, $N^*((\alpha-\delta,\alpha))\le K_3+K_4.$
\end{proof}

Example~\ref{exatbb} demonstrates how Theorem~\ref{tbsufficient} can be applied to prove the boundedness of the turnpike function near an irregular point. Example~\ref{exatbbounded} demonstrates that Conditions~\hyperref[thm:BB-]{$\mathrm{B_-}$} and \hyperref[thm:BB+]{$\mathrm{B_+}$} are essential in Theorem~\ref{tbsufficient}. We note that Theorem~\ref{tbsufficient} also provides sufficient conditions for $\limsup_{\beta\uparrow\alpha}N(\alpha)<\infty$ and $\limsup_{\beta\downarrow\alpha}N(\alpha)<\infty$ respectively. We recall that $\alpha\in\P((0,1))$ and consider the following conditions.



\begin{abc1}\label{thm:Ca-} There exists $K\in\N^+$ such that $D_n(\alpha)\cap D(\alpha-)\ne\emptyset$ for all $n\ge K$.\end{abc1}

\begin{abc2}\label{thm:Ca+} There exists $K\in\N^+$ such that $D_n(\alpha)\cap D(\alpha+)\ne\emptyset$ for all $n\ge K$.\end{abc2}

Condition~\hyperref[thm:Ca-]{$\mathrm{A_-}$} means that starting from some $K\in\N^+$, for each horizon $n\ge K$, there exists a decision rule $\phi^{(n)},$ which is both first-step-optimal for the $n$-horizon problem  and left-optimal. The similar interpretation holds for Condition~\hyperref[thm:Ca+]{$\mathrm{A_+}$}. Condition~\hyperref[conA]{$\mathrm{A}$} implies both Conditions~\hyperref[thm:Ca-]{$\mathrm{A_-}$} and \hyperref[thm:Ca+]{$\mathrm{A_+}$}. In view of the definition of the turnpike integer $N(\alpha),$ the converse is true, if $\alpha$ is a non-touching  point and both $D(\alpha-)$ and $D(\alpha+)$ are singletons. If $\alpha$ is a regular point, then both Conditions~\hyperref[thm:Ca-]{$\mathrm{A_-}$} and \hyperref[thm:Ca+]{$\mathrm{A_+}$} hold, but not necessarily Condition~\hyperref[conA]{$\mathrm{A}$}; see \cite[Appendix C]{TSENG1990287}. The following theorem states necessary conditions for the boundedness of lower limits of the turnpike function near an irregular point.
\begin{thm}
\label{tbnecessary} Let $\alpha\in\P((0,1)).$ Then
\begin{enumerate}[(a)]
	\item if $\liminf_{\beta\uparrow\alpha}N(\alpha)<\infty$, then Condition~\hyperref[thm:Ca-]{$\mathrm{A_-}$} holds at $\alpha;$
	\item if $\liminf_{\beta\downarrow\alpha}N(\alpha)<\infty$, then Condition~\hyperref[thm:Ca-]{$\mathrm{A_+}$} holds at $\alpha.$
\end{enumerate}
\end{thm}
The proof of Theorem~\ref{tbnecessary} is similar to the proof of \cite[Theorems 3 and 4]{doi:10.1287/moor.2017.0912}. We provide it for completeness.

    \begin{proof}{Proof of Theorem~\ref{tbnecessary}.}  Suppose Condition \hyperref[thm:Ca-]{$\mathrm{A_-}$} does not hold at $\alpha.$ Then there exists an increasing sequence of positive integers $\{n_i\}_{i=1}^\infty\subset\N^+$ such that $D_{n_i}(\alpha)\cap D(\alpha-)=\emptyset$. Thus, for each $i\in\N^+$ there exists $\psi^{(i)}\in D_{n_i}(\alpha)\setminus D(\alpha-)$ and $y^i\in\X$ such that $V_{n_i,\alpha}(y^i)=T_{\alpha}^{\psi^{(n_i)}}V_{n_i-1,\alpha}(y^i)>T_{\alpha}^{\phi}V_{n_i-1,\alpha}(y^i)$ for all $\phi\in D(\alpha-).$ Since $\F$ and $\X$ are finite, there exists a subsequence $n_{i_j}$ and $\psi\in\F\setminus D(\alpha-),x\in\X$ such that $\phi^{(n_{i_j})}=\psi$ and $y^{n_{i_j}}=x.$ Then $V_{n_{i_j},\alpha}(x)=T_{\alpha}^{\psi}V_{{n_{i_j}}-1,\alpha}(x)>T_{\alpha}^{\phi}V_{n_{i_j}-1,\alpha}(x)$ for all $\phi\in D(\alpha-).$ Continuity in $\beta$ of the function $T_\beta^\phi V_{n-1,\beta}(x)$ stated in Proposition~\ref{partition}(b) and the finiteness of $\F$  imply that for each $j\in\N^+$ there exists some $\delta_j>0$ such that $T_\beta^\psi V_{n_{i_j}-1,\beta}(x)>T_\beta^\phi V_{n_{i_j}-1,\beta}(x)$ for all $\phi\in D(\alpha-)$ and for all $\beta\in(\alpha-\delta_i,\alpha)$. This implies $N^*((\alpha-\delta_j,\alpha))>n_{i_j}$ for all $j\in\N^+$, which contradicts $\liminf_{\beta\uparrow\alpha}N(\alpha)<\infty$. This proves (a), and the proof of (b) is similar.
    \end{proof}


\begin{defi}\label{defsbp}\emph{
    We say that a break point $\alpha\in(0,1)$ is a \textbf{simple break point} if it is not touching, and both $D(\alpha-)$ and $D(\alpha+)$ are singletons.}
\end{defi}

As explained following the introduction of Conditions~\hyperref[thm:BB+]{$\mathrm{B_+}$}, if $\alpha$ is a simple break point, then Conditions~\hyperref[thm:BB-]{$\mathrm{B_-}$} and \hyperref[thm:BB+]{$\mathrm{B_+}$} are equivalent. The following corollary states that, if $\alpha$ is a simple break point satisfying Condition~\hyperref[thm:BB-]{$\mathrm{B_-}$}, then the turnpike function has either a left or a right limit at $\alpha.$

\begin{cor} Let $\alpha$ be a simple break point. If Condition~\hyperref[thm:BB-]{$\mathrm{B_-}$} holds at $\alpha$, then either $\liminf_{\beta\uparrow\alpha}N(\beta)=\limsup_{\beta\uparrow\alpha}N(\beta)$ or $\liminf_{\beta\downarrow\alpha}N(\beta)=\limsup_{\beta\downarrow\alpha}N(\beta).$
\end{cor}
\begin{proof}
    Suppose neither is true. By Theorem~\ref{finitetpi}, this implies $\liminf_{\beta\uparrow\alpha}(\beta)<\infty,$ $\liminf_{\beta\downarrow\alpha}(\beta)<\infty,$ and $\limsup_{\beta\to\alpha}N(\beta)=\infty.$ By Theorem~\ref{tbnecessary} this implies both Conditions~\hyperref[thm:Ca-]{$\mathrm{A_-}$} and \hyperref[thm:Ca-]{$\mathrm{A_+}$} hold at $\alpha,$ which implies Condition~\hyperref[conA]{$\mathrm{A}$} holds at $\alpha.$ By Theorem~\ref{tbsufficient} this implies that there exists $\delta>0$ such that $N((\alpha-\delta,\alpha+\delta))<\infty,$ which is a contradiction.
\end{proof}

\section{Turnpike Properties for Small Discount Factors}\label{secsmalltp}
This section studies properties of the turnpike function for small discount factors, for which Proposition~\ref{smallpar} estimates the sets of first-step-optimal decision rules and optimal decision rules, and Theorem~\ref{ndboundm} provides an upper bound of the turnpike function. For finite MDPs, \cite[Theorem 2]{CB} implies that for every $\epsilon>0$ the value iteration algorithm finds an $\epsilon$-optimal policy after one iteration if the positive discount factor is sufficiently small. Theorem~\ref{ndboundm}(a) shows that, if the discount factor is smaller than $\Delta_L,$ where $\Delta_L$ is a positive number defined by formulae~\eqref{FXdef}$-$\eqref{Deldef}, then the number of value iterations for finding an optimal policy is bounded by the number of states. Example~\ref{exanbm} shows that this upper bound is sharp.

Let us consider the definitions of maximum absolute values in \eqref{eqR}. The bounds derived in this paper are more favorable if these maximum absolute values are smaller. For a given MDP, we can minimize the maximum absolute values by modifying the reward by subtracting constant vectors $F_1\pmb{1}$ and $F_2\pmb{1}$ from $r(\cdot,\cdot)$ and $s(\cdot),$ respectively, where
    \[F_1:=\frac{1}{2}\left\{\max_{x\in\X,a\in A(x)}r(x,a)+\min_{x\in\X,a\in A(x)}r(x,a)\right\},\qquad F_2:=\frac{1}{2}\left\{\max_{x\in\X}s(x)+\min_{x\in\X}s(x)\right\}.\]
    By processing these subtractions, the objective functions of each policy for infinite-horizon problems and each $n$-horizon problem are subtracted by the vectors $\frac{F_1}{1-\alpha}\pmb{1}$ and $\left[\frac{F_1(1-\alpha^n)}{1-\alpha}+\alpha^nF_2\right]\pmb{1}$ respectively. Hence, the sets of optimal policies remain the same for both finite-horizon problems and infinite-horizon problems. The new maximum absolute values, denoted by $R_1^*,$ $R_2^*$ and $R^*$ respectively, become
\begin{equation}
    \begin{aligned}\label{eqRb}
         R_1^* & := \frac{1}{2}\left\{\max_{x\in\X,a\in A(x)}r(x,a)-\min_{x\in\X,a\in A(x)}r(x,a)\right\} = R_1-\abs*{F_1}\le R_1,\\
         R_2^* & := \frac{1}{2}\left\{\max_{x\in\X}s(x)-\min_{x\in\X}s(x)\right\} = R_2-\abs*{F_2}\le R_2,\qquad R^*:=\max\left\{R_1^*,R_2^*\right\}\le R.
   \end{aligned}
   \end{equation}
    We say that one-step rewards are \textbf{balanced} if $R_1=R_1^*,$ and all rewards are balanced if $R_1=R_1^*$ and $R_2=R_2^*.$  
The following lemma shows that balancing rewards does not change the properties of discount factors, turnpike functions, and sets of optimal policies.

\begin{lem}\label{lembalance} For a given MDP, a given discount factor $\alpha\in[0,1)$ and a given interval $\mathcal{I}\in\I,$
    \begin{enumerate}[(a)]
        \item balancing all rewards does not change the points, sets,  partitions in Definitions~\ref{defMnDnD}, \ref{def+-}, \ref{defirr}, \ref{defirrfinite}, \ref{defpartition}, \ref{deftpi}, \ref{defdiscon} as well as the value of the turnpike integer defined in Definition~\ref{deftp};
        \item the sets $D(\alpha),$ $D(\alpha-),$ $D(\alpha+),$ the set of points $\{a_i\}_{i=0}^l$ in Proposition~\ref{partition}(c), and the canonical partition do not depend on the terminal reward vector $s,$ and therefore balancing one-step rewards changes neither these sets nor the canonical partition.
    \end{enumerate}
\end{lem}
\begin{proof}
    (a) follows from the definition of balanced rewards; (b) follows from equations \eqref{eqvaluesum2}, \eqref{eqvaluesum3} and from Definition~\ref{defMnDnD}.
\end{proof}

Let $\F_{-1}:=\F.$ For $n\in\N,$ let
\begin{equation}
	\begin{aligned}\label{FXdef}
		\F_n:&=\left\{\phi\in\F_{n-1}: P^n(\phi)r(\phi)(x)\ge P^n(\psi)r(\psi)(x)\ \text{for all\ } \psi\in\F_{n-1}\ \text{and for all\ }x\in\X\right\},\\
		\X_n:&=\left\{x\in\X: P^n(\phi)r(\phi)(x)\ne P^n(\psi)r(\psi)(x)\ \text{for\ some\ }\phi,\psi\in\F_{n-1}\right\}.
	\end{aligned}
    \end{equation}
	We observe that $\F_{n-1}\supset\F_n\ne\emptyset$ and $\{P^n(\phi)r(\phi)\}_{\phi\in\F_n}$ is a singleton for each $n\in\N$; see \eqref{defaltern1}, \eqref{defaltern2} for explicit definitions of $\F_n$ and $\X_n.$ Also we note that $\{\F_t\}_{t=-1}^\infty$ is a non-increasing sequence of finite sets. Let $H\in\N$ be the number of elements in the set $\{n\in\N^+: \F_n\ne\F_{n-1}\}.$ We define $L_0:=0,$ and for $i=1,\ldots,H$ if $H>0,$ let $L_i$ be the step number of $i$-th decrease of the set $\F_n,$ 
    \begin{equation}
    \begin{aligned}\label{Ldef}
        L_i:&=\min\{n\in\N^+: j>L_{n-1}\ \text{and\ } \F_n\ne\F_{n-1}\},\quad i=1,2,\ldots,H;\\ 
        L:&=L_H=\min\left\{K\in\N: \text{$\F_n=\F_{n+1}$ for all $n\ge K$}\right\}<\infty.
    \end{aligned}
    \end{equation}
    Let $L_{H+1}:=\infty.$ 
    For $n\in\N$ we note that $\X_n=\emptyset$ iff $\F_n=\F_{n-1}$. Therefore, if $\X_n\ne\emptyset,$ then $n=L_i$ for some $0\le i\le H$, and the converse is true if $n\in\N^+.$ Also $\X_n=\emptyset$ for all $n\ge L+1.$ Clearly $H\le L.$ We have that $L_0=0,$ $L_H=L,$  $\F_{L_i}\ne\F_{L_{i+1}},$ and $\F_{L_i+k}=\F_{L_{i+1}-1},$ if $0\le i\le H-1$ and $0\le k\le L_{i+1}-L_i-1.$

    By definition, if $\phi,\psi\in\F_{m-1},$ then $P^n(\phi)r(\phi)=P^n(\psi)r(\psi)$ for all $n=0,1,\ldots,m-1,$ where $m$ is the number of states, which, by Proposition~\ref{thmla}, implies $P^n(\phi)r(\phi)=P^n(\psi)r(\psi)$ for all $n\in\N.$ Therefore, $\F_n=\F_{n+1}$ for all $n\ge m-1,$ which implies $L\le m-1.$
    

    \begin{lem}\label{lempower} Let $0\le i\le H,$ $0\le j<L_{i+1},$ and $\phi^{(1)},\ldots,\phi^{(j)}\in\F_{L_i}.$ Then
    \[\left[\prod_{k=1}^jP(\phi^{(k)})\right]r(\phi^{(j)})=P^{j}(\phi)r(\phi)\qquad {\rm for\  all\ } \phi\in\F_{L_i}.\]
    \end{lem}
    \begin{proof}
        This is by the definition of $\{\F_{L_k}\}_{k=0}^H$ and by the definitions of $L_0,L_1,\ldots,L_H,L_{H+1}.$ 
    \end{proof}

    \begin{lem}\label{lemFpower} Consider $0\le i\le H$ and $\pi=(\phi_0,\phi_1,\ldots)$ such that $\{\phi_q\}_{q=0}^{L_{i+1}-L_k-1}\subset\F_{L_k}$ for all $k=0,1,\ldots,i.$ Then $P_t(\pi)r_t(\pi)=P^t(\phi)r(\phi)$ for all $t=0,1,\ldots,L_{i+1}-1$ and for all $\phi\in\F_{L_i}.$
    \end{lem}

    \begin{proof} 
				The proof is by induction. If $i=0$ then the statement follows from Lemma~\ref{lempower}. Suppose the statement is true for all $i=0,1,\ldots,j-1,$ where $1\le j\le H.$ Let $i=j$ and $\{\phi_q\}_{q=0}^{L_{j+1}-L_k-1}\subset\F_{L_k}$ for all $k=0,1,\ldots,j.$ Then in particular $\{\phi_q\}_{q=0}^{L_j-L_k-1}\subset\{\phi_q\}_{q=0}^{L_{j+1}-L_k-1}\subset\F_{L_k}$ for all $k=0,1,\ldots,j-1,$ which by the inductive assumption of the statement for $i=j-1$ implies $P_t(\pi)r_t(\pi)=P^t(\phi)r(\phi)$ for all $t=0,1,\ldots,L_j-1$ and for all $\phi\in\F_{L_{j-1}}\supset\F_{L_j}.$ If $L_j\le t\le L_{j+1}-1,$ we let $\tilde{\pi}=(\psi_0,\psi_1,\ldots)=(\phi_{t-L_j+1},\phi_{t-L_j+2},\ldots)$ and note that $\{\psi_q\}_{q=0}^{L_j-L_k-1}=\{\phi_{q+t-L_j+1}\}_{q=0}^{L_j-L_k-1}=\{\phi_q\}_{q=t-L_j+1}^{t-L_k}\subset\{\phi_q\}_{q=0}^{L_{j+1}-L_k-1}\subset\F_{L_k}$ for all $k=0,1,\ldots,j-1.$ The inductive assumption for $i=j-1$ implies $P_{L_j-1}(\tilde{\pi})r_{L_j-1}(\tilde{\pi})=P^{{L_j-1}}(\phi)r(\phi)$ for all $\phi\in\F_{L_{j-1}}\supset\F_{L_j}.$ We also have $\{\phi_q\}_{q=0}^{t-L_j}\subset\{\phi_q\}_{q=0}^{L_{j+1}-L_j-1}\subset\F_{L_j}.$ Thus, $P_t(\pi)r_t(\pi)=P^t(\phi)r(\phi)$ for all $t=0,1,\ldots,L_j-1$ and for all $\phi\in\F_{L_{j-1}}\supset\F_{L_j},$ and
                for $t=L_j,L_j+1,\ldots,L_{j+1}-1,$
                \begin{equation*}
				\begin{aligned}
	  P_t(\pi)r_t(\pi)&=\left[\prod_{k=0}^{t-L_j}P(\phi_k)\right]\left[\prod_{k=t-L_j+1}^{t-1}P(\phi_k)\right]r(\phi_t)=\left[\prod_{k=0}^{t-L_j}P(\phi_k)\right]P_{L_j-1}(\tilde{\pi})r_{L_j-1}(\tilde{\pi})\\
                    &=\left[\prod_{k=0}^{t-L_j}P(\phi_k)\right]P^{{L_j-1}}(\phi)r(\phi)=P^t(\phi)r(\phi) \quad\text{for all } \phi\in\F_{L_j},
				\end{aligned}
                \end{equation*}
				where the last equality follows from Lemma~\ref{lempower}. This completes the induction. 
			\end{proof}

		Let $C_{-1}=\infty.$ For $n\in\N,$ let
		\begin{equation}\label{Cdef}
		    C_n:=\begin{cases}
			C_{n-1} &\text{if}\ \X_n=\emptyset,\\
			\min\left\{C_{n-1},\min\limits_{\substack{x\in\X_n \\ \psi\in\F_{n-1}\setminus\F_n}}\left\{P^n(\phi)r(\phi)(x)-P^n(\psi)r(\psi)(x): \phi\in\F_n\right\}\right\} &\text{if}\ \X_n\ne\emptyset.
		\end{cases}
		\end{equation}
		Recall that $R_1^*,R^*$ are defined in \eqref{eqRb}. For $n\in\N,$ let
		\begin{equation}\label{Deldef}
		    \begin{cases}
		        \displaystyle\Delta_n:=\frac{C_n}{2R^*+C_n}\ \ \text{and}\ \ \tilde{\Delta}_n:=\frac{C_n}{2R_1^*+C_n}  &\text{if $C_n<\infty$},\\
          \displaystyle\Delta_n=\tilde{\Delta}_n:=1 &\text{if $C_n=\infty$}.
		    \end{cases}
		\end{equation}
		We note that $\{C_t\}_{t=-1}^\infty,$ $\{\Delta_t\}_{t=-1}^\infty,$ and $\{\tilde{\Delta}_t\}_{t=-1}^\infty$ are three non-increasing sequences of positive values. Also, $\Delta_n\le\tilde{\Delta}_n$ for all $n\in\N,$ and for $n\ge L$ we have $C_n=C_{n+1},$ $\Delta_n=\Delta_{n+1},$ and $\tilde{\Delta}_n=\tilde{\Delta}_{n+1}.$ 
        In view of \eqref{eqRb} and \eqref{Cdef}, $C_n\le C_{L_1}\le 2R_1^*\le 2R^*$ if $n\ge L_1.$
       Therefore, if $\X_i\ne\emptyset$ for some $i\in\N,$ then $\Delta_n\le\tilde{\Delta}_n\le0.5$ for all $n\ge i.$ If $\X_i=\emptyset$ for all $i\in\N,$ then
     all decision rules are equivalent for the infinite-horizon problem, meaning that $C_n=
     \infty$ and therefore $\displaystyle\Delta_n=\tilde{\Delta}_n=1$ for all $n\in\N.$ 
    
	   The following proposition shows that, for small discount factors, all optimal decision rules and all first-step-optimal decision rules for large horizon problems belong to one of the sets from $\{\F_t\}_{t=0}^H.$ 

			\begin{pro}\label{smallpar} Let $0\le i\le H.$ Then
            \begin{enumerate}[(a)]
                \item $D_n((0,\Delta_{L_i}))\subset\F_{L_i}$ for all $n\ge L_i+1;$
                \item $D((0,\tilde{\Delta}_{L_i}))\subset\F_{L_i}.$
            \end{enumerate}				
			\end{pro}
			\begin{proof} Let us prove (a) by induction. By Lemma~\ref{lembalance}(a) we can assume that all rewards are balanced, and $R=R^*.$ Let $i=0.$ If $\F_0=\F,$ then $\F_{L_0}=\F_{0}=\F,$ and statement (a) holds for an arbitrary subset of $
            \F$ in its left-hand side. Suppose $\F_0\ne\F.$ Then $\X_0\ne\emptyset.$ Let $\phi\in\F_0,$ $\psi\in\F\setminus\F_{0}\ne\emptyset,$ $\alpha\in(0,\Delta_0),$ and $n\ge1.$ Then by the definitions of $\F_0,$ $\X_0,$ and $C_0,$ by the triangle inequality, and by Lemma~\ref{lembound}(a), there exists  $x\in\X_0$ such that
            \begin{equation}
				\begin{aligned}T_\alpha^\phi V_{n-1,\alpha}(x)-T_\alpha^{\psi}V_{n-1,\alpha}(x)&=r(\phi)(x)-r(\psi)(x)+\alpha\left[P(\phi)-P(\psi)\right]V_{n-1,\alpha}(x)\\
					&\ge C_0-\frac{2\alpha R^*}{1-\alpha}>0,\label{smallsubineq0}
				\end{aligned}
                \end{equation}
                where the last inequality follows from $0<\alpha<\Delta_0=C_0/(2R^*+C_0).$ Since $\psi\in\F\setminus\F_{0}$ is arbitrary, \eqref{smallsubineq0} implies that $D_n((0,\Delta_0))\subset\F_0$ for all $n\ge1.$ Thus, (a) is true for $i=0.$ Suppose (a) is true for $i=0,1,\ldots j-1$ for some $1\le j\le H.$ Let $i=j.$ Let $\alpha\in(0,\Delta_{L_j})$ and $n\ge L_j+1.$ Therefore, $(n-1)-(L_j-L_k-1)\ge L_k+1$ for  $k=0,1,\ldots,j-1.$ By Lemma~\ref{lemrl} and by the inductive assumption, there exists $\pi=(\phi_0,\phi_1,\ldots)\in M_{n-1}(\alpha)$ such that $\{\phi_t\}_{t=0}^{L_j-L_k-1}\subset\F_{L_k}$ for all $k=0,1,\ldots,j-1.$ Let $\phi\in\F_{L_{j-1}}.$ Then by Lemma~\ref{lemFpower} and by formula~\eqref{eqvaluesum1},
				\begin{align}
					&\quad T_\alpha^\phi V_{n-1,\alpha}=T_\alpha^\phi v_{n-1,\alpha}^\pi=r(\phi)+\alpha P(\phi)v_{n-1,\alpha}^\pi\nonumber\\
					&=r(\phi)+\alpha P(\phi)\sum_{t=0}^{L_j-1}\alpha^tP_t(\pi)r_t(\pi)+\alpha^{L_j+1}P(\phi)P_{L_j}(\pi)V_{n-L_j-1,\alpha}\label{smallsubb}\\
&=\sum_{t=0}^{L_j}\alpha^tP^t(\phi)r(\phi)+\alpha^{L_j+1}P_{L_j+1}(\phi^\pi)V_{n-L_j-1,\alpha}\qquad\text{for all } \phi\in\F_{j-1}.\nonumber
				\end{align}
                We note that $\F_{L_j}\ne\F$ since $j\ge 1.$ Let $\phi\in\F_{L_j}\subset\F_{L_{j-1}}$ and $\psi\in\F_{L_{j-1}}\setminus\F_{L_j}\ne\emptyset.$ Then $\X_{L_j}\ne\emptyset.$ By equation~\eqref{smallsubb}, by the definitions of $\F_{L_j},$ $\F_{L_{j-1}},$ and $C_{L_j},$ by the triangle inequality, and by Lemma~\ref{lembound}(a), there exists  $x\in\X_{L_j}$ such that
                \begin{equation}
				\begin{aligned}&\quad T_\alpha^\phi V_{n-1,\alpha}(x)-T_\alpha^{\psi}V_{n-1,\alpha}(x)\\
					&=\alpha^{L_j}\left[P^{L_j}(\phi)r(\phi)-P^{L_j}(\psi)r(\psi)\right](x)+\alpha^{L_j+1}\left[P_{L_j+1}(\phi^\pi)-P_{L_j+1}(\psi^\pi)\right]V_{n-L_j-1,\alpha}(x)\\
					&\ge\alpha^{L_j}\left(C_{L_j}-\frac{2\alpha R^*}{1-\alpha}\right)>0.\label{smallsubineq}
				\end{aligned}
                \end{equation}
				Since $D_n((0,\Delta_{L_j}))\subset D_n((0,\Delta_{L_{j-1}}))\subset\F_{L_{j-1}}$ by the inductive assumption, and since $\alpha\in(0,\Delta_{L_j}),$ $n\ge L_j+1,$ $\psi\in\F_{L_{j-1}}\setminus\F_{L_j}$ are arbitrary, \eqref{smallsubineq} implies that $D_n((0,\Delta_{L_j}))\subset\F_{L_j}$ for all $n\ge L_j+1.$ Therefore, (a) is proved by induction.
				
				Let us prove (b). If $i=0$ and $\F_0=\F$ then statement (b) holds for an arbitrary subset of $\F$ in its left-hand side. Suppose $i\ne0$ or $\F_0\ne\F.$ Let $L_{-1}:=-1.$ By Lemma~\ref{lembalance}(b), let us assume without loss of generality that $s=\pmb{0}.$ Then $R_1^*=R^*=R$ and $\Delta_{L_i}=\tilde{\Delta}_{L_i}$ for all $i=0,1,\cdots,H.$ By letting $i=j$ and taking $n\to\infty$ in formulae \eqref{smallsubineq0} and \eqref{smallsubineq}, we conclude that for each $i=0,1,\ldots,H,$  if $\alpha\in(0,\tilde{\Delta}_{L_i}),$ and $\phi\in\F_{L_i},\psi\in\F_{L_{i-1}}\setminus\F_{L_i},$ then there exists $x\in\X_{L_i}$ such that $T_\alpha^\phi V_\alpha(x)>T_\alpha^{\psi}V_\alpha(x).$ By the optimality equation \eqref{eqoptimal} and since $\F=\F_{L_{-1}}\supset\F_{L_0}\supset\F_{L_1}\supset\ldots\supset\F_{L_H},$ this implies that $D((0,\tilde{\Delta}_{L_i}))\subset\F_{L_i}$ for all $i=0,1,\ldots,H.$ 
			\end{proof}

            Recall from Definition~\ref{def+-} that there exists $\delta\in (0,1)$ such that $D(0+)=D((0,\delta)).$ Statement (b) of the following corollary provides a lower estimate of  $\delta.$
            
			\begin{cor}\label{FLD0} (a) $\F_0=D(0);$ (b) $\F_L=D(0+)=D((0,\tilde{\Delta}_L)).$
			\end{cor}
			\begin{proof}
				Let us prove (a). By definition, 
                \[\F_0=\{\phi\in\F:r(\phi)(x)\ge r(\psi)(x)\ \text{for all $\psi\in\F$ and for all $x\in\X$}\}.\]
                In view of~\eqref{eqTalpha}, $\phi\in\F_0$ iff $T_0^\phi V_0=T_0V_0,$ which is $\phi\in D(0)$ by Definition~\ref{defMnDnD}.
                
                Let us prove (b). By letting $i=H$ in Proposition~\ref{smallpar}(b), we have $D((0,\tilde{\Delta}_L))\subset\F_L,$ which implies $D(0+)\subset D((0,\tilde{\Delta}_L))\subset\F_L.$  On the other hand, if $\psi\in\F_L$ and $\phi\in D(0+)\subset\F_L,$ then $P^{t}(\phi)r(\phi)=P^t(\psi)r(\psi)$ for all $t\in\N$ by the definition of $\F_L.$ Then by Theorem~\ref{coreq}(a) and by Definitions~\ref{defMnDnD} and \ref{def+-}, there exists $\epsilon>0$ such that $v_\alpha^\psi=v_\alpha^\phi=V_\alpha$ for all $\alpha\in[0,\epsilon],$ which implies $\psi\in D(0+).$ Hence, $\F_L\subset D(0+),$ and therefore $\F_L=D(0+)=D((0,\tilde{\Delta}_L)).$
			\end{proof}

			The following corollary shows that, the smallest positive irregular point, if exists, is bounded below by $\tilde{\Delta}_L.$ Example~\ref{exaa1bound} shows that this bound can hold in the form of an equality.
            
			\begin{cor}\label{a1bound} $a_1\ge\tilde{\Delta}_L.$
			\end{cor}
			\begin{proof}
				This inequality follows from Corollary~\ref{FLD0}(b) and  Proposition~\ref{partition}(c).  
			\end{proof}

			The following theorem provides a uniform upper bound of the turnpike function for small discount factors. Example~\ref{exanbm} shows that two inequalities in statement (a) are sharp. Statement (b) implies that, if the terminal reward is $\pmb{0},$ then the turnpike function is equal to the constant $L+1$ for small discount factors. This is true because, if $s$ is a vector with constant coordinates, then $Ps=s$ for every stochastic matrix $P,$ and therefore the assumption in Theorem~\ref{ndboundm}(b) holds. 
            
            \begin{thm}\label{ndboundm}\begin{enumerate}[(a)]
			    \item $N^*([0,\Delta_L))\le L+1\le m;$
                    \item $N(\alpha)=L+1$ for all $\alpha\in(0,\Delta_L)$ if $\{P^t(\phi)s\}_{\phi\in\F_0}$ is a singleton for every $0\le t\le L.$ 
			\end{enumerate}
			\end{thm} 
			\begin{proof}
				Let us prove (a). Let $\alpha\in(0,\Delta_L).$ By Corollary~\ref{a1bound} all points in $(0,\tilde{\Delta}_L)$ are regular. Hence, by Theorems~\ref{utt}, \ref{smallpar}(a) and by Corollary~\ref{FLD0}(b) we have $D_n(\alpha)\subset\F_L=D((0,\tilde{\Delta}_L))=D(\alpha)$ for all $n\ge L+1,$ which implies $N(\alpha)\le L+1.$ Therefore, $N^*((0,\Delta_L))\le L+1\le m,$ where the last inequality is proven before Lemma~\ref{lempower}. We also recall that $N(0)=1$ by the definition of a turnpike integer.
				
				Let us prove (b). If $L=0$ then (b) follows from (a). Let $L\ge1$ and $\alpha\in(0,\Delta_L).$ Then $H\ge1.$ Since $(L-1)-(L-L_j-2)=L_j+1$ for all $j=0,1,\ldots,H-1,$ by Lemma~\ref{lemrl} and by Proposition~\ref{smallpar}(a) there exists $\pi=(\phi_0,\phi_1,\ldots)\in M_{L-1}(\alpha)$ such that $\{\phi_t\}_{t=0}^{L-L_j-2}\subset\F_{L_j}$ for all $j=0,1,\ldots,H-1.$ Let $\phi\in\F_{L_{H-1}}.$ We denote $\phi_{-1}=\phi$ and let $\tilde{\pi}=(\psi_0,\psi_1,\ldots)=(\phi_{-1},\phi_0,\phi_1,\ldots)=\phi^\pi.$ Then $\{\psi_t\}_{t=0}^{L-L_j-1}=\{\phi_t\}_{t=-1}^{L-L_j-2}\subset\F_{L_j}$ for all $j=0,1,\ldots,H-1.$ By Lemma~\ref{lemFpower} this implies $P_t(\phi^\pi)r_t(\phi^\pi)=P_t(\tilde{\pi})r_t(\tilde{\pi})=P^t(\phi)r(\phi)$ for all $t=0,1,\ldots,L-1.$ Also $P_L(\phi^\pi)s=P^L(\phi)s$ by $D(0)=\F_0=\F_{L_0}$ from Corollary~\ref{FLD0} and by the assumption that $\{P^t(\tilde{\phi})s\}_{\tilde{\phi}\in\F_0}$ is a singleton for every $0\le t\le L.$ Therefore, since $\phi\in\F_{L_{H-1}},$ 
				\begin{equation}\label{TLs}
					T_\alpha^\phi V_{L-1,\alpha}=v_{L,\alpha}^{\phi^\pi}=\sum_{t=0}^{L-1}\alpha^tP_t(\phi^\pi)r_t(\phi^\pi)+\alpha^LP_L(\phi^\pi)s=\sum_{t=0}^{L-1}\alpha^tP^t(\phi)r(\phi)+\alpha^LP^L(\phi)s.
				\end{equation} 
                
    			Now let $\phi\in D_L(\alpha)$ and $\psi\in\F_{L_{H-1}}\setminus\F_L.$ Since $L\ge L_{H-1}+1$ and $\Delta_L\le \Delta_{L_{H-1}},$ Proposition~\ref{smallpar}(a) implies $\phi\in\F_{L_{H-1}}.$ Then equation~\eqref{TLs}, the assumption that $\{P^t(\tilde{\phi})s\}_{\tilde{\phi}\in\F_0}$ is a singleton for every $0\le t\le L,$ the definition of $\F_{L_{H-1}}$ and the definition of $D_L(\alpha)$ imply $T_\alpha^\phi V_{L-1,\alpha}=T_\alpha^\psi V_{L-1,\alpha}=V_{L,\alpha}.$ Hence, $\psi\in D_L(\alpha),$ but by Proposition~\ref{smallpar}(b) we have $\psi\notin D(\alpha)$ because $0<\alpha<\Delta_{L}\le\tilde{\Delta}_L$ and $\psi\notin\F_L.$ Therefore, $D_L(\alpha)\not\subset D(\alpha),$ which implies $N(\alpha)\ge L+1.$ And by (a) we conclude that $N(\alpha)=L+1.$
			\end{proof}

			In view of Theorem~\ref{ndboundm}(a), $\Delta_L$ can be the biggest value of $\delta>0$ such that $N^*((0,\delta))<\infty;$ see Example~\ref{exaa1bound}.

            For each $x\in\X$ let $r_0(x)=r(x):=\max_{a\in A(x)}r(x,a),$ $A_0(x):=\argmax_{a\in A(x)}r(x,a),$ and
    \begin{equation}\label{defaltern1}
        r_n(x):=\max_{a\in A_{n-1}(x)}\sum_{y\in\X}p(y|x,a)r_{n-1}(y),\quad  A_n(x):=\argmax_{a\in A_{n-1}(x)}\sum_{y\in\X}p(y|x,a)r_{n-1}(y),\quad n\in\N^+.
    \end{equation}
     Then alternatively, the sets $\{\F_n\}_{n\ge0},$ $\{\X_n\}_{n\ge0},$ and $\{C_n\}_{n\ge0}$ can be written as
     \begin{equation}\label{FXalt}
        \F_n=\left\{\phi\in\F:\phi(x)\in A_n(x)\ \text{for all}\ x\in\X\right\},\quad    \X_n=\left\{x\in\X:A_n(x)\ne A_{n-1}(x)\right\},
    \end{equation}
        \begin{equation}\label{defaltern2}
        C_n=\begin{cases}
        C_{n-1} &\text{if\ \ $\X_n=\emptyset$},\\
    \displaystyle\min\left\{C_{n-1}, \min_{x\in\X_n}\{r_n(x)-\max_{a\in A_{n-1}(x)\setminus A_n(x)}\sum_{y\in\X}p(y|x,a)r_{n-1}(y)\}\right\} &\text{if\ \ $\X_n\ne\emptyset$}.
    \end{cases}
    \end{equation} 
    We denote by $r_n$ the corresponding vectors in $\R^m$ and view $A_n(\cdot):\X\to\A$ as set-valued functions. The following algorithm computes $\{A_n\}_{n=0}^L,$ which define $\{\F_n\}_{n=0}^L,$ $L,$ $H$, $\{L_i\}_{i=1}^H,$ $\{C_t\}_{t=0}^L,$ $\{\Delta_t\}_{t=0}^L,$ and $\{\tilde{\Delta}_t\}_{t=0}^L.$ In view of Theorem~\ref{ndboundm}(a), the algorithm conducts  at most $m$ iterations. For each iteration the number of arithmetic operations is at most $O(mq).$ Thus, the total number of arithmetic operations is $O(m^2q).$ We notice that, in view of Corollary~\ref{FLD0}(b), every deterministic policy $\phi\in\F_L,$ where $\phi(x)\in A_L(x)$ for all $x\in\X,$ is infinite-horizon optimal for all discount factors $\alpha\in [0,\tilde{\Delta}_L],$ where $L,$ $\tilde{\Delta}_L,$ and $A_L(x),$ $x\in\X,$ are outputs of the algorithm.

\begin{algorithm}[H]
\SetAlgoLined
\LinesNumbered
\caption{Computes $L,$ $H,$ $\{A_n\}_{n=0}^L,$ $\{L_i\}_{i=1}^H,$ $\{C_t\}_{t=0}^L,$ $\{\Delta_t\}_{t=0}^L,$ $\{\tilde{\Delta}_t\}_{t=0}^L.$}\label{algsmall}
\KwIn{$\X,$ $\{A(x)\}_{x\in\X},$ $\{r(x,a)\}_{x\in\X,\ a\in A(x)},$ $\{p(y|x,a)\}_{x,y\in\X,\ a\in A(x)},$ $\{s(x)\}_{x\in\X}.$}
\KwOut{$L,$ $H,$ $\{A_n\}_{n=0}^L,$ $\{L_i\}_{i=1}^H,$ $\{C_t\}_{t=0}^L,$ $\{\Delta_t\}_{t=0}^L,$ $\{\tilde{\Delta}_t\}_{t=0}^L.$}
Set $n=1$ and $H=L=L_0=0.$ Compute $R^*,$ $R_1^*,$ $r_0,$ $A_0,$ $C_0,$ $\Delta_0,$ and $\tilde{\Delta}_0$\;
 \While{$n<m$ and there exists some $x\in\X$ such that $A_{n-1}(x)$ is not a singleton}{
 \For{each $x\in\X$}{
 $r_n(x)\leftarrow\max_{a\in A_{n-1}(x)}\sum_{y\in\X}p(y|x,a)r_{n-1}(y)$\;
 $A_n(x)\leftarrow\argmax_{a\in A_{n-1}(x)}\sum_{y\in\X}p(y|x,a)r_{n-1}(y)$\;}
 $\X_n\leftarrow\{x\in\X:A_n(x)\ne A_{n-1}(x)\}$\;
  \uIf{$\X_n=\emptyset$}{
   $C_n\leftarrow C_{n-1}$\;
   }
   \Else{\lFor{each $x\in\X_n$}
   {$C^{x}\leftarrow\max_{a\in A_{n-1}(x)\setminus A_n(x)}\sum_{y\in\X}p(y|x,a)r_{n-1}(y)$}
   $C\leftarrow\min_{x\in\X_n}\left\{r_n(x)-C^{x}\right\},$ $C_n\leftarrow\min\{C_{n-1},C\},$ $H\leftarrow H+1,$ $L_H\leftarrow n$\;}
  \leIf{$C_n=\infty$}{
   $\Delta_n\leftarrow1,$ $\tilde{\Delta}_n\leftarrow1$
   }{$\Delta_n\leftarrow\frac{C_n}{2R^*+C_n},$ $\tilde{\Delta}_n\leftarrow\frac{C_n}{2R_1^*+C_n}$
  }
  $n\leftarrow n+1,$ $L\leftarrow L_H$\;
 }
\end{algorithm}

        Recall from Corollary~\ref{0regular} that, if $0$ is a regular point then there exist some $\delta>0$ such that $N(\alpha)=1$ for all $\alpha\in[0,\delta).$ The following corollary provides a concrete value of such a $\delta$.
			
			\begin{cor}\label{0regularb} If $0$ is a regular point, then $N^*([0,\Delta_0))=1.$ 
			\end{cor}
			\begin{proof}
				If $0$ is a regular point, then $D(0)=D(0+)$ by Definition~\ref{defirr}, which implies $L=0$ by Corollary~\ref{FLD0}. The result follows from Theorem~\ref{ndboundm}(a).
			\end{proof}

   The following corollary states that, if an MDP is deterministic and the maximum one-step rewards of all the states are distinct, then for small discount factors the value iteration algorithm converges within at most two iterations.
			
			\begin{cor} If an MDP is deterministic and $r(x)\ne r(z)$ for all $x,z\in\X$ such that $x\ne z,$ then $L\le1$ and $N^*([0,\Delta_L))\le2.$ 
			\end{cor}
			\begin{proof}
				If $m=1,2$ then the result is by the fact that $L\le m-1$ and by Theorem~\ref{ndboundm}(a). We assume $m\ge3$. Suppose there are some $\phi\in\F_L$ and $\psi\in\F_0\setminus\F_L$ such that $P(\phi)r(\phi)=P(\psi)r(\psi).$ Since $\phi,\psi\in\F_0,$ we have $r(\phi)=r(\psi)=r.$ Since the MDP is deterministic and $\phi,\psi$ are different, there exist some $1\le i\le m$ and some $1\le j_1,j_2\le m$ with $j_1\ne j_2$ such that $P_{i,j_1}(\phi)=P_{i,j_2}(\psi)=1$. Then $r(x^{j_1})=r(x^{j_2})$ from the $i$-th row of $P(\phi)r(\phi)=P(\psi)r(\psi),$ which contradicts the assumption that $\{r(x^i)\}_{i=1}^m$ are distinct. This implies $L\le1,$ and the result follows from Theorem~\ref{ndboundm}(a).
			\end{proof}

			The following corollary strengthens Proposition~\ref{uttirr}(a).
            
			\begin{cor}\label{uttirrno0} Let $\mathcal{I}:=[a,b]\subset [0,1)$ such that $\mathcal{I}$ contains at least one irregular point, that is, $\P(\mathcal{I})\ne\emptyset.$ Then
				\begin{enumerate}[(a)]
					\item if $
					\P(\mathcal{I})=\{0\}$, then there exists $K<\infty$ such that $D_n(\mathcal{I})\subset D(\P(\mathcal{I}))$ for all $n\ge K$;
					\item if $\P(\mathcal{I})\setminus\{0\}\ne\emptyset$, then there exists $K<\infty$ such that 
					$D_n(\mathcal{I})\subset D(\P(\mathcal{I})\setminus\{0\})$ for all $n\ge K.$
				\end{enumerate}
			\end{cor}
			
			\begin{proof}
				(a) follows from Proposition~\ref{uttirr}(a). Let us prove (b). If $0\notin\P(\mathcal{I}),$ then this follows from Proposition~\ref{uttirr}(a). Let $0\in\P(\mathcal{I})$. Note that $a_1\in\P(\mathcal{I})\setminus\{0\}.$ Then $D_n((0,\Delta_L))=\F_L=D(0+)=D(a_1-)\subset D(a_1)$ for all $n\ge m$ by Corollary~\ref{smallpar}(a), by Theorem~\ref{ndboundm}(a), and by Corollary~\ref{FLD0}(b). Therefore, (b) is true in view of Proposition~\ref{uttirr}(a). 
			\end{proof}

			The following corollary strengthens Corollary~\ref{ndlimit}.
            
			\begin{cor}\label{ndlimitg}
				If $\alpha$ is a limit point of $\mathbb{D}([0,1))$, then either $\alpha\in\P((0,1))$ or $\alpha=1.$
			\end{cor}
			\begin{proof}
				Theorem~\ref{ndboundm}(a) and Theorem~\ref{finitetpi} imply $\alpha\ne0.$ The result follows from Corollary~\ref{ndlimit}.
			\end{proof}

\section{Examples}\label{secexamples}
In this section, we present 
examples showing different kinds of behaviors of the set-valued functions $M_n(\alpha),$ $D_n(\alpha),$ $D(\alpha)$ and the turnpike function $N(\alpha)$ that are discussed in this paper. All examples in this section present deterministic MDPs. We use arrows to indicate the succeeding state when an action is taken at a state, and the numbers on arrows are rewards of corresponding actions. In each example, either the terminal reward vector $s={\bf 0}$ or it is defined in the corresponding example. 


\begin{exa}\label{exaDn+-} \emph{\textbf{\boldmath{For an $n$-horizon-first-step break point $\alpha$, while $M_n(\alpha-)\cap M_n(\alpha+)=\emptyset,$ it is possible that $D_n(\alpha-)\cap D_n(\alpha+)\ne\emptyset$. Also, in this example, the converse of Proposition~\ref{Nd}(c) is not true.}} Consider the following MDP; see Figure~\ref{fig:exaDn+-}. There are in total four different decision rules, denoted by $\phi^{(1)},\phi^{(2)},\phi^{(3)},\phi^{(4)}$, which differ by their stochastic transition matrices
\begin{figure}[ht!]
	\centerline{\includegraphics[scale=1.59]{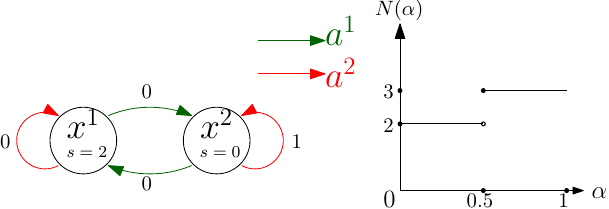}}
	\caption{MDP model and graph of the turnpike function $N(\alpha)$ for Example~\ref{exaDn+-}.}
	\label{fig:exaDn+-}
\end{figure}
\begin{equation*}\begin{aligned}
	P(\phi^{(1)})&=\begin{bmatrix}
		1 & 0 \\
		1 & 0
	\end{bmatrix},\quad P(\phi^{(2)})=\begin{bmatrix}
		1 & 0 \\
		0 & 1
	\end{bmatrix},\quad P(\phi^{(3)})=\begin{bmatrix}
		0 & 1 \\
		1 & 0
	\end{bmatrix},\quad P(\phi^{(4)})=\begin{bmatrix}
		0 & 1 \\
		0 & 1
	\end{bmatrix},
\end{aligned}\end{equation*}
where the first and second rows represent states $x^1$ and $x^2$ respectively. Rewards and terminal rewards are
\begin{equation*}\begin{aligned}
	r(\phi^{(1)})&=\begin{bmatrix}
		0 \\ 0
	\end{bmatrix},\qquad r(\phi^{(2)})=\begin{bmatrix}
		0 \\ 1
	\end{bmatrix},\qquad r(\phi^{(3)})=\begin{bmatrix}
		0 \\ 0
	\end{bmatrix},\qquad r(\phi^{(4)})=\begin{bmatrix}
		0 \\ 1
	\end{bmatrix},\qquad s=\begin{bmatrix}
		2 \\ 0
	\end{bmatrix}.
\end{aligned}\end{equation*}
It can be shown by calculations that
\begin{equation*}
    \begin{cases}
    \displaystyle V_{1,\alpha}=T_\alpha^{\phi^{(2)}}s=\begin{bmatrix}
			2\alpha \\ 1
		\end{bmatrix},\quad V_{2,\alpha}=T^{\phi^{(4)}}_\alpha\begin{bmatrix}
			2\alpha \\ 1
		\end{bmatrix}=\begin{bmatrix}
			\alpha \\ 1+\alpha
		\end{bmatrix}, &\text{if $0\le\alpha<0.5,$}\\
  \\
   \displaystyle V_{1,\alpha}=T_\alpha^{\phi^{(1)}}s=\begin{bmatrix}
			2\alpha \\ 2\alpha
		\end{bmatrix},\quad V_{2,\alpha}=T_\alpha^{\phi^{(2)}}\begin{bmatrix}
			2\alpha \\ 2\alpha
		\end{bmatrix}=T_\alpha^{\phi^{(4)}}\begin{bmatrix}
			2\alpha \\ 2\alpha
		\end{bmatrix}=\begin{bmatrix}
			2\alpha^2 \\ 1+2\alpha^2
		\end{bmatrix}, &\text{if $0.5\le\alpha<1.$}
    \end{cases}
\end{equation*}
We conclude that $D_2(0.5-)=\{\phi^{(4)}\}$ and $D_2(0.5)=D_2(0.5+)=\{\phi^{(2)},\phi^{(4)}\}$. Since $D_2(0.5-)\ne D_2(0.5+)$, then $\alpha=0.5$ is a $2$-horizon-first-step break point. However, $D_2(0.5-)\cap D_2(0.5+)=\{\phi^{(4)}\}\ne\emptyset$; $M_2(0.5-)=\{(\phi^{(4)},\phi^{(2)},\ldots)\}$, $M_2(0.5+)=\{(\phi^{(4)},\phi^{(1)},\ldots),(\phi^{(2)},\phi^{(1)},\ldots)\},$ and thus $M_2(0.5-)\cap M_2(0.5+)=\emptyset$. Furthermore, it can be shown by induction that for $n\ge 3$
\begin{equation*}
    \begin{cases}
    V_\alpha=v_\alpha^{\phi^{(4)}}=\begin{bmatrix}
			\displaystyle \frac{\alpha}{1-\alpha} \\ \\ \displaystyle \frac{1}{1-\alpha}
		\end{bmatrix},\quad V_{n,\alpha}=T_\alpha^{\phi^{(4)}}V_{n-1,\alpha}=\begin{bmatrix}
			\displaystyle\frac{\alpha-\alpha^n}{1-\alpha} \\ \\ \displaystyle\frac{1-\alpha^n}{1-\alpha}
		\end{bmatrix},\ & \text{if $0\le\alpha<0.5,$}\\
   V_\alpha=v_\alpha^{\phi^{(4)}}=\begin{bmatrix}
			\displaystyle \frac{\alpha}{1-\alpha} \\ \\ \displaystyle \frac{1}{1-\alpha}
		\end{bmatrix},\quad V_{n,\alpha}=T_\alpha^{\phi^{(4)}}V_{n-1,\alpha}=\begin{bmatrix}
			\displaystyle\frac{\alpha-\alpha^{n-1}}{1-\alpha}+2\alpha^n \\ \\ \displaystyle\frac{1-\alpha^{n-1}}{1-\alpha}+2\alpha^n
		\end{bmatrix},\ & \text{if $0.5\le\alpha<1$}.
    \end{cases}
\end{equation*}
Hence, $D_n([0,1))=D([0,1))=\{\phi^{(4)}\}$ for $n\ge3,$ $N(\alpha)=2$ for $\alpha\in[0,0.5),$ and $N(\alpha)=3$ for $\alpha\in[0.5,1);$  see Figure~\ref{fig:exaDn+-} for the graph of $N(\alpha$). We see that $(0,1)$ is a partition interval, $N(0.5)=3,$ $\phi^{(4)}\in D((0,1)),$ and $\phi^{(4)}\in D_2(\alpha)=D_{N(0.5)-1}(\alpha)$ for all $\alpha\in[0.5,1).$ However, $0.5\notin \mathbb{D}^+((0,1)).$ Therefore, the converse of Proposition~\ref{Nd}(c) is not true.}
\end{exa}

\begin{exa}\label{exaisolate} \emph{\textbf{\boldmath{A turnpike function $N(\alpha)$ can be neither left-continuous nor right-continuous at some regular point in $(0,1)$.}} Consider the following MDP; see Figure~\ref{fig:exaisolate}.
\begin{figure}[ht!]
	\centerline{\includegraphics[scale=1.13]{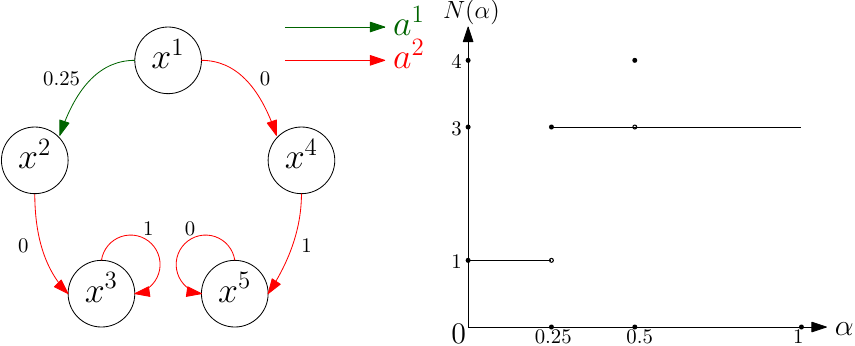}}
	\caption{MDP model and graph of the turnpike function $N(\alpha)$ for Example~\ref{exaisolate}.}
	\label{fig:exaisolate}
\end{figure}
There are only two decision rules $\phi,\psi$ which only differ at state $x^1,$ and $\phi(x^1)=a^1,\psi(x^1)=a^2$. Thus, we only need to consider the objective functions at state $x^1$. We have
\begin{equation*}
	v_\alpha^{\phi}(x^1) = 0.25+\sum_{t=2}^\infty\alpha^t= 0.25+\frac{\alpha^2}{1-\alpha};\ \ v_\alpha^{\psi}(x^1) =\alpha.
\end{equation*}
Since $v_\alpha^{\phi}(x^1)\ge0.25+\alpha^2\ge\alpha=v_\alpha^{\psi}(x^1),$ we have $v_\alpha^{\phi}(x^1)>v_\alpha^{\psi}(x^1)$ for all $\alpha\in[0,1).$  Therefore, $D([0,1))=\{\phi\}$. 
For $n=1$ we have $T_\alpha^{\phi}V_{0,\alpha}(x^1)-T_\alpha^{\psi}V_{0,\alpha}(x^1)=0.25>0;$ for $n=2$ we have $T_\alpha^{\phi}V_{1,\alpha}(x^1)-T_\alpha^{\psi}V_{1,\alpha}(x^1)=0.25-\alpha$; for $n=3$ we have $T_\alpha^{\phi}V_{2,\alpha}(x^1)-T_\alpha^{\psi}V_{2,\alpha}(x^1)=(0.25+\alpha^2)-\alpha=(\alpha-0.5)^2\ge 0,$
and the equality holds only at $\alpha=0.25$; for $n\ge4$ we have $T_\alpha^{\phi}V_{n-1,\alpha}(x^1)-T_\alpha^{\psi}V_{n-1,\alpha}(x^1)=(0.25+\sum_{t=2}^n\alpha^t)-\alpha=(\alpha-0.5)^2+\sum_{t=3}^n\alpha^t>0;$ see Figure~\ref{fig:exaisolate} for
the graph of $N(\alpha).$ Theorems~\ref{Ndiff}$-$\ref{ndclass} can be verified at points $\alpha=0.25$ and $\alpha=0.5,$ where $N(\alpha)$ is discontinuous. Also, $\alpha=0.5$ is a point where $N(\alpha)$ is neither left-continuous nor right-continuous. This does not violate the conclusions of Theorem~\ref{finitetpi} since $\{0.5\}=[0.5,0.5]$ is also an interval.}
\end{exa}

\begin{exa}\label{exaub1} \emph{\textbf{\boldmath{A turnpike function $N(\alpha)$ can be unbounded near $\alpha=1.$}} Consider the following MDP; see Figure~\ref{fig:exaub1}.
\begin{figure}[ht!]
	\centerline{\includegraphics[scale=1.2]{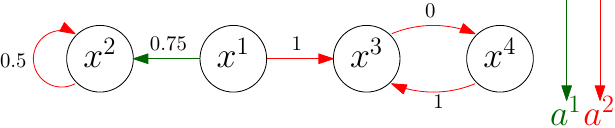}}
	\caption{MDP model for Example~\ref{exaub1}.}
	\label{fig:exaub1}
\end{figure}
There are only two decision rules $\phi,$ $\psi$ which only differ at $x^1$ such that $\phi(x^1)=a^1,\psi(x^1)=a^2$. We have
\[v_\alpha^{\phi}(x^1)-v_\alpha^{\psi}(x^1)=\left(0.75+0.5\sum_{t=1}^\infty\alpha^t\right)-\sum_{t=0}^\infty\alpha^{2t}=\frac{\alpha-1}{4(1+\alpha)}<0,\ \ \text{for $\alpha\in[0,1)$}.\]
Hence, $D([0,1))=\{\psi\}$. For $n\in\N^+$ and for $\alpha\in (0,1),$
\begin{equation*}
\begin{aligned}
    &T_\alpha^{\phi}V_{2n-1,\alpha}(x^1)-T_\alpha^{\psi}V_{2n-1,\alpha}(x^1)=\left(0.75+0.5\sum_{t=1}^{2n-1}\alpha^t\right)-\sum_{t=1}^{n}\alpha^{2t-2}=\frac{2\alpha^{2n}+\alpha-1}{4(1+\alpha)},\\
    &T_\alpha^{\phi}V_{2n,\alpha}(x^1)-T_\alpha^{\psi}V_{2n,\alpha}(x^1)=\left(0.75+0.5\sum_{t=1}^{2n}\alpha^t\right)-\sum_{t=1}^{n}\alpha^{2t}=\frac{\alpha-1-2\alpha^{2n+1}}{4(1+\alpha)}<0. 
\end{aligned}
\end{equation*}
Let $\alpha(n)$ be the unique real root of $2\alpha^{2n}+\alpha-1=0$ in $(0,1)$ for  $n\in\N^+.$ Then
\[2\left[\alpha(n)\right]^{2(n+1)}+\alpha(n)-1<2\left[\alpha(n)\right]^{2n}+\alpha(n)-1=0=2\left[\alpha(n+1)\right]^{2(n+1)}+\alpha(n+1)-1,\]
which implies $\alpha(n)<\alpha(n+1)$ and $N(\alpha)=2n+1$ for all $n\in\N^+,\alpha\in[\alpha(n),\alpha(n+1)).$ By Corollary~\ref{ndlimit} $\alpha(n)\to1$ as $n\to\infty$. Therefore, $N(\alpha)$ is unbounded near $\alpha=1.$}
\end{exa}

\begin{exa}\label{exalirb}
    \emph{\textbf{\boldmath{A turnpike function can be unbounded on either side of $\alpha$ and be bounded on the other side of $\alpha,$ where $\alpha$ is an irregular point.}} Consider the following MDP; see Figure~\ref{fig:exalirb}. There are only two decision rules $\phi,\psi$ which only differ at $x^1,$ and $\phi(x^1)=a^1,\psi(x^1)=a^2$.
\begin{figure}[ht!]
	\centerline{\includegraphics[scale=1.2]{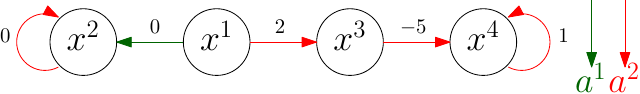}}
	\caption{MDP model for Example~\ref{exalirb}.}
	\label{fig:exalirb}
\end{figure}
We have
    \[v_\alpha^\phi(x^1)-v_\alpha^\psi(x^1)=0-\left(2-5\alpha+\sum_{t=2}^\infty\alpha^t\right)=-\frac{(2\alpha-1)(3\alpha-2)}{1-\alpha}.\]
Hence, $D([0,\frac{1}{2}))=D((\frac{2}{3},1))=\left\{\psi\right\},$ $D((\frac{1}{2},\frac{2}{3}))=\left\{\phi\right\},$ and $D(\frac{1}{2})=D(\frac{2}{3})=\left\{\phi,\psi\right\}.$ There are two irregular points $\frac{1}{2}$ and $\frac{2}{3}$ which are both break points. Also,
\begin{equation}\label{eqlr}
    T_\alpha^\phi V_{n-1,\alpha}(x^1)-T_\alpha^\psi V_{n-1,\alpha}(x^1)=0-\left(2-5\alpha+\sum_{t=2}^{n-1}\alpha^t\right)=\frac{\alpha^n-(2\alpha-1)(3\alpha-2)}{1-\alpha},\qquad n\ge2.
\end{equation}
Let $f_n(\alpha):=\alpha^n-(2\alpha-1)(3\alpha-2)$ for each $n\ge2.$ Then $f_n(0)=-2<0,$ $f_n(\frac{1}{2})=\frac{1}{2^n}>0,$ $f_n(\frac{2}{3})=(\frac{2}{3})^n>0,$ and $f_n(1)=0.$ The first graph of Figure~\ref{fig:exalirbb} shows the graph of $f_n(\alpha).$ Thus, there exists $\alpha_n\in(0,\frac{1}{2})$ and $\beta_n\in(\frac{2}{3},1)$ such that $f_n(\alpha_n)=f_n(\beta_n)=0,$ which implies $T_\alpha^\phi V_{n-1,\alpha}(x^1)=T_\alpha^\psi V_{n-1,\alpha}(x^1)$ for $\alpha=\frac{1}{2},\frac{2}{3}.$ This means $D_n(\alpha_n)=D_n(\beta_n)=\{\phi,\psi\}\not\subset D(\alpha_n)=D(\beta_n)=\{\psi\},$ and therefore $N(\alpha_n)>n$ and $N(\beta_n)>n.$ Theorem~\ref{utt} implies that $\alpha_n\to\frac{1}{2}$ and $\beta_n\to\frac{2}{3}.$ 
Moreover, for $n\ge2,$ note that $T_\alpha^\phi V_{n-1,\alpha}(x^1)-T_\alpha^\psi V_{n-1,\alpha}(x^1)=v_\alpha^\phi(x^1)-v_\alpha^\psi(x^1)+\frac{\alpha^n}{1-\alpha},$ and hence $T_\alpha^\phi V_{n-1,\alpha}(x^1)>T_\alpha^\psi V_{n-1,\alpha}(x^1)$ for $\alpha\in(\frac{1}{2},\frac{2}{3}).$ So, $D_n((\frac{1}{2},\frac{2}{3}))=\{\phi\}=D((\frac{1}{2},\frac{2}{3}))$ for $n\ge2,$ which implies $N^*((\frac{1}{2},\frac{2}{3}))\le2.$ The second graph in Figure~\ref{fig:exalirbb} is the actual graph of $N(\alpha).$
\begin{figure}[ht!]
	\centerline{\includegraphics[scale=0.81]{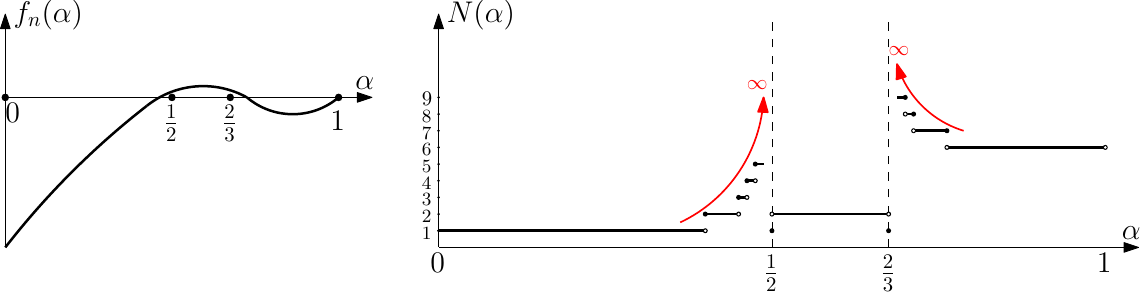}}
	\caption{Graph of the turnpike function $N(\alpha)$ for Example~\ref{exalirb}.}
	\label{fig:exalirbb}
\end{figure}
}
\end{exa}

To verify the next two examples, we first introduce the following two propositions. Recall that simple break points are defined in Definition~\ref{defsbp}. If $\alpha$ is a simple break point, and Condition \hyperref[conA]{$\mathrm{A}$} holds at $\alpha,$ the following proposition shows that this condition can be verified within a finite number of computations. 

\begin{pro}\label{propoA} Let $\alpha\in\P((0,1))$ and both $D(\alpha-)=\{\phi\},D(\alpha+)=\{\psi\}$ are singletons. 
\begin{enumerate}[(a)]
\item If both Conditions \hyperref[thm:Ca-]{$\mathrm{A_-}$} and \hyperref[thm:Ca+]{$\mathrm{A_+}$} hold at $\alpha$, then 
     there exists $K\ge N(\alpha)-1$ such that for all $t\in\N$
     \begin{equation} \label{march17} P^t(\phi)(V_{\alpha}-V_{K,\alpha})=P^t(\psi)(V_{\alpha}-V_{K,\alpha}).\end{equation}
     
    \item  If there exists $K\ge N(\alpha)-1$ such that equality \eqref{march17} holds for all $t=0,1,\ldots,\min\{G(\phi,V_{\alpha}-V_{K,\alpha}),G(\psi,V_{\alpha}-V_{K,\alpha})\}-1,$ then \eqref{march17} holds for all $t\in\N,$
 and, in addition, if  $\alpha$ is a simple break point, that is, $D(\alpha)=\{\phi,\psi\},$ then Condition \hyperref[conA]{$\mathrm{A}$} holds at $\alpha.$
\end{enumerate}
\end{pro}
\begin{proof}
    (a). Since $D(\alpha-)=\{\phi\}$ and $D(\alpha+)=\{\psi\},$ both Conditions \hyperref[thm:Ca-]{$\mathrm{A_-}$} and \hyperref[thm:Ca+]{$\mathrm{A_+}$} hold at $\alpha$ if and only if there exists $\tilde{K}\in\N^+$ such that $D_n(\alpha)=\{\phi,\psi\}$ for $n\ge \tilde{K}.$ By \eqref{eqoptimal}, this means
    \begin{align}\label{Vconvergeeq}
    \left(T_{\alpha}^\phi\right)^{n-\tilde{K}+1}V_{\tilde{K}-1,\alpha}=\left(T_{\alpha}^\psi\right)^{n-\tilde{K}+1}V_{\tilde{K}-1,\alpha}=V_{n,\alpha},\qquad
    n\ge\tilde{K}-1.
\end{align}
    By 
    \eqref{eqvaluesum1} and \eqref{eqoptimal},
\begin{equation}
\begin{aligned}
    \left(T_{\alpha}^\phi\right)^{n-\tilde{K}+1}&V_{\tilde{K}-1,\alpha}=\sum_{t=0}^{n-\tilde{K}}\alpha^tP^t(\phi)r(\phi)+\alpha^{n-\tilde{K}+1}P^{n-\tilde{K}+1}(\phi)V_{\tilde{K}-1,\alpha}\\
	&=\left[I-\alpha^{n-\tilde{K}+1}P^{n-\tilde{K}+1}(\phi)\right]\left[I-\alpha P(\phi)\right]^{-1}r(\phi)+\alpha^{n-\tilde{K}+1}P^{n-\tilde{K}+1}(\phi)V_{\tilde{K}-1,\alpha}\\
	& = \label{Vconverge} V_{\alpha}-\alpha^{n-\tilde{K}+1}P^{n-\tilde{K}+1}(\phi)\left(V_{\alpha}-V_{\tilde{K}-1,\alpha}\right),\qquad
    n\ge \tilde{K}-1.
\end{aligned}
\end{equation}
where $V_{\alpha}=v_{\alpha}^\phi=[I-\alpha P(\phi)]^{-1}r(\phi)$ follows from \eqref{eqvaluesum3} and $\phi\in D(\alpha).$ Similarly, equations \eqref{Vconverge} hold if $\phi$ is replaced by $\psi,$ which is
\begin{equation}\label{Vconvergeend}
	\left(T_{\alpha}^\psi\right)^{n-\tilde{K}+1}V_{\tilde{K}-1,\alpha}=V_{\alpha}-\alpha^{n-\tilde{K}+1}P^{n-\tilde{K}+1}(\psi)\left(V_{\alpha}-V_{\tilde{K}-1,\alpha}\right),\qquad
    n\ge \tilde{K}-1.
\end{equation}
Since $\alpha\ne0$, equations \eqref{Vconvergeeq}$-$\eqref{Vconvergeend} imply \eqref{march17} with $K=\tilde{K}-1.$ Statement (a) is proved.

(b). The first part follows from Proposition~\ref{thmla}. Suppose, in addition, $D(\alpha)=\{\phi,\psi\}.$ Let $K\ge N(\alpha)-1$ such that \eqref{march17} holds and let $\tilde{K}=K+1.$ Then equations \eqref{Vconverge} and \eqref{Vconvergeend} hold, which implies the first equality in \eqref{Vconvergeeq}. The second equality in \eqref{Vconvergeeq} also holds by
%
%
%
%
$D(\alpha)=\{\phi,\psi\},$ by $\tilde{K}\ge N(\alpha),$ and by the definition of $N(\alpha).$ 
Therefore, \eqref{Vconvergeeq} holds, which implies $D_n(\alpha)=\{\phi,\psi\}=D(\alpha)$ for all $n\ge\tilde{K},$ and (b) is proved.
\end{proof}


\begin{pro}\label{propoB} Let $\alpha\in\P((0,1)).$ Then 
\begin{enumerate}[(a)]
    \item Condition \hyperref[thm:BB-]{$\mathrm{B_-}$} is violated if there exist $\phi\in D(\alpha-)$ and $\psi\in D(\alpha)\setminus D(\alpha-)$ such that $\frac{d}{d\beta}v_\beta^{\phi}\big|_{\beta=\alpha}=\frac{d}{d\beta}v_\beta^{\psi}\big|_{\beta=\alpha};$
    \item Condition \hyperref[thm:BB+]{$\mathrm{B_+}$} is violated if there exist $\phi\in D(\alpha+)$ and $\psi\in D(\alpha)\setminus D(\alpha+)$ such that $\frac{d}{d\beta}v_\beta^{\phi}\big|_{\beta=\alpha}=\frac{d}{d\beta}v_\beta^{\psi}\big|_{\beta=\alpha};$
    \item if the function $V_\alpha$ is differentiable at $\alpha,$ then both Conditions \hyperref[thm:BB-]{$\mathrm{B_-}$} and \hyperref[thm:BB+]{$\mathrm{B_+}$} are violated.
\end{enumerate}
\end{pro}
\begin{proof}
    Let us prove (a). Suppose there exist $\phi\in D(\alpha-)$ and $\psi\in D(\alpha)\setminus D(\alpha-)$ such that $\psi\in D(\alpha)\setminus D(\alpha-)$ and $\frac{d}{d\beta}v_\beta^{\phi}\big|_{\beta=\alpha}=\frac{d}{d\beta}v_\beta^{\psi}\big|_{\beta=\alpha}.$ Let $\pi$ be the deterministic policy $\psi.$ Then for $\beta\in[0,1)$
    \begin{equation*}
    \begin{aligned}
        v_\beta^{\phi^\pi}-v_\beta^{\psi^\pi}=r(\phi)+\beta P(\phi)v_\beta^{\psi}-v_\beta^\psi&=r(\phi)+\beta P(\phi)v_\beta^{\phi}+\beta P(\phi)(v_\beta^{\psi}-v_\beta^{\phi})-v_\beta^\psi\\
        &=v_\beta^\phi-v_\beta^\psi+\beta P(\phi)(v_\beta^{\psi}-v_\beta^{\phi})=\left[I-\beta P(\phi)\right](v_\beta^{\phi}-v_\beta^{\psi}),
    \end{aligned}
    \end{equation*}
    where the first and the third equalities follow from  \eqref{eqvaluesum3}. Hence,
    \[\frac{d}{d\beta}\left(v_\beta^{\phi^\pi}-v_\beta^{\psi^\pi}\right)\Big|_{\beta=\alpha}=- P(\phi)(v_{\alpha}^{\phi}-v_{\alpha}^{\psi})+\left[I-\alpha P(\phi)\right]\frac{d}{d\beta}(v_\beta^{\phi}-v_\beta^{\psi})\big|_{\beta=\alpha}=\pmb{0},\]
    where $v_{\alpha}^{\phi}-v_{\alpha}^{\psi}=\pmb{0}$ since $\{\phi,\psi\}\subset D(\alpha)$ and $\frac{d}{d\beta}(v_\beta^{\phi}-v_\beta^{\psi})\big|_{\beta=\alpha}=\pmb{0}$ by assumption. Therefore,
    \[\sup_{\pi\in\Pi(D(\alpha))}\frac{d}{d\beta}\left(v_\beta^{\phi^\pi}-v_\beta^{\psi^\pi}\right)\Big|_{\beta=\alpha}(x)\ge0,\qquad x\in\X.\]
    Condition \hyperref[thm:BB-]{$\mathrm{B_-}$} is violated at $\alpha.$ This proves (a). The proof of (b) is similar.
    
    Let us prove (c). If $\alpha$ is a break point, then the differentiability of $V_\alpha$ at $\alpha$ 
    implies $\frac{d}{d\beta}v_\beta^{\phi}\big|_{\beta=\alpha}=\frac{d}{d\beta}v_\beta^{\psi}\big|_{\beta=\alpha}$ for all $\phi\in D(\alpha-)\subset D(\alpha)\setminus D(\alpha+)$ and $\psi\in D(\alpha+)\subset D(\alpha)\setminus D(\alpha-).$ If $\alpha\in\P((0,1))$ is not a break point, then it must be a touching point with $D(\alpha-)=D(\alpha+).$ The $\frac{d}{d\beta}v_\beta^{\phi}\big|_{\beta=\alpha}=\frac{d}{d\beta}v_\beta^{\psi}\big|_{\beta=\alpha}$ for all $\phi\in D(\alpha-)=D(\alpha+)$ and all $\psi\in D(\alpha)\setminus D(\alpha-)=D(\alpha)\setminus D(\alpha+),$ because otherwise there exists $x\in\X$ such that $v_\beta^\phi(x)<v_\beta^\psi(x)$ for some $\beta$ close to $\alpha,$ contradicting to $\phi\in D(\alpha-)=D(\alpha+).$ In both cases, (c) follows from (a) and (b).
\end{proof}

\begin{exa}\label{exatbbounded} \emph{\textbf{\boldmath{The assumptions of Theorem~\ref{tbnecessary} (and therefore \cite[Theorems 3]{doi:10.1287/moor.2017.0912}) are not sufficient to guarantee the boundedness of a turnpike function at either side of a break point.}} Consider the following MDP; see Figure~\ref{fig:exatbbounded}.
\begin{figure}[ht!]
	\centerline{\includegraphics[scale=1.2]{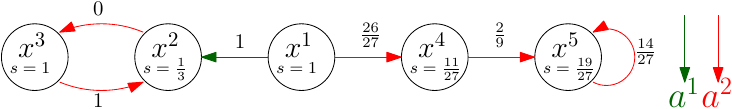}}
	\caption{MDP model for Example~\ref{exatbbounded}.}
	\label{fig:exatbbounded}
\end{figure}
There are only two decision rules $\phi$ and $\psi$ which only differ at state $x^1$ with $\phi(x^1)=a^1$ and $\psi(x^1)=a^2.$ Thus,
\begin{equation}\label{eqtan}
v_\alpha^{\phi}(x^1)- v_\alpha^{\psi}(x^1) =\left(\sum_{t=0}^\infty\alpha^{2t}\right)-\left(\frac{26}{27}+\frac{2}{9}\alpha+\frac{14}{27}\alpha^2\sum_{t=0}^\infty\alpha^t\right)=\frac{(1-2\alpha)^3}{27(1-\alpha^2)}.
\end{equation}
Hence, $\alpha_*=0.5$ is the only irregular point. It is a simple break point, and $D([0,0.5))=\{\phi\},$ $D((0.5,1))=\{\psi\}.$ Note that $N(0.5)=1$ since $D(0.5)=\F=\{\phi,\psi\},$ and
\[V_{0.5}-V_{0,0.5}=V_{0.5}-s=
{\left[
	\begin{matrix}
		\frac{4}{3}&\frac{2}{3}&\frac{4}{3}&\frac{20}{27}&\frac{28}{27}
	\end{matrix}
	\right]}^{\mathrm{T}}
-
{\left[
	\begin{matrix}
		1&\frac{1}{3}&1&\frac{11}{27}&\frac{19}{27}
	\end{matrix}
	\right]}^{\mathrm{T}}
=\frac{1}{3}
{\left[
	\begin{matrix}
		1&1&1&1&1
	\end{matrix}
	\right]}^{\mathrm{T}}\]
is an eigenvector of eigenvalue $1$ for both $P(\phi)$ and $P(\psi).$ Therefore, Condition \hyperref[conA]{$\mathrm{A}$} holds at $\alpha_*=0.5$ by Proposition~\ref{propoA}(b)  with $K=0.$ For $2n$-horizon problems, we have
\begin{equation*}
\begin{aligned}
	T_{\alpha}^{\phi}V_{2n-1,\alpha}(x^1)-T_{\alpha}^{\psi}V_{2n-1,\alpha}(x^1) & = \left(\sum_{t=0}^{n-1}\alpha^{2t}+\alpha^{2n}\right)-\left(\frac{26}{27}+\frac{2}{9}\alpha+\frac{14}{27}\alpha^2\sum_{t=0}^{2n-3}\alpha^t+\frac{19}{27}\alpha^{2n}\right)\\
	& = \frac{(1-2\alpha)^2+(4\alpha-5)\alpha^{2n}}{27(1-\alpha^2)}(1-2\alpha).
\end{aligned}
\end{equation*}
Let $f_n(\alpha):=(1-2\alpha)^2+(4\alpha-5)\alpha^{2n}.$ Since $f_n(0.5)=-\frac{3}{4^n}<0$, by continuity, we have that for each $n\in\N^+$ there exists $\delta_n>0$ such that $f_n(\alpha)<0$ if $\alpha\in (0.5-\delta_n,0.5+\delta_n)$. Hence, we have that $T_{\alpha}^{\phi}V_{2n-1,\alpha}(x^1)<T_{\alpha}^{\psi}V_{2n-1,\alpha}(x^1),$ if $\alpha\in (0.5-\delta_n,0.5),$ and $T_{\alpha}^{\phi}V_{2n-1,\alpha}(x^1)>T_{\alpha}^{\psi}V_{2n-1,\alpha}(x^1),$ if $\alpha\in (0.5,0.5+\delta_n)$, which implies $N^*((0.5-\delta_n,0.5))>2n$ and $N^*((0.5,0.5+\delta_n))>2n$. Therefore, $N(\alpha)$ is unbounded at both sides of the simple break point $\alpha_*=0.5,$ even though Condition~\hyperref[conA]{$\mathrm{A}$} holds at this point. The graph of $N(\alpha)$ is provided by Figure~\ref{fig:Niu}.
\begin{figure}[ht!]
	\centerline{\includegraphics[scale=0.8]{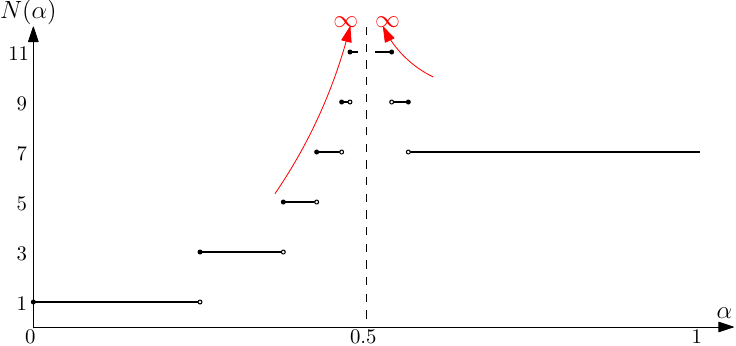}}
	\label{fig:Niu}
\end{figure}
In fact, both Conditions \hyperref[thm:BB-]{$\mathrm{B_-}$} and \hyperref[thm:BB+]{$\mathrm{B_+}$} are violated at $\alpha_*=0.5$ (these conditions are equivalent at $\alpha_*$ since it is a simple break point). Therefore, the assumptions of Theorem~\ref{tbnecessary} do not guarantee the boundedness of $N(\alpha)$ near $\alpha=0.5.$ Violations of both Conditions \hyperref[thm:BB-]{$\mathrm{B_-}$} and \hyperref[thm:BB+]{$\mathrm{B_+}$} can be verified by Proposition~\ref{propoB}: differentiating \eqref{eqtan} yields
\[\frac{d}{d\alpha}\left[v_\alpha^{\phi}(x^1)-v_\alpha^{\psi}(x^1)\right]\bigg|_{\alpha_*=0.5}=\frac{d}{d\alpha}\left[\frac{(1-2\alpha)^3}{27(1-\alpha^2)}\right]\bigg|_{\alpha_*=0.5}=0.\]
If $x\in\X\setminus\{x^1\},$ then $v^\phi(x)=v^\psi(x),$ and therefore $\frac{d}{d\alpha}\left[v_\alpha^{\phi}(x)- v_\alpha^{\psi}(x)\right]\Big|_{\alpha_*=0.5}=0.$ }
\end{exa}

        \begin{rem}\label{rem0v}\emph{As shown in Figure~\ref{fig:re}, the MDP from Example~\ref{exatbbounded} can be transformed into a similar MDP  with all terminal rewards  equal to $0.$}
\begin{figure}[ht!]
	\centerline{\includegraphics[scale=1]{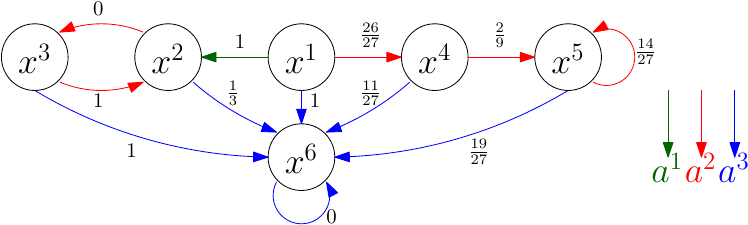}}
	\caption{MDP model for Remark~\ref{rem0v}.}
	\label{fig:re}
	\end{figure}
\end{rem}

\begin{exa}\label{exatbb} \emph{\textbf{\boldmath{This example demonstrates how Theorem~\ref{tbsufficient} can be used to prove the boundedness of $N(\alpha)$ near a break point.}} Consider the following MDP model; see Figure~\ref{fig:exatbb}.
\begin{figure}[ht!]
	\centerline{\includegraphics[scale=1.1]{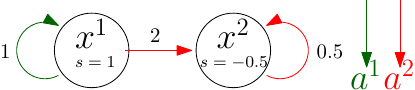}}
	\caption{MDP model for Example~\ref{exatbb}.}
	\label{fig:exatbb}
\end{figure}
There are only two decision rules $\phi$ and $\psi,$ which only differ at state $x^1,$ where $\phi(x^1)=a^1$ and $\psi(x^1)=a^2.$ Thus,
\begin{equation*}
	v_\alpha^{\phi}(x^1)=\frac{1}{1-\alpha},\qquad v_\alpha^{\psi}(x^1)=2+\frac{0.5\alpha}{1-\alpha}=\frac{2-1.5\alpha}{1-\alpha},\qquad v_\alpha^{\phi}(x^1)-v_\alpha^{\psi}(x^1)=\frac{1.5\alpha-1}{1-\alpha},
\end{equation*}
which implies that $\alpha_*=\frac{2}{3}$ is the unique irregular point, and it is a simple break point with $D([0,\frac{2}{3}))=\{\psi\}$ and $D((\frac{2}{3},1))=\{\phi\}.$ Note that $N(\frac{2}{3})=1$ since $D(\frac{2}{3})=\F=\{\phi,\psi\},$ and $V_{\frac{2}{3}}-s=
	{\left[
		\begin{matrix}
			3&1.5
		\end{matrix}
		\right]}^{\mathrm{T}}
	-
	{\left[
		\begin{matrix}
			1&-0.5
		\end{matrix}
		\right]}^{\mathrm{T}}
	=
	{\left[
		\begin{matrix}
			2&2
		\end{matrix}
		\right]}^{\mathrm{T}}$
is an eigenvector of eigenvalue $1$ for both $P(\phi)$ and $P(\psi),$ and therefore, by Proposition~\ref{propoA}(b), Condition \hyperref[conA]{$\mathrm{A}$} holds at $\alpha_*=\frac{2}{3}.$  Since it is a simple break point, Conditions \hyperref[thm:BB-]{$\mathrm{B_-}$} and \hyperref[thm:BB+]{$\mathrm{B_+}$} are equivalent. We verify Condition \hyperref[thm:BB+]{$\mathrm{B_+}$}. Note that state $x^2$ is always absorbing. For  $\pi\in\Pi,$ let $k_\pi\in\N^+\cup\{\infty\}$ be the first time when the policy $\phi^\pi$ uses action $a^2$ at state $x^1$ if $x_0=x^1.$ Then
\begin{equation*}
\begin{aligned}
	&\quad\frac{d}{d\alpha}\left(v_\alpha^{\phi^\pi}-v_\alpha^{\psi^\pi}\right)(x^1) = \frac{d}{d\alpha}\left[\sum_{t=0}^{k_\pi-1}\left(\alpha^t+\frac{2-1.5\alpha}{1-\alpha}\alpha^{k_\pi}\right)-\left(\frac{2-1.5\alpha}{1-\alpha}\right)\right]=1.5\left(\frac{1-\alpha^{k_\pi}}{1-\alpha}\right).\\
	& = \frac{d}{d\alpha}\left[\frac{1-\alpha^{k_\pi}}{1-\alpha}(1.5\alpha-1)\right]\bigg|_{\alpha=\frac{2}{3}}=1.5\left(\frac{1-\alpha^{k_\pi}}{1-\alpha}\right)\bigg|_{\alpha=\frac{2}{3}}\\
	& = 4.5\left[1-{\left(\frac{2}{3}\right)}^{k_\pi}\right]\quad\implies \inf_{\pi\in\Pi(D(\frac{2}{3}))}\frac{d}{d\alpha}\left(v_\alpha^{\phi^\pi}-v_\alpha^{\psi^\pi}\right)(x^1)\Big|_{\alpha=\frac{2}{3}}=1.5>0.
\end{aligned}
\end{equation*} 
Since Conditions \hyperref[conA]{$\mathrm{A}$}, \hyperref[thm:BB-]{$\mathrm{B_-}$} and \hyperref[thm:BB+]{$\mathrm{B_+}$} hold, by Theorem~\ref{tbsufficient} $N(\alpha)$ is bounded on both sides of $\alpha_*=\frac{2}{3}.$}  
\end{exa}

%

\begin{exa}\label{exanbm} \emph{\textbf{\boldmath{In this example, $N(\alpha)=m$ for all $\alpha\in(0,1)$, where $m$ is the number of states. This shows that the upper bound of the turnpike function for small discount factors in Theorem~\ref{ndboundm}(a) is sharp.}} Consider the  MDP presented in Figure~\ref{fig:exanbm}. The broken dots in the figure mean $p(x^{i+1}|x^i,a^2)=1$ and $r(x^i,a^2)=1$ for $i=1,2,\ldots,m-1.$
\begin{figure}[ht!]
	\centerline{\includegraphics[scale=1.2]{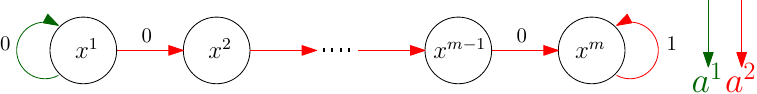}}
	\caption{MDP model for Example~\ref{exanbm}.}
	\label{fig:exanbm}
\end{figure}
There are only two decision rules $\phi,\psi$ which only differ at $x^1$ with $\phi(x^1)=a^1$ and $\psi(x^1)=a^2.$ It is easy to see that $v_\alpha^{\phi}(x^1)=0,$ $v_\alpha^{\psi}(x^1)=\frac{\alpha^{m-1}}{1-\alpha}$, and thus $D(\alpha)=\{\psi\}$ for all $\alpha\in(0,1)$. We have
\[T_\alpha^{\phi}V_{n-1,\alpha}(x^1)=\begin{cases}
	0, &\text{if $n\le m,$} \\
    \frac{\alpha^n}{1-\alpha}, &\text{if $n\ge m+1,$}
\end{cases}\qquad T_\alpha^{\psi}V_{n-1,\alpha}(x^1)=\begin{cases}
	0, &\text{if $n\le m-1,$} \\
	\frac{\alpha^{n-1}}{1-\alpha}, &\text{if $n\ge m.$}
\end{cases}\]
This implies that for all $\alpha\in(0,1)$ we have $D_n(\alpha)=\{\phi,\psi\}$ for $n\le m-1$ and $D_n(\alpha)=\{\psi\}$ for $n\ge m.$ Therefore, $N(\alpha)=m$ for all $\alpha\in(0,1)$.}
\end{exa}

\begin{exa}\label{exaa1bound} \emph{\textbf{\boldmath{In this example $a_1=\Delta_L=\tilde{\Delta}_L=0.5,$ and the inequality stated in Corollary~\ref{a1bound} holds in the form of the equality. Also, $\Delta_L$ can be the largest value of $\delta>0$ such that $N^*((0,\delta))<\infty,$ while $N^*((0,\Delta_L))<\infty$ in view of Theorem~\ref{ndboundm}(a).}} Consider the following MDP; see Figure~\ref{fig:exaa1bound}.
\begin{figure}[ht!]
	\centerline{\includegraphics[scale=1.2]{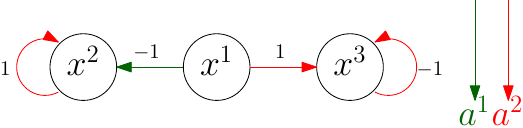}}
	\caption{MDP model for Example~\ref{exaa1bound}.}
	\label{fig:exaa1bound}
\end{figure}
There are only two decision rules $\phi$ and $\psi$ which only differ at state $x^1$ with $\phi(x^1)=a^1$ and $\psi(x^1)=a^2.$ We have
\[v_\alpha^{\phi}(x^1)- v_\alpha^{\psi}(x^1)=\left(-1+\sum_{t=1}^\infty\alpha^t\right)-\left(1-\sum_{t=1}^\infty\alpha^t\right)=\frac{2(2\alpha-1)}{1-\alpha}.\]
Hence, $a_1=0.5$ is the only irregular point,  $D([0,0.5))=\{\psi\},$ and $D((0.5,1))=\{\phi\}.$ Note that the rewards in this MDP are already balanced with $R^*=R_1^*=1.$ By formulae \eqref{FXdef}$-$\eqref{Cdef}, we  have $L=0$ and $C_0=2.$   By formula \eqref{Deldef}, $\Delta_0=\tilde{\Delta}_0=0.5.$ Therefore, the equality of Corollary~\ref{a1bound} is attained. Furthermore,
\[T_\alpha^{\phi}V_{n-1,\alpha}(x^1)-T_\alpha^{\psi}V_{n-1,\alpha}(x^1)=\left(-1+\sum_{t=1}^{n-1}\alpha^t\right)-\left(1-\sum_{t=1}^{n-1}\alpha^t\right)=\frac{2(2\alpha-1-\alpha^n)}{1-\alpha}.\]
Let $f_n(\alpha)=2\alpha-1-\alpha^n.$ For $n\ge3$ we have $f_n'(\alpha)=2-n\alpha^{n-1}$ with $f_n'(0)=2>0$ and $f_n'(1)=2-n<0,$ and $f_n''(\alpha)=-n(n-1)\alpha^{n-2}<0.$ This implies that for each $n\ge3$ the function $f_n(\alpha)$ first increases and then decreases in $[0,1).$ Since $f_n(0.5)<0$ and $f_n(1)=0,$ we conclude that for each $n\ge3$ the function $f_n(\alpha)$ has a unique zero in $(0.5,1),$  which we denote by $\alpha(n).$ We set $\alpha(1)=\alpha(2)=1.$ Then $T_\alpha^{\phi}V_{n-1,\alpha}(x^1)<T_\alpha^{\psi}V_{n-1,\alpha}(x^1)$ for all $\alpha\in(0,\alpha(n)),$ which implies $N(\alpha)=1$ for all $\alpha\in[0,0.5],$ but $N(\alpha)>n$ for all $\alpha\in(0.5,\alpha(n)],n\ge3.$ Therefore, $N^*((0,0.5))=1,$ while $N^*((0,\alpha))=\infty$ for $\alpha\in(0.5,1).$}
\end{exa}

\section{Conclusion}\label{seccon}
This paper continues the study of turnpike functions for discounted MDPs initiated in 
\cite{doi:10.1287/moor.2017.0912}. Turnpike functions are important for applications because they provide the number of iterations by VIA guaranteeing computations of optimal policies. This number can change drastically with small changes of discount factors near irregular points, and Algorithm~\ref{algequi} detects in strongly polynomial time whether a discount factor is irregular. For small discount factors with an upper bound $\Delta_L,$ Theorem~\ref{ndboundm}(a) guarantees that the turnpike function is bounded by the number of states. Based on this result, Algorithm~\ref{algsmall} solves an MDP in strongly polynomial time for all small discount factors with an upper bound $\tilde{\Delta}_L\ge\Delta_L.$ Definitions of $\Delta_L$ and $\tilde{\Delta}_L$ are provided in \eqref{Deldef}.

The results of this paper are based on the assumption that the sets of states and actions are finite. For MDPs with infinite action or state sets, optimal policies may not exist, and, if they exist, finding an optimal policy by the VIA can require an infinite number of iterations, and therefore a turnpike integer can be infinite. While we hope that some relevant results can be established for more general MDPs, this paper is restricted to MDPs with finite numbers of states and actions.


        While Algorithm~\ref{algsmall} computes a deterministic policy optimal for all discount factors in some neighborhood of 0,   it would be natural to design an algorithm for finding an optimal deterministic policy for discount factors in some left and right neighborhoods of $\alpha\in (0,1).$ If $\alpha$ is a regular point, then Proposition~\ref{partition}(c) shows that the set of optimal policies for $\alpha$ is also the set of optimal policies for all discount factors from the same partition interval, and the question becomes how to find the endpoints of this partition interval. 
If $\alpha$ is an irregular point, Proposition~\ref{upperhemi} guarantees that any optimal policy for discount factors close to $\alpha$ is also optimal for $\alpha,$ and the question becomes how to detect which optimal policies for $\alpha$ are optimal for which side of $\alpha.$ If $\alpha$ is the largest break point, then a Blackwell-optimal policy \cite{BW,BWD} is optimal for all discount factors from the interval $[\alpha,1).$  So, studying properties of break points and developing algorithms for finding them may lead to the development of methods for solving parametric dynamic programming problems for various values of discount factors.

\nomenclature{\(\N\)}{set of natural numbers including 0}
\nomenclature{\(\N^+\)}{set of positive integers}
\nomenclature{\(\R\)}{set of real numbers}
\nomenclature{\(\X\)}{set of states}
\nomenclature{\(q_x\)}{number of available actions at state $x$}
\nomenclature{\(A(x)\)}{set of available actions at state $x$}
\nomenclature{\(\A\)}{set of all actions}
\nomenclature{\(m\)}{number of states}
\nomenclature{\(q\)}{number of state-action pairs}
\nomenclature{\(s(x)\)}{terminal reward at state $x$}
\nomenclature{\(r(x,a)\)}{one-step reward when action $a$ is applied at state $x$}
\nomenclature{\(p(y\lvert x,a)\)}{transition probability from state $x$ to state $y$ when action $a$ is applied}
\nomenclature{\(\phi\)}{decision rule (deterministic policy)}
\nomenclature{\(\pi\)}{Markov policy}
\nomenclature{\(\F\)}{set of decision rules (deterministic policies)}
\nomenclature{\(\Pi\)}{set of Markov policies}
\nomenclature{\(2^X\)}{set of all subsets of a set $X$}
\nomenclature{\(\Pi(S)\)}{set of Markov policies whose decision rules at each step are from $S$}
\nomenclature{\(\phi^\pi\)}{policy that uses the decision rule $\phi$ at the initial step and then follows the policy $\pi$}
\nomenclature{\(P(\phi)\)}{stochastic matrix of transition probabilities defined by a deterministic policy $\phi$}
\nomenclature{\(P_t(\pi)\)}{stochastic matrix of transitions between steps 0 and $t$ for a  Markov policy $\pi$}
\nomenclature{\(v_{n,\alpha}^{\pi}(x)\)}{expected $n$-horizon $\alpha$-discounted total reward starting from state $x$ under policy $\pi$}
\nomenclature{\(v_\alpha^{\pi}(x)\)}{expected infinite horizon $\alpha$-discounted total reward starting from state $x$ under policy $\pi$}
\nomenclature{\(V_{n,\alpha}(x)\)}{$n$-horizon value function at state $x$}
\nomenclature{\(V_\alpha(x)\)}{infinite-horizon value function at state $x$}
\nomenclature{\(T_\alpha^\phi\)}{total reward operator defined in \eqref{eqTalpha}}
\nomenclature{\(T_\alpha\)}{optimality operator defined in \eqref{eqTalpha}}
\nomenclature{\(\mathbf{0}\)}{column vectors in $\R^m$ with all entries equal to $0$}
\nomenclature{\(\mathbf{1}\)}{column vectors in $\R^m$ with all entries equal to $1$}
\nomenclature{\(R_1\)}{defined in \eqref{eqR}}
\nomenclature{\(R_2\)}{defined in \eqref{eqR}}
\nomenclature{\(R\)}{defined in \eqref{eqR}}
\nomenclature{\(R_1^*\)}{defined in \eqref{eqRb}}
\nomenclature{\(R_2^*\)}{defined in \eqref{eqRb}}
\nomenclature{\(R^*\)}{defined in \eqref{eqRb}}
\nomenclature{\(\norm{\cdot}\)}{infinity norm of matrices}
\nomenclature{\(G(\phi,\mathbf{v})\)}{defined in \eqref{eqG}}
\nomenclature{\(M_n(\alpha)\)}{set of Markov optimal policies for an $n$-horizon problem with discount factor $\alpha$}
\nomenclature{\(D_n(\alpha)\)}{set of first-step optimal policies for an $n$-horizon problem with discount factor $\alpha$}
\nomenclature{\(D(\alpha)\)}{set of optimal policies for an infinite-horizon problem with discount factor $\alpha$}
\nomenclature{\(M_n(E)\)}{set of Markov policies that are optimal for some discount factor $\alpha\in E$}
\nomenclature{\(D_n(E)\)}{set of deterministic policies that are first-step-optimal for some discount factor $\alpha\in E$}
\nomenclature{\(D(E)\)}{set of deterministic policies that are optimal for some discount factor $\alpha\in E$}
\nomenclature{\(N(\alpha)\)}{turnpike integer at discount factor $\alpha$ (turnpike function)}
\nomenclature{\(N^*(E)\)}{supremum of turnpike integers in $E$}
\nomenclature{\(M_n(\alpha-)\)}{set of left optimal Markov policies for an $n$-horizon problem with discount factor $\alpha$}
\nomenclature{\(M_n(\alpha+)\)}{set of right optimal Markov policies for an $n$-horizon problem with discount factor $\alpha$}
\nomenclature{\(D_n(\alpha-)\)}{set of first-step-left optimal decision rules for an  $n$-horizon problem with discount factor $\alpha$}
\nomenclature{\(D_n(\alpha+)\)}{set of first-step-right optimal decision rules for an $n$-horizon problem with discount factor $\alpha$}
\nomenclature{\(D(\alpha-)\)}{set of left optimal decision rules for an infinite-horizon problem with discount factor $\alpha$}
\nomenclature{\(D(\alpha+)\)}{set of left optimal decision rules for an infinite-horizon problem with discount factor $\alpha$}
\nomenclature{\(a_i\)}{the $i$-th nonzero irregular point}
\nomenclature{\(b_i^{(n)}\)}{the $i$-th nonzero   $n$-horizon-first-step irregular point}
\nomenclature{\(a_i^{(n)}\)}{the $i$-th nonzero $n$-horizon irregular point}
\nomenclature{\(\I\)}{set of nonempty subintervals of $[0,1)$ including singletons}
\nomenclature{\(\mathcal{I}\)}{subinterval of $[0,1)$ (can be a singleton)}
\nomenclature{\(\mathbb{D}^-(\mathcal{I})\)}{set of point in $\mathcal{I}$ where the turnpike function is not left-continuous}
\nomenclature{\(\mathbb{D}^+(\mathcal{I})\)}{set of point in $\mathcal{I}$ where the turnpike function is not right-continuous}
\nomenclature{\(\mathbb{D}(\mathcal{I})\)}{set of points in $\mathcal{I}$ at which the turnpike function is discontinuous}
\nomenclature{\(\Hat{\mathbb{D}}(\mathcal{I})\)}{set of points in $\mathcal{I}$ at which the turnpike function is neither left-continuous nor right-continuous}
\nomenclature{\(\P(\mathcal{I})\)}{set of irregulars point in $\mathcal{I}$}
\nomenclature{\(\P_n(\mathcal{I})\)}{set of $n$-horizon-first-step irregular point}
\nomenclature{\(\P_n^T(\mathcal{I})\)}{set of $n$-horizon-first-step touching point}
\nomenclature{\(\P_n^B(\mathcal{I})\)}{set of $n$-horizon-first-step break point}
\nomenclature{\(\P\)}{set of irregular points}
\nomenclature{\(\mu\)}{Lebesgue measure on $[0,1]$}
\nomenclature{\(\F_n\)}{defined in \eqref{FXdef} and \eqref{FXalt}}
\nomenclature{\(\X_n\)}{defined in \eqref{FXdef} and \eqref{FXalt}}
\nomenclature{\(L\)}{defined in \eqref{Ldef}}
\nomenclature{\(L_n\)}{defined in \eqref{Ldef}}
\nomenclature{\(H\)}{the smallest natural number such that $L_H=L$}
\nomenclature{\(C_n\)}{defined in \eqref{Cdef} or \eqref{defaltern2}}
\nomenclature{\(\Delta_n\)}{defined in \eqref{Deldef}}
\nomenclature{\(\tilde{\Delta}_n\)}{defined in \eqref{Deldef}}
\nomenclature{\(r_n(x)\)}{defined in \eqref{defaltern1}}
\nomenclature{\(A_n(x)\)}{defined in \eqref{defaltern1}}
\printnomenclature

\end{document}